\documentclass[10pt,a4paper]{article}
\usepackage{bbm}

\usepackage{amsmath}  
\usepackage{amsfonts}
\usepackage{graphicx}
\usepackage{amsmath}
\usepackage{amssymb}
\usepackage{latexsym}
\usepackage{amsmath, amsfonts,amssymb, amsthm, euscript,makeidx,color,mathrsfs}
\usepackage{enumerate}
\usepackage{mathtools}

\usepackage{lipsum}
\usepackage{amsmath}
\usepackage{amssymb}
\usepackage{dsfont}
\usepackage{mathrsfs}
\usepackage{graphicx}
\usepackage{subfigure}
\usepackage{color}
\usepackage{epstopdf}

\usepackage{algorithmic}
\usepackage{amsopn}

\usepackage[colorlinks,linkcolor=cyan,anchorcolor=green,citecolor=red]{hyperref}

\def\sqr#1#2{{\vcenter{\vbox{\hrule height.#2pt
              \hbox{\vrule width.#2pt height#1pt \kern#1pt \vrule width.#2pt}
              \hrule height.#2pt}}}}
\def\5n{\negthinspace \negthinspace \negthinspace \negthinspace \negthinspace }
\def\4n{\negthinspace \negthinspace \negthinspace \negthinspace }
\def\3n{\negthinspace \negthinspace \negthinspace }
\def\2n{\negthinspace \negthinspace }
\def\1n{\negthinspace }

\newcommand{\lint}{\mathalpha{\ltimes}}
\newcommand{\rint}{\mathalpha{\rtimes}}

\def\dbE{\mathbb{E}}
\def\dbF{\mathbb{F}}

\def\dbL{\mathbb{L}}

\def\dbP{\mathbb{P}}

\def\dbR{\mathbb{R}}
\def\dbS{\mathbb{S}}

\def\sL{\mathscr{L}}

\def\={\buildrel \triangle \over =}

\def\ds{\displaystyle}

\def\ns{\noalign{\ss}}
\def\a{\alpha}
\def\b{\beta}

\def\d{\delta}

\def\z{\zeta}

\def\l{\lambda}
\def\m{\mu}

\def\si{\sigma}
\def\t{\tau}
\def\f{\varphi}
\def\th{\theta}
\def\o{\omega}
\def\p{\phi}
\def\i{\infty}
\def\G{\Gamma}
\def\D{\Delta}
\def\Th{\Theta}

\def\O{\Omega}
\def\cA{{\cal A}}
\def\cB{{\cal B}}
\def\cC{{\cal C}}

\def\cE{{\cal E}}
\def\cF{{\cal F}}
\def\cG{{\cal G}}

\def\cM{{\cal M}}

\def\cP{{\cal P}}

\def\cR{{\cal R}}
\def\cS{{\cal S}}
\def\cT{{\cal T}}

\def\BA{{\bf A}}

\def\BC{{\bf C}}

\def\BX{{\bf X}}

\def\BBf{\boldsymbol\varphi}

\def\no{\noindent}

\def\ss{\smallskip}
\def\ms{\medskip}

\def\q{\quad}
\def\qq{\qquad}

\def\blan{\big\langle}
\def\bran{\big\rangle}

\def\rf{\eqref}

\def\h{\widehat}
\def\wt{\widetilde}
\def\ti{\tilde}
\def\cd{\cdot}
\def\cds{\cdots}

\def\ae{\hbox{\rm a.e.}}

\def\les{\leqslant}
\def\ges{\geqslant}

\def\({\Big (}
\def\){\Big )}
\def\[{\Big[}
\def\]{\Big]}

\def\bde{\begin{definition}\label}
\def\ede{\end{definition}}
\def\be{\begin{equation}}
\def\bel{\begin{equation}\label}
\def\ee{\end{equation}}
\def\bt{\begin{theorem}\label}
\def\et{\end{theorem}}
\def\bc{\begin{corollary}\label}
\def\ec{\end{corollary}}
\def\bl{\begin{lemma}\label}
\def\el{\end{lemma}}
\def\bp{\begin{proposition}\label}
\def\ep{\end{proposition}}
\def\bas{\begin{assumption}\label}
\def\eas{\end{assumption}}
\def\br{\begin{remark}\label}
\def\er{\end{remark}}
\def\bex{\begin{example}\label}
\def\ex{\end{example}}
\def\ba{\begin{array}}
\def\ea{\end{array}}
\def\bpf{\begin{proof}}
\def\epf{\end{proof}}
\def\ed{\end{document}}

\def\square#1{\vbox{\hrule\hbox{\vrule height#1%
     \kern#1\vrule}\hrule}}
\def\rectangle#1#2{\vbox{\hrule\hbox{\vrule height#1%
     \kern#2\vrule}\hrule}}

\font\tenbb=msbm10 \font\sevenbb=msbm7 \font\fivebb=msbm5

\newfam\bbfam
\scriptscriptfont\bbfam=\fivebb \textfont\bbfam=\tenbb
\scriptfont\bbfam=\sevenbb

\newtheorem{theorem}{\hskip 1.3em Theorem}[section]
\newtheorem{definition}[theorem]{\hskip 1.3em Definition}
\newtheorem{proposition}[theorem]{\hskip 1.3em Proposition}
\newtheorem{corollary}[theorem]{\hskip 1.3em Corollary}
\newtheorem{lemma}[theorem]{\hskip 1.3em Lemma}
\newtheorem{remark}[theorem]{\hskip 1.3em Remark}
\newtheorem{example}[theorem]{\hskip 1.3em Example}

\newtheorem{assumption}[theorem]{\hskip 1.3em Assumption}

\makeatletter
   \renewcommand{\theequation}{%
            \thesection.\arabic{equation}}
   \@addtoreset{equation}{section}
\makeatother

\begin{document}

\title{\bf Closed-loop solvability of delayed control problems:\qq A stochastic Volterra system approach
\thanks{This work is supported by the China National Key Research and Development Program (2024YFA1012800),
the National Natural Science Foundation of China (62433020, 12371449, 12501618, 12471419).
}
}
\author{\normalsize Weijun Meng\thanks{\it School of Mathematics $\&$ Statistics, Nanjing University of Science and Technology, Nanjing 210094, China, E-mail: mengwj@njust.edu.cn}
\quad Tianxiao Wang\thanks{\it School of Mathematics, Sichuan University, Chengdu 610065, China, E-mail: wtxiao2014@scu.edu.cn}
 \thanks{\it Corresponding author.}
\quad Ji-Feng Zhang\thanks{\it School of Automation and Electrical Engineering, Zhongyuan University of Technology, Zhengzhou 450007; State Key Laboratory of Mathematical Sciences, Academy of Mathematics and Systems Science, Chinese Academy of Sciences, Beijing 100190, and school of Mathematical Sciences, University of Chinese Academy of Sciences, Beijing 100149, China, E-mail: jif@iss.ac.cn}
}

\ms

\maketitle

\noindent{\bf Abstract:}\quad A general and new stochastic linear quadratic optimal control problem is studied, where the coefficients are allowed to be time-varying, and both state delay and control delay can appear simultaneously in the state equation and the cost functional. The closed-loop outcome control of this delayed problem is given by a new Riccati system whose solvability is carefully established. To this end, a novel method is introduced to transform the delayed problem into a control problem driven by a stochastic Volterra integral system without delay. This method offers several advantages: it bypasses the difficulty of decoupling the forward delayed state equation and the backward anticipated adjoint equation, avoids the introduction of infinite-dimensional spaces and unbounded control operators, and ensures that the closed-loop outcome control depends only on past state and control, without relying on future state or complex conditional expectation calculations. Finally, several particular important stochastic systems are discussed. It is found that the model can cover a class of stochastic integro-differential systems, whose closed-loop solvability has not been available before.

\vspace{2mm}

\ms


\vspace{2mm}

\ms


\rm

\section{Introduction}

\subsection{Delayed optimal control problems}

In this paper, given $0\les t_0<T$ and the constant delay time $\d > 0$, we consider an optimal control problem where the state equation is described as follows:
\bel{state}\left\{\ba{ll}
\ns\ds  dx(t)=\big[A_1(t)x(t)+A_2(t) y(t)+A_3(t)z(t)+B_1(t)u(t)+B_2(t)\nu(t)+B_3(t)\m(t)+b(t)\big]dt\\
\ns\ds\qq\qq+\big[C_1(t)x(t)+C_2(t) y(t)+C_3(t)z(t)+D_1(t)u(t)+\si(t)\big]dW(t),\q t\in(t_0,T),\\
\ns\ds x(t)=\xi(t-t_0),\ \ \ u(t)=\varsigma(t-t_0),\ \ t\in[t_0-\d,t_0],
\ea\right.
\ee
and the cost functional is defined as:
\begin{eqnarray}
J(t_0,\xi(\cd),\varsigma(\cd);u(\cd))\1n=\mathbb{E}\int_{t_0}^T\2n\bigg[
x(t)^\top Q_1(t)x(t)+y(t)^\top Q_2(t)y(t)+z(t)^\top Q_3(t)z(t)
\notag\\
+u(t)^\top R_1(t)u(t)+\nu(t)^\top R_2(t)\nu(t)
\bigg]dt.\qq\qq\q\label{cost}
\end{eqnarray}
Here $(\Omega,\mathcal{F},\dbF,\dbP)$ is a complete filtered probability space, $\dbE[\cd]$ denotes the expectation, and the filtration $\dbF=\{\mathcal{F}_t\}_{t\ges 0}$ is generated by a one-dimensional standard Brownian motion $\{W(t)\}_{t\ges0}$.
In the above, $x(\cd)\in\dbR^n$ is state, $u(\cd)\in \dbR^m$ is control, $y(\cd)$ and $\nu(\cd)$ are the pointwise delays of the state and the control, respectively, $z(\cd)$ and $\m(\cd)$ are the corresponding (extended) distributed delays. They are defined as follows:
\bel{yz-nu-mu}\3n\ba{ll}
\ns\ds
y(t)\equiv x(t-\d),\ \ z(t)\equiv\int_{t_0}^t F(t,s)x(s)ds,\ \nu(t)\equiv u(t-\d),\ \m(t)\equiv\int_{t_0}^t\ti F(t,s)u(s)ds,\q t\in(t_0,T).
\ea\ee
In addition, $\xi(\cd) $ and $\varsigma(\cd)$ are called the initial trajectories of the state and the control, respectively. The conditions satisfied by coefficients $F,\ti F,A_i,B_i,C_i,D_1,Q_i,R_1,R_2,b$, $\si$ above will be specified in Section 2, which ensure the well-posedness of the cost functional \rf{cost}, $i=1,2,3$. The optimal control problem is stated as follows:

\ms

\textbf{Problem (P).} To find a control $u^*(\cdot)$ such that (\ref{state}) is satisfied and (\ref{cost}) is minimized, i.e.,
\vskip-2mm
$$\ba{ll}
J(t_0,\xi(\cd),\varsigma(\cd);u^*(\cd))  =\inf\limits_{u(\cd)\in\,L_{\dbF}^2(t_0,T;\dbR^m)}J(t_0,\xi(\cd),\varsigma(\cd);u(\cd)) =V(t_0,\xi(\cd),\varsigma(\cd)),
\ea$$
\vskip-1mm
\no where $L_{\dbF}^2(t_0,T;\dbR^m)$ is the Hilbert space consisting of $\mathbb{F}$-adapted processes $\phi(\cd)$ such that $\mathbb{E}\int_0^T|\phi(t)|^2dt<\infty$. Any $u^*(\cdot)\in L_{\dbF}^2(t_0,T;\dbR^m)$ that achieves the above infimum is called an {\it open-loop optimal control} for $(t_0,\xi(\cd),\varsigma(\cd))$ and the corresponding solution $x^*(\cdot)$ is called the {\it open-loop optimal trajectory}. $(x^*(\cdot),u^*(\cdot))$ is called an {\it open-loop optimal pair}. $V(t_0,\xi(\cd),\varsigma(\cd))$ is called the value function.
In the special case, when $b(\cd)$ and $\si(\cd)$ vanish, we denote the corresponding delayed linear quadratic (LQ) problem, the cost functional, and the value function by Problem (P$_0$), $J_0(t_0,\xi(\cd),\varsigma(\cd);u(\cd))$ and $V_0(t_0,\xi(\cd),\varsigma(\cd))$, respectively.

\ss

As to the above \rf{state}, it is an extension of the following stochastic differential delay equation (SDDE):
\bel{state-distributed}\left\{\ba{ll}
\ns\ds  dx(t)=\big[A_1(t)x(t)+A_2(t) y(t)+A_3(t)\wt z(t)+B_1(t)u(t)+B_2(t)\nu(t)+B_3(t)\wt \m(t)+\wt b(t)\big]dt\\
\ns\ds\qq\qq+\big[C_1(t)x(t)+C_2(t) y(t)+C_3(t)\wt z(t)+D_1(t)u(t)+\wt \si(t)\big]dW(t),\q t\in(t_0,T),\\
\ns\ds x(t)=\xi(t-t_0),\ \ \ u(t)=\varsigma(t-t_0),\ \ t\in[t_0-\d,t_0],\\
\ns\ds \wt z(t)\equiv\int_{t-\d}^t G_1(t,s)x(s)ds,\ \wt \m(t)\equiv\int_{t-\d}^t G_2(t,s)u(s)ds,\q t\in(t_0,T),
\ea\right.
\ee
%
%
where $y(\cd)$ and $\nu(\cd)$ are defined in \rf{yz-nu-mu}, and $G_1$, $G_2$ are bounded. In fact,
$$\ba{ll}
\ns\ds \wt z(t)=\wt G_1(t) +\int_{t_0}^{t} \h G_1(t,s) x(s)ds,\ \ \wt \mu(t)=\wt G_2(t) +\int_{t_0}^{t} \h G_2(t,s) u(s)ds,
\ea
$$
where
$$\ba{ll}
\ns\ds \wt G_1(t)\equiv \int_{t-\d}^{t_0} G_1(t,s)\xi(s-t_0)ds {\bf 1}_{[t_0,t_0+\d)}(t),\\
\ns\ds \wt G_2(t)\equiv  \int_{t-\d}^{t_0} G_2(t,s)\varsigma(s-t_0)ds {\bf 1}_{[t_0,t_0+\d)}(t),\\
\ns\ds \h G_1(t,s)\equiv \big[{\bf 1}_{[t_0,t_0+\d)}(t)+{\bf 1}_{[t_0+\d,T]}(t){\bf 1}_{[t-\d,t)}(s) \big]G_1(t,s),\\
\ns\ds \h G_2(t,s)\equiv  \big[{\bf 1}_{[t_0,t_0+\d)}(t)+{\bf 1}_{[t_0+\d,T]}(t){\bf 1}_{[t-\d,t)}(s) \big]G_2(t,s).
\ea
$$
Therefore, \rf{state-distributed} can be seen as a special case of \rf{state} with
$$\ba{ll}
\ns\ds b=\wt b+A_3\wt G_1+B_3\wt G_2,\ \ \si=\wt \si+C_3\wt G_1,  \ \ F=\h G_1,\ \
\ti F=\h G_2.
\ea
$$
Based on this fact, in the following we name \rf{state} \it the (extended) controlled SDDE. \rm More details about SDDE \rf{state-distributed} can be referred to \cite{Mohammed-1984,Mohammed-1998}.  Notice that here we allow $t_0+\d\les T$.

\subsection{Motivations}

In the real world, many challenges across disciplines like economics, finance, aerospace, and network communication can be framed as optimal control problems \cite{Carmona-Fouque-Mousavi-Sun-2018, Chen-Wu-2020, Federico-2011}. Moreover, the evolution of certain phenomena hinges not only on present conditions but also on their historical trajectories. Consequently, optimal control problem of stochastic control systems containing state delay and control delay, like the above \rf{state-distributed},
is an important issue in control theory.
The relevant optimal control problems have attracted enormous attention of the optimization and engineering communities in the last
decades. We refer to the monographes by Bensoussan--Da Prato--Delfour--Mitter \cite{Bensoussan-Prato-Delfour-Mitter-2007}, Chang \cite{Chang-book}, Hale \cite{Hale-1977-book}, Kolmanovskii--Myshkis \cite{Kolmanovskii-Myshkis}, Kolmanovskii--Shaikhet \cite{Kolmanovskii-Shaikhet}, Meng--Shi--Yong \cite{Meng-Shi-Yong-2023}, Zhang--Xie \cite{Zhang-Xie-2007-book}.

\ss

To motivate the current study, let us make some careful discussions about the existing papers from the viewpoints of frameworks, methodologies, and conclusions.

{\bf I.} The existing frameworks on stochastic delayed systems seem scattered and decentralized, many of which cannot fully cover each other.

$\bullet$  Time-invariant versus time-varying coefficients: Such a difference happens even for deterministic controlled systems. For example, \cite{Lee-You-1989, Lefebvre-Miller-2021, Liang-Xu-Zhang-2018, Ma-Xu-Zhang-2022,  Ma-Qi-Li-Zhang-2022, Meng-Shi-Zhang-Zhao-2025, You-1999, Zhang-Xu-2017} and the relevant papers therein are devoted to the time-invariant case, while \cite{Alekal-Brunovsky-Chyung-Lee-1971, Bagchi-1978,Basar-Moon-Yong-Zhang-2023-CDC, Chen-Wu-Yu-2012, Delfour-Mitter-1972, Eller-Aggarwal-Banks, Lindquist-1973, Zhang-Xiong-Shi-2021} and the relevant papers therein are concerned with the time-varying one. The appearance of both cases in the literature is attributed to the complexity of delayed problems and the methodology limitation developed accordingly.

$\bullet$  State delay versus control delay: We found that the inclusion of delayed terms indeed brings significant differences to the corresponding study. When only state variables contain the delayed terms, we refer to \cite{Delfour-1977, Delfour-Mitter-1972,Duan-2023,Eller-Aggarwal-Banks} for the deterministic case and \cite{DeFeo-Swiech,Flandoli-1990,Liang-Xu-Zhang-2018, Lindquist-1973} for the stochastic case. When delayed terms are only added to control variables, we refer to \cite{Ichikawa-1982,You-1983} for the deterministic case and \cite{Lefebvre-Miller-2021,Liang-Xu-Zhang-2018, Wang-Zhang-2013, Zhang-Xu-2017} for the stochastic case. In fact, the reason, why models with state delay or control delay are discussed separately in the literature, lies in that the methods for handling them are fundamentally different.

$\bullet$ Pointwise delay versus distributed delay: We found that delay types also affect the analogue study. For example, \cite{Chen-Wu-Yu-2012, Eller-Aggarwal-Banks,Lefebvre-Miller-2021, Liang-Xu-Zhang-2018, Zhang-Xu-2017} are concerned with the pointwise delay in state or control, while \cite{Chen-Wu-2020,DeFeo-Swiech} are devoted to the case with distributed delay. The poinwise delay is usually more challenging than the distributed delay because the distributed delay can be dealt with using the derivative formula or the It\^o formula, but the pointwise delay cannot.

$\bullet$ State (or control) dependence versus constants in diffusion terms: In the literature, we found that whether diffusion terms depend on the state/control variables, or even their delayed terms, brings essential differences. We refer to \cite{Fabbri-Federico-2014,DeFeo-Swiech, Lefebvre-Miller-2021,Zhang-Xu-2017} for the former case while \cite{Bagchi-1978,Carmona-Fouque-Mousavi-Sun-2018, Chen-Wu-2020} for the simple constant case. In our opinion, the related reason is that this dependence occurs within the stochastic integral terms, rendering the previously used methods for handling the Lebesgue integral terms ineffective.

$\bullet$ Delayed terms in cost functionals or not: Eventually, let us point out the different frameworks caused by the dependence of delayed terms in cost functionals. Even for deterministic controlled systems, pointwise delay in cost functionals makes weight operators unbounded, while distributed delay makes their processing highly dependent on the given weight coefficients (see the Introduction part of \cite{Lee-You-1989} for detailed explanations). When cost functionals contain state/control delay terms, we refer to \cite{Lee-You-1989,Meng-Shi-Wang-Zhang-2023, Meng-Shi-Yong-2023,Zhang-Xu-2017},
while most other existing studies do not touch this topic.

Based on the above classifications of the models, we pose the first question:

\it {\bf (Q1)}: Is it possible to provide the study in general framework to cover the above models? \rm

\ms

{\bf II.} After careful observations of the existing papers, we found that there are several approaches treating LQ problems for stochastic delay systems. We list them out and make some analyses as follows.

$\bullet$ Maximum principle: It is a very natural method to study optimal control for stochastic delayed systems. We refer to \cite{Chen-Zhang-2023,Xu-Shi-Zhang-2018,Zhang-Xu-2017} for more details. To obtain the explicit forms of optimal controls, one has to decouple the resulting forward-backward stochastic systems, which is the challenging part. A success example is
\cite{Zhang-Xu-2017}, where the system only contains pointwise control delay, and a class of Riccati equation is established by decoupling. The difficulty in decoupling lies in that the delay effect disrupts the classical decoupling relation and makes It\^o formula fail to work \cite{Basar-Moon-Yong-Zhang-2023-CDC}.

$\bullet$  Dynamic programming: This method is very effective in obtaining feedback forms of optimal controls, see \cite{Alekal-Brunovsky-Chyung-Lee-1971, Bensoussan-Prato-Delfour-Mitter-2007,Larssen-2002, Lefebvre-Miller-2021,Liang-Xu-Zhang-2018} for more details. However, the main difficulty for delayed problems is the lack of Markovianity, which necessitates addressing the inherently infinite-dimensional Hamilton-Jacobi-Bellman (HJB) equations at the cost of increased complexity.

$\bullet$  State space method: This method is to lift state equations to Banach or Hilbert spaces (depending on the regularity of the data). It is developed in the deterministic case by \cite{Delfour-Mitter-1972,Ichikawa-1982} and then extended in the stochastic case by \cite{ DeFeo-Federico-Swiech-2024,DeFeo-Swiech, Flandoli-1990, Meng-Shi-Zhang-Zhao-2025}.
This method always comes at a cost of moving problems to infinite dimension. Another limitation in the stochastic case lies in the appearance of unbounded operators in state equations \cite{Meng-Shi-Zhang-Zhao-2025}.

$\bullet$ Variation of constants formula method: By using this method, one can transform the original delayed system into another integral system of Volterra type without delay which is solved by fundamental solution. We refer to \cite{Huang-Li-Wang-2016,Lee-You-1989,Lee-You-1990,  Meng-Shi-Wang-Zhang-2023} for more details. This works well in deterministic controlled systems. However, when diffusion terms contain state delay terms, it becomes unclear to go through the whole procedure since anticipated stochastic integrals may be involved.

Based on the above analysis, different methods are used to handle different models, but each has its own drawbacks. Therefore, we raise the second question:

\it {\bf (Q2)}: Is it possible to come up with a new approach to bypass the disadvantages in the existing  papers? \rm

\ms

{\bf III.} Last but not the least, let us make some discussions about the obtained conclusions in terms of feedback controls and Riccati systems in the literature. We separate them into the following cases.

$\bullet$ Riccati systems versus Fredholm systems: To obtain feedback controls, Riccati systems are one of the most popular and important tools. We refer to \cite{Alekal-Brunovsky-Chyung-Lee-1971, Bagchi-1978, Carmona-Fouque-Mousavi-Sun-2018, Chen-Wu-2011, Chen-Zhang-2023, Delfour-1977-1, Delfour-1986, Lefebvre-Miller-2021, Meng-Shi-Yong-2023, Meng-Shi-Zhang-Zhao-2025, Uchida-Shimemura-Kubo-Abe-1988, Zhang-Xu-2017} and the related papers therein. On the other hand, according to the works \cite{Lee-You-1989,Lee-You-1990}, Fredholm systems also play important roles. In other words, the helpful systems for establishing the feedbacks are by no means unique.

$\bullet$ The challenges of Riccati systems: This point can be explained in the following manners. Firstly, the successful introduction of Riccati systems only hold in specific settings, e.g., the time-invariant case \cite{Ichikawa-1982, Lefebvre-Miller-2021,Meng-Shi-Zhang-Zhao-2025, Zhang-Xu-2017}, the only state delay case \cite{Chen-Zhang-2023,Eller-Aggarwal-Banks, Liang-Xu-Zhang-2018}, the only control delay case \cite{Carmona-Fouque-Mousavi-Sun-2018, Liang-Xu-Zhang-2018, Zhang-Xu-2017}. Secondly, even for the papers containing Riccati systems, verifying the coincidences among them seems quite involved and technical, such as \cite{Chen-Zhang-2023} and \cite{Meng-Shi-Zhang-Zhao-2025}, or \cite{Lefebvre-Miller-2021} and \cite{Zhang-Xu-2017}. Thirdly, even though Riccati systems are derived, they are quite challenging to further discuss their solvability. Along this line we mention the works of \cite{Chen-Zhang-2023, Meng-Shi-Zhang-Zhao-2025} for some progress in the stochastic setting. Fourthly, as to the papers obtaining the solvability, different works require different assumptions \cite{Alekal-Brunovsky-Chyung-Lee-1971,Delfour-1977}. Hence, it seems difficult to give a unified precondition.

$\bullet$ The challenges of feedback controls: In the literature, open-loop optimal controls and closed-loop optimal controls are two main notions for representing the explicit forms. For the former one, there are some relevant papers \cite{Chen-Wu-2011,Meng-Shi-Zhang-Zhao-2025,Xu-Shi-Zhang-2018, Zhang-2021,Zhang-Xiong-Shi-2021}. Since it is not easy to decouple forward-backward adjoint systems, the closed-loop representations of open-loop optimal controls are not successfully given until some recent works of \cite{Chen-Zhang-2023, Meng-Shi-Zhang-Zhao-2025, Zhang-Xu-2017}.
As to the later one, there are some works treating the closed-loop solvability. We refer to \cite{Meng-Shi-Yong-2023, Meng-Shi-Zhang-Zhao-2025}. However, their closed-loop solvability depends on the transformed problems, thus lacks generality, and their methods are not applicable to time-varying coefficients.

To sum up the above, different models have different versions of Riccati systems, among which the difference analysis is lacking and the relevant solvability is not obtained. Therefore, we come up with the third question:

\it {\bf (Q3)}: Is it possible to derive a unified Riccati system and feedback control for the general model \rf{state} and \rf{cost}, and verify the consistency with the existing papers? \rm

\ms

Now let us return back to Problem (P) and discuss another special case where $A_2, B_2,C_2,Q_2,R_2=0$. Here we then arrive at the LQ optimal control problem for (a class of) stochastic integro-differential systems (or stochastic differential systems with memory).
This type of optimal control problems naturally emerges in many different application scenarios, and is particularly common when studying the optimal performances of systems in response to specific inputs. In such cases, the systems' responses do not occur immediately but rather appear after a certain period of time.
However, in contrast with the previous discussion on SDDEs, for optimal control problems
involving systems with memory, just a few isolated results are available at the moment. We refer to
Carlier--Tahraoui \cite{Carlier-Tahraoui,Carlier-Tahraoui-2008} for the finite dimensional case  and  Cannarsa--Frankowska--Marchini \cite{Cannarsa-Frankowska-Marchini} for the infinite dimensional study.
In the LQ framework, Pandolfi \cite{Pandolfi-2018} studied a simple integro-differential model in the finite-dimensional space and obtained the optimal synthesis
via Riccati-type equations for the first time (commented by \cite{Acquistapace-Bucci-2025}).
Other enhanced results appeared recently in \cite{Gong-Wang-2025,Pandolfi-2025}. For the extension to the infinite dimension, we refer to  Acquistapace--Bucci \cite{Acquistapace-Bucci,Acquistapace-Bucci-2025}. These results hold in the deterministic controlled system framework. Due to the lack of relevant study on stochastic systems, we raise the fourth question:

\it {\bf (Q4)}: It is possible to provide some systematic study in stochastic integro-differential systems to fill the gaps left by the existing literature? \rm

\subsection{Contributions and novelties}

In this paper, our goal is to treat the general framework with the state equation \rf{state} and the cost functional \rf{cost} by employing new methodologies, and give positive answers to the aforementioned four questions. To begin with, we propose a general definition of the closed-loop solvability for Problem (P). To show its sufficiency, we separate the procedures into four parts. We firstly adapt the transformation procedures, introduced in Meng--Shi--Wang--Zhang \cite{Meng-Shi-Wang-Zhang-2023}, into our framework and end up with an LQ problem for a stochastic Volterra integral system without delay. Secondly, by borrowing the ideas developed in Hamaguchi--Wang \cite{Hamaguchi-Wang-2024-II} and Gong--Wang \cite{Gong-Wang-2025}, we introduce and discuss a class of Riccati systems and backward stochastic adjoint systems. Then, by the previous Riccati systems, we explicitly construct the desired closed-loop strategy, the resulting
closed-loop outcome control, and prove its optimality. Finally, we make detailed comparisons/coincidences with the existing study.
The contributions and innovations of this paper are summarized as follows.

\ss

\begin{itemize}
  \item A general yet new stochastic LQ optimal control problem is studied. On the one hand, it covers stochastic differential delay system where the coefficients are time-varying, both state delay and control delay can appear in the drift term, the diffusion term, and the cost functional. On the other hand, it covers a class of stochastic integro-differential equations where memory terms enter into the drift, diffusion terms, and the cost functional. This gives a nice response to (Q1), (Q4).

  \item A new transformation method is employed for the above general framework, and the advantages are shown in three aspects. Firstly, it avoids the complicated decoupling procedures for forward-backward systems. Secondly, it can handle time-varying coefficients case where traditional infinite dimensional lifting methods fail. Thirdly, in our approach there are no unbounded control operators, and there is no need to use infinite-dimensional analysis theories such as operator semigroups.   This gives a nice response to (Q2).

  \item By our main conclusions, we find the following four advantages and new facts. Firstly, we give the solvability of the Riccati system corresponding to Problem (P) and their coincidences with the existing literature in particular cases. Secondly, we explicitly construct the unique closed-loop outcome control, which does not rely on the future state and avoids complex tools of conditional expectations. Then, we only impose integrability conditions on the coefficients of Problem (P), and do not require continuity or even differentiability assumptions. Finally, even for the deterministic control system \rf{state}, we present, for the first time, results regarding the integro-differential part and the cost functional \rf{cost} with pointwise/distributed delays. This gives a nice response to (Q3), (Q4).

\end{itemize}

\ss

The rest part is organized as follows. In Section 2, we formulate the control problem studied in the paper. In Section 3, we transform the original delayed control system into a Volterra integral control system without delay, and then, study the closed-loop solvability of the original delayed optimal control problem. In Section 4, we discuss several important cases to compare our results with the previous ones. In Section 5, we give some concluding remarks. Finally, we provide the proofs of the main results in Appendix.

\section{Preliminary}

In this section, we give some preliminaries and formulate the control problem studied in this paper.

\ss

Let $T>0$ be a given finite time duration, define
$\triangle_2\left(0, T\right)\equiv\left\{(t, s) \in\left(0, T\right)^2 \mid T>t>s>0\right\}$, 
$\Box_3\left(0, T\right)\equiv\left\{\left(s_1, s_2, t\right) \in\left(0, T\right)^3 \mid t<(s_1 \wedge s_2)\right\}$. $I$ is the identity matrix with appropriate dimension.
$\dbS^n$ is the set of all $n\times n$ symmetric matrices. Next we define the following spaces which will be used in this paper.
Denote by $L^\i(0,T;\dbR^n)$ the Banach space consisting of $\dbR^n$-valued variables $\phi(\cd)$ such that $\sup\limits_{0\les t\les T}|\phi(t)|<\i$,
by $L^p(0,T;\dbR^n)$ the Banach space consisting of $\dbR^n$-valued variables $\phi(\cd)$ such that $\int_0^T|\phi(t)|^pdt<\i$, where $p$ is an integer.
Denote by $L^2_{\cF_t}(\Omega;\dbR^n)$ the Hilbert space consisting of $\dbR^n$-valued $\mathcal{F}_t$-measurable random variables $\p$ such that $\mathbb{E}|\p|^2<\infty$,  by $L^2_{\dbF}(\Omega;C([0,T];\dbR^n))$ the Banach space consisting of $\dbR^n$-valued $\mathbb{F}$-adapted continuous processes $\phi(\cd)$ such that $\mathbb{E}\big[\sup\limits_{0\les t\les T}|\phi(t)|^2\big]<\infty$,
%
by $L_{\mathbb{F}}^{2,p}\left(\triangle_2\left(0,T\right); \mathbb{R}^{n\times n}\right)$ the Banach space of $\mathbb{R}^{n \times n}$-valued and measurable processes $\p$ on $\triangle_2\left(0, T\right)$ such that $\p(t,\cdot)$ is $\mathbb{F}$-progressively measurable on $\left(0,t\right)$ for each $t\in\left(0,T\right)$, and $
\mathbb{E}\big[\int_{0}^T\left(\int_{0}^t|\p(t, s)|^p ds\right)^{\frac{2}{p}} dt\big]^{\frac 1 2}<\i.$ For $p=2$, we simply denote $L_{\mathbb{F}}^2\left(\triangle_2\left(0,T\right); \mathbb{R}^{n \times n}\right)\equiv L_{\mathbb{F}}^{2,2}\left(\triangle_2\left(0,T\right); \mathbb{R}^{n \times n}\right)$.
Denote by $\mathscr{L}^2\left(\triangle_2\left(0,T\right); \mathbb{R}^{n \times n}\right)$ the set of $\p \in L^2\left(\triangle_2\left(0,T\right); \mathbb{R}^{n \times n}\right)$ satisfying
$\underset{t\in\left(0,T\right)}{\operatorname{ess}\sup} \left(\int_t^T|\p(s, t)|^2 ds\right)^{\frac 1 2}\3n<\2n\infty$, and for any $\varepsilon>0$, there exists a finite partition $\left\{a_i\right\}_{i=0}^m$ of $\left(0,T\right)$ with $0=a_0<a_1<\cdots<a_m=T$ such that $
\underset{t \in\left(a_i, a_{i+1}\right)}{\operatorname{ess} \sup }\3n\left(\int_t^{a_{i+1}}\2n|\p(s,\1nt)|^2 ds\right)^{\frac 1 2}\2n<\2n\varepsilon$,
for each $i \in\{0,1, \ldots, m\1n-\1n1\}$.
Denote by
$L^{2,2,1}\left(\Box_3(0, T) ; \mathbb{R}^{n \times n}\right)$ the Banach space of $\mathbb{R}^{n \times n}$-valued deterministic functions $\p$ on $\Box_3\1n\left(0,\1n T\right)$ such that $(\int_{0}^T\3n \int_{0}^T\2n(\int_{0}^{s_1 \1n\wedge\1n s_2}\1n\left|\p\1n\left(s_1,\1n s_2,\1n t\right)\right| \1nd t)^2 d s_1 d s_2\1n)^{1 / 2}\2n<\2n\infty.$

\ms

For any $0<t_0<T$, consider the following SDDE:
\bel{SDDE-}\left\{\ba{ll}
\ns\ds  d\bar x(t)=[\bar A_1(t)\bar x(t)+\bar A_2(t) \bar y(t)+\bar A_3(t)\bar z(t)+\bar b(t)]dt\\
\ns\ds\qq\qq+[\bar C_1(t)\bar x(t)+\bar C_2(t) \bar y(t)+\bar C_3(t)\bar z(t)+\bar \si(t)]dW(t),\q t\in(t_0,T),\\
\ns\ds \bar x(t)=\bar \xi(t-t_0),\ \ t\in[t_0-\d,t_0],
\ea\right.
\ee
where $\bar\xi(\cd)\in C([-\d,0];\dbR^n)$ is the initial trajectory of the state, $\d>0$ is the constant delay time,
$$\ba{ll}
\ns\ds
\bar y(t)\equiv \bar x(t-\d),\ \ \bar z(t)\equiv \int_{t_0}^t \bar F(t,s)\bar x(s)ds,\q t\in(t_0,T).
\ea$$
As to the coefficients of \rf{SDDE-}, we impose the following conditions:

\ms

(H1) $\bar A_1(\cd), \bar A_2(\cd),\bar A_3(\cd)\in L^2(0,T;\dbR^{n\times n})$, $\bar C_1(\cd), \bar C_2(\cd),\bar C_3(\cd)\in L^\i(0,T;\dbR^{n\times n})$,

\ \
$\bar F(\cd,\cd)\in L^{2,1}(\D_2(0,T);\dbR^{n\times n})$, $\bar b(\cd)\in L_\dbF^2(\O;L^1(0,T;\dbR^n))$, $\bar\si(\cd)\in L^2_\dbF(0,T;\dbR^n)$.

\ms

The following proposition, whose proof is given in the Appendix, guarantees its solvability. In contrast with e.g. \cite{Chen-Wu-2011,Chen-Wu-Yu-2012,Chen-Zhang-2023, Zhang-Xiong-Shi-2021} where $\bar A_i(\cd)$ is assumed to be bounded, $i=1,2,3$, we slightly relax them into proper integrable conditions.

\bp{Solvability of SDDE}
Let {\rm (H1)} hold. Then, SDDE \rf{SDDE-} admits a unique solution $\bar x(\cd)\in L^2_{\dbF}\big(\Omega;C([t_0,T];\dbR^n)\big)$.
\ep

\ss
For later usefulness, for any $0\les t_0<T$, we consider the linear SVIE:
\begin{eqnarray}
\ds \check X(t)=\check\f(t)+\int_{t_0}^t\bigg[\check A(t,s)X(s)+\check b(t,s)\bigg]ds +\int_{t_0}^t\bigg[\check C(t,s)X(s)+\check \si(t,s)\bigg]dW(s),\q t\in(t_0,T), \label{SVIE-}
\end{eqnarray}
with the coefficients satisfying:

\ms

(H2) $\check\f(\cd)\in L^2_{\mathbb{F}}(0,T;\dbR^n)$, $\check A(\cd,\cd)\in L^2(\D_2(0,T);\dbR^{n\times n})$, $\check C(\cd,\cd)\in\sL^2\left(\triangle_2(0,T); \mathbb{R}^{n\times n}\right)$,

 \ \ \ $\check b(\cd,\cd)\in L_{\mathbb{F}}^{2,1}\left(\triangle_2(0,T); \mathbb{R}^n\right)$, $\check \sigma(\cd,\cd)\in L_{\mathbb{F}}^2\left(\triangle_2(0,T); \mathbb{R}^n\right)$.

\ms

Then, we present the following result as Lemma A.3 in \cite{Hamaguchi-Wang-2024-I}.
\bp{Solvability of SVIE}
Let {\rm (H2)} hold. Then, SVIE \rf{SVIE-} admits a unique solution $\check X(\cd)\in L^2_{\dbF}(t_0,T;\dbR^n)$.
\ep

\ss

Based on the above preparations, let us return back to the state equation \rf{state}, the coefficients of which satisfy the following assumptions.
$${\bf(A1)}\left\{\ba{ll}
\ns\ds A_1(\cd),A_2(\cd),A_3(\cd)\in L^2(0,T;\dbR^{n\times n}),\  B_1(\cd),B_2(\cd),B_3(\cd)\in L^\i(0,T;\dbR^{n\times m}),\ b(\cd)\in L_\dbF^2(\O;L^1(0,T;\dbR^n)),\\
\ns\ds C_1(\cd),C_2(\cd),C_3(\cd)\in L^\i(0,T;\dbR^{n\times n}),\q\ \  D_1(\cd)\in L^\i(0,T;\dbR^{n\times m}),\qq\q \si(\cd)\in L^2_\dbF(0,T;\dbR^n),\\
\ns\ds F(\cd,\cd)\in L^{\i}(\D_2(0,T);\dbR^{n\times n}), \qq\ \ \ \   \ti{F}(\cd,\cd)\in L^{\i}(\D_2(0,T);\dbR^{n\times m}), \ \ \xi(\cd)\in C([-\d,0];\dbR^n),\\
\ns\ds Q_1(\cd),Q_2(\cd),Q_3(\cd)\in L^\i(0,T;\dbS^n),\q\ \ R_1(\cd),R_2(\cd)\in L^\i(0,T;\dbS^{m}),\ \ \varsigma(\cd)\in L^2(-\d,0;\dbR^m).
\ea\right.$$
In terms of Proposition \ref{Solvability of SDDE}, SDDE \rf{state} admits a unique solution, the cost functional \rf{cost} is well-defined, and hence, it becomes natural to pose Problem (P) in the Introduction.

%
%
%

\ms

Due to the particular LQ structure, we are interested in the closed-loop optimal control. To this end,
we first look at the closed-loop strategy. For any given $t_0\in[0,T)$, define
$$
\ba{ll}
\ns\ds \dbL\equiv L^2(t_0,T;\dbR^{m\times n}) \times L^2(\D_2(t_0,T);\dbR^{m\times n})\times L^\i(t_0,T;\dbR^{m\times n}) \times L^2(\D_2(t_0,T);\dbR^{m\times m})
\times L_\dbF^2(t_0,T;\dbR^m).
\ea
$$
In the following, for any $(K_1(\cd),K_2(\cd,\cd),K_3(\cd),K_4(\cd,\cd),v(\cd))\in \dbL$, we call it a \emph{closed-loop strategy} on $[t_0,T]$. Later we will use this closed-loop strategy, which does not depend on the initial data but only the given coefficients, to construct a closed-loop control on $[t_0,T]$. To introduce the closed-loop state, let us consider the SDDE:
\bel{closed-loop system}\5n\left\{\ba{ll}
\ns\ds  dx(t)=\big[A_1(t)x(t)+A_2(t) y(t)+A_3(t)z(t)+B_1(t)u(t)+B_2(t)\nu(t)+B_3(t)\m(t)+b(t)\big]dt\\
\ns\ds\qq\qq+\big[C_1(t)x(t)+C_2(t) y(t)+C_3(t)z(t)+D_1(t)u(t)+\si(t)\big]dW(t),\q t\in(t_0,T),\\
\ns\ds x(t)=\xi(t-t_0),\ \ \ u(t)=\varsigma(t-t_0),\ \ t\in[t_0-\d,t_0],\\
\ns\ds u(t)=K_1(t)x(t)+\int_{t_0}^tK_2(t,s)x(s)ds+K_3(t)x(t-\d)
+\int_{t_0}^tK_4(t,s)u(s)ds+v(t),\ \ t\in(t_0,T).
\ea\right.
\ee
\no We call \rf{closed-loop system} a \emph{closed-loop system} under $(K_1(\cd),K_2(\cd,\cd),K_3(\cd),K_4(\cd,\cd), v(\cd))$ corresponding to $(t_0,\xi(\cd),\varsigma(\cd))$, and call $x(\cd)$, $u(\cd)$ the corresponding \emph{closed-loop state} and \emph{closed-loop outcome control}, respectively. In the system \rf{closed-loop system}, $K_1(\cd)$, $K_2(\cd,\cd)$, $K_3(\cd)$, $K_4 (\cd,\cd)$ represent respectively the gain coefficients of the current state, the distributed state delay, the pointwise state delay and the distributed control delay.
The proof of the following result is given in the Appendix.

\bp{solvability of closed-loop system}
Under Assumption {\rm (A1)}, for any given $t_0\in[0,T)$ and $(K_1(\cd),K_2(\cd,\cd),K_3(\cd),K_4(\cd,\cd), v(\cd))$ $\in \dbL$, the closed-loop system \rf{closed-loop system} admits a unique solution $x(\cd)\in L^2_{\dbF}\big(\Omega;C([t_0,T];\dbR^n)\big)$ on $[t_0,T]$, and there exists a constant $L>0$ such that
\begin{eqnarray}
\dbE\sup\limits_{t_0\les t\les T}| x(t)|^2+\dbE\int_{t_0}^T| u(t)|^2dt\qq\qq\qq\qq\qq\qq\qq\qq\qq\qq\q
\notag\\
 \les\1n L\bigg\{\1n\sup\limits_{t_0-\d\les t\les t_0}
\1n| \xi(t-t_0)|^2\1n+\3n\int_{t_0-\d}^{t_0}\2n| \varsigma(t-t_0)|^2dt+\dbE\1n\int_{t_0}^T\3n\big(| b(t)|^2\1n+\1n| \si(t)|^2\1n+\1n|v(t)|^2\big)dt
\bigg\}.\label{pre-estimate-of-controlled-SDDE}
\end{eqnarray}
\ep

At last, let us introduce the following notions.

\bde{def optimal closed-loop strategy}
For any given $t_0\in[0,T)$, the closed-loop strategy $(K_1^*(\cd),K_2^*(\cd,\cd),K_3^*(\cd),K_4^*(\cd,\cd), v^*(\cd))\in \dbL$ is called an \emph{optimal closed-loop strategy} of Problem {\rm (P)} if
\bel{equiv-iii}
J\big(t_0,\xi,\varsigma;(K_1^*,K_2^*,K_3^*, K_4^*, v^*)\big) \les J\left(t_0,\xi,\varsigma;u\right),
\ee
for any $\left(\xi,\varsigma\right) \in C([-\d,0];\dbR^n)\times L^2(-\d,0;\dbR^m)$ and any control $u(\cd)\in L_\dbF^2(t_0,T;\dbR^m)$. If there (uniquely) exists an optimal closed-loop strategy on $[t_0,T]$, Problem {\rm (P)} is said to be \emph{(uniquely) closed-loop solvable on $[t_0,T]$}.
\ede

 \br{definition of closed-loop}
Inspired by \cite{Alekal-Brunovsky-Chyung-Lee-1971, Delfour-Mitter-1972, Ichikawa-1982,Sun-Li-Yong-2016, Sun-Yong-2014,Zhang-Xu-2017}, the closed-loop strategy for deterministic Problem {\rm (P)} was introduced in \cite{Meng-Shi-Yong-2023}, and then extended to the stochastic setting in \cite{Meng-Shi-Zhang-Zhao-2025}. Here our version is more general due to the appearance of $K_3(\cd)$. In addition, we do not have the transformed equivalent problems appeared in the definitions of the closed-loop strategies in \cite{Meng-Shi-Yong-2023,Meng-Shi-Zhang-Zhao-2025}.
\er

\section{The closed-loop solvability of Problem (P)}

In this section, we firstly transform the delayed Problem (P) into another optimal control problem driven by a (finite dimensional) SVIE without delay. Based on such a transformation, we then introduce and discuss the desired Riccati system for Problem (P). Eventually, we construct an explicit closed-loop strategy and prove its optimality in the sense of Definition \ref{def optimal closed-loop strategy}.

\subsection{Problem transformation}

To begin with, let us take a closer look at the state $x(\cd)$.
Since $\varsigma(\cd)$ is the initial trajectory of the control, for either $t\in\big[t_0,(t_0+\d)\wedge T\big]$ or
$t\in \big((t_0+\d)\wedge T,t_0+\d\big]$ we get
$$\ba{ll}
\ns\ds \int_{t_0}^tB_2(s)u(s-\d)ds=\int_{t_0}^{t} B_2(s)\varsigma(s-\d-t_0)ds
=\int_{t_0}^{t_0+\d}{\bf 1}_{[t_0,t)}(s) B_2(s)\varsigma(s-\d-t_0)ds.
\ea$$
For $t\in\big(t_0+\d,(t_0+\d)\vee T\big]$, we have
$$\ba{ll}
\ns\ds  \int_{t_0}^tB_2(s)u(s-\d)ds=\int_{t_0}^{t_0+\d } B_2(s)u(s-\d)ds +\int_{t_0+\d }^{t} B_2(s)u(s-\d)ds \\
\ns\ds =\int_{t_0}^{t_0+\d } {\bf 1}_{[t_0,t)}(s)B_2(s)\varsigma(s-\d-t_0)ds
+\int_{t_0}^{t} {\bf 1}_{[t_0 ,t-\d)}(s)B_2(s+\d)u(s)ds.
\ea
$$
To sum up, for $t\in[t_0,T]$, we have
$$\ba{ll}
\ns\ds \int_{t_0}^tB_2(s)u(s-\d)ds=\int_{t_0}^{t_0+\d } {\bf 1}_{[t_0,t)}(s)B_2(s)\varsigma(s-\d-t_0)ds
+\int_{t_0}^{t} {\bf 1}_{[t_0 ,t-\d)}(s)B_2(s+\d)u(s)ds.
\ea
$$
As to the $\mu(\cd)$ term in $x(\cd)$, by the Fubini theorem, for $t\in[t_0,T]$, we obtain
$$\ba{ll}
\ns\ds \int_{t_0}^t B_3(s)\mu(s)ds=\int_{t_0}^tB_3(s)\int_{t_0}^s\tilde{F}(s,r)u(r)drds =\int_{t_0}^t\int_r^tB_3(s)\tilde{F}(s,r)dsu(r)dr.
\ea$$
Thus we deduce
\begin{eqnarray}
x(t)=\xi(0)+\int_{t_0}^{t_0+\d }B_2(s)\varsigma(s-\d-t_0) {\bf1}_{[t_0,t)}(s)ds+\int_{t_0}^t\Big[A_1(s)x(s)+A_2(s)y(s)
\qq\qq\notag\\
+A_3(s)z(s)+\Big(B_1(s)+B_2(s+\d){\bf1}_{[t_0,t-\d)}(s)
+\int_s^tB_3(r)\tilde{F}(r,s)dr\Big)u(s)+b(s)\Big]ds\ \ \notag\\
+\int_{t_0}^t\big[C_1(s)x(s)+C_2(s)y(s)+C_3(s)z(s)+D_1(s)u(s)
+\sigma(s)\big]dW(s),\q\!  t\in[t_0,T].\ \label{deduced state}
\end{eqnarray}
Next let us turn to the term of $y(\cd)$.
For $t\in\big[t_0+\d,(t_0+2\d)\wedge T\big]$ and $t\in\big((t_0+2\d)\wedge T,t_0+2\d\big]$, we have
$$\ba{ll}
\ns\ds  \int_{t_0}^{t-\d}B_2(s)u(s-\d)ds=\int_{t_0}^{t-\d} B_2(s)\varsigma(s-\d-t_0)ds =
\int_{t_0}^{t_0+\d} {\bf 1}_{[t_0,t-\d)}(s)B_2(s)\varsigma(s-\d-t_0)ds.
\ea
$$
For $t\in\big(t_0+2\d,(t_0+2\d)\vee T\big]$, we have
$$\ba{ll}
\ns\ds  \int_{t_0}^{t-\d}B_2(s)u(s-\d)ds=\int_{t_0}^{t_0+\d} B_2(s)u(s-\d)ds +\int_{t_0+\d}^{t-\d} B_2(s)u(s-\d)ds \\
\ns\ds =\int_{t_0}^{t_0+\d} {\bf 1}_{[t_0,t-\d)}(s)B_2(s)\varsigma(s-\d-t_0)ds
+\int_{t_0}^{t} {\bf 1}_{[t_0,t-2\d)}(s)B_2(s+\d)u(s)ds.
\ea
$$
Therefore, for $t\in[t_0+\d,T]$,
$$\ba{ll}
\ns\ds \int_{t_0}^{t-\d}B_2(s)u(s-\d)ds{\bf1}_{(t_0+\d,\i)}(t) =\Big[\int_{t_0}^{t_0+\d} {\bf 1}_{[t_0,t-\d)}(s)B_2(s)\varsigma(s-\d-t_0)ds\\
\ns\ds \qq\qq\qq\qq\qq\qq\qq\qq
+\int_{t_0}^{t} {\bf 1}_{[t_0,t-2\d)}(s)B_2(s+\d)u(s)ds\Big] {\bf1}_{(t_0+\d,\i)}(t).
\ea
$$
As to the $\mu(\cd)$ term in $y(\cd)$,
for $t\in[t_0+\d,T]$, it follows from the Fubini theorem that
$$\ba{ll}
\ns\ds \int_{t_0}^{t-\d}B_3(s)\mu(s)ds=\int_{t_0}^{t-\d}B_3(s)\int_{t_0}^s\tilde{F}(s,r)u(r)drds =\int_{t_0}^{t-\d}B_3(r)\int_{t_0}^r\tilde{F}(r,s)u(s)dsdr\\
\ns\ds
=\int_{t_0}^{t-\d}\int_s^{t-\d}B_3(r)\tilde{F}(r,s)dru(s)ds.
\ea$$
Hence for the pointwise state delay, we have
\bel{pointwise delay of state}\ba{ll}
\ns\ds
y(t)=\xi(t-\d-t_0){\bf1}_{[t_0,t_0+\d]}(t)+{\bf1}_{(t_0+\d,\i)} (t)\Big\{
\xi(0)
+\int_{t_0}^{(t_0+\d)\wedge(t-\d)}B_2(s)\varsigma(s-t_0-\d)ds\Big\} \\
\ns\ds \qq\q+\int_{t_0}^t{\bf1}_{[t_0,t-\d)}(s)\Big(A_1(s)x(s)+A_2(s) y(s)+A_3(s)z(s)
+\Big[B_1(s)+B_2(s+\d){\bf1}_{[t_0,t-2\d)}(s)\\
\ns\ds\qq\q +\int_s^{t-\d} B_3(r)\ti F(r,s)dr\Big]u(s)+b(s)\Big)ds +\int_{t_0}^t{\bf1}_{[t_0,t-\d)}(s)\big[C_1(s)x(s)+
C_2(s) y(s)
\\
\ns\ds\qq\q +C_3(s)z(s)+D_1(s)u(s)+\si(s)\big]dW(s),\q t\in[t_0,T].
\ea
\ee
Eventually, let us treat the term of $z(\cd)$.
For $t\in[t_0,T]$,
\begin{eqnarray}
z(t)=\int_{t_0}^tF(t,s)x(s)ds=\int_{t_0}^tF(t,s)
\Big[\xi(0)+\int_{t_0}^s\big(A_1(r)x(r)+A_2(r)y(r) +A_3(r)z(r)\notag\\
+B_1(r)u(r)+B_2(r)\nu(r)+B_3(r)\mu(r)+b(r)\big)dr\notag \qq\qq\qq\qq\qq\q\ \\
+\int_{t_0}^s\big(C_1(r)x(r)+C_2(r)y(r) +C_3(r)z(r)+D_1(r)u(r)
+\sigma(r)\big)dW(r)\Big]ds.\label{Z1}
\end{eqnarray}
By Fubini theorem,
$$\ba{ll}
\ns\ds \int_{t_0}^t\int_{t_0}^sF(t,s)B_2(r)\nu(r)drds =\int_{t_0}^t\Big[\int_r^tF(t,s)ds\Big] B_2(r)u(r-\d)dr.
\ea$$
Notice that for $t\in[ t_0,t_0+\d)$,
$$\ba{ll}
\ns\ds \int_{t_0}^t\int_{t_0}^sF(t,s)B_2(r)\nu(r)drds =\int_{t_0}^{ t}\Big[\int_r^tF(t,s)ds\Big]B_2(r)\varsigma(r-\d-t_0)dr,
\ea$$
while for $t\in\big[t_0+\d,(t_0+\d)\vee T\big)$, one has
$$\ba{ll}
\ns\ds \int_{t_0}^t\int_{t_0}^sF(t,s)B_2(r)\nu(r)drds =\int_{t_0}^t\Big[\int_r^tF(t,s)ds\Big] B_2(r)u(r-\d)dr\\
\ns\ds =\int_{t_0}^{t_0+\d }\Big[\int_r^tF(t,s)ds\Big]B_2(r)\varsigma(r-\d-t_0)dr
+\int_{t_0}^{t-\d}\Big[\int_{r+\d}^tF(t,s)ds\Big]B_2(r+\d) u(r)dr.
\ea$$
Therefore, for $t\in[t_0,T]$, we obtain that
\begin{eqnarray}
\int_{t_0}^t\int_{t_0}^sF(t,s)B_2(r)\nu(r)drds
\qq\qq\qq\qq\qq\qq\qq\qq\qq\qq\qq\qq\qq\qq \notag\\
=\int_{t_0}^{(t_0+\d)\wedge t }\Big[\int_r^tF(t,s)ds\Big]B_2(r)\varsigma(r-\d-t_0)dr
+\int_{t_0\wedge (t-\d)}^{t-\d}\Big[\int_{r+\d}^tF(t,s)ds\Big]B_2(r+\d) u(r)dr.\label{Z2}
\end{eqnarray}
As to the $\mu(\cd)$ term in $z(\cd)$, by the Fubini theorem again, for $t\in[t_0,T]$,
\begin{eqnarray}
\int_{t_0}^t\int_{t_0}^sF(t,s)B_3(r)\mu(r)drds
=\int_{t_0}^t\int_{t_0}^sF(t,s)B_3(r)\Big[\int_{t_0}^r\tilde{F}(r,\a)
u(\a)d\a\Big] drds\qq\qq\qq\q\notag
\\
=\int_{t_0}^t\int_{t_0}^s\int_\a^s F(t,s)B_3(r)\tilde{F}(r,\a)
u(\a)drd\a ds=\int_{t_0}^t\Big[\int_\a^t\int_\a^sF(t,s) B_3(r)\tilde{F}(r,\a)drds\Big] u(\a)d\a.\label{Z3}
\end{eqnarray}
Hence, from \rf{Z1}--\rf{Z3} and by adding some indicative functions, we derive
\bel{distributed delay of state}\!\!\ba{ll}
\ns\ds z(t)=\!\int_{t_0}^t\! F(t,s)\Big(\xi(0)+\!\int_{t_0}^{t_0+\d}\! B_2(r)\varsigma(r-\d-t_0){\bf1}_{[0,s)}(r)dr
\Big)ds
\!+\!\int_{t_0}^t\cE(t,s)
\Big[A_1(s)x(s)\!+\! A_2(s) y(s)
\\
\ns\ds\qq +A_3(s)z(s)\Big]ds\!+\!\int_{t_0}^t\Big[\!\int_s^tF(t,r)\Big(B_1(s)\!+\! B_2(s+\d){\bf1}_{[t_0,r-\d)}(s)
\!+\!\int_s^rB_3(\th)\ti F(\th,s)d\th\Big)dr\Big] u(s)ds \\
\ns\ds\qq+\int_{t_0}^t \cE(t,s)b(s)ds \!+\!\int_{t_0}^t\cE(t,s)
\big[C_1(s)x(s)\!+\!C_2(s) y(s)
\!+\!C_3(s)z(s)\!+\!D_1(s)u(s)\!+\!\si(s)\big]dW(s).
\ea
\ee
\vskip-2mm
\no Here and next, for $T\ges t>s\ges t_0$, we define $\cE(\cd,\cd)$ and  some other notations as follows:
\begin{eqnarray}
X(t)\equiv\begin{bmatrix}\ba{c}
x(t)    \\
y(t)   \\
z(t)
\ea\end{bmatrix},\qq\qq\q \cE(t,s)\equiv\int_s^tF(t,r)dr{\bf1}_{[0,t]}(s),\ \
\qq\qq\qq\qq\q\ \ \ \notag\\
\f(t)\equiv\begin{bmatrix}\begin{array}{c}
 \xi(0)+\int_{t_0}^{t_0+\d} B_2(s)\varsigma(s-t_0-\d){\bf 1}_{[t_0,t)}(s)ds   \\[1mm]
 \xi(t-t_0-\d){\bf 1}_{[t_0,t_0+\d]}(t)+[\xi(0)+\int_{t_0}^{(t_0+\d)\wedge(t-\d)} B_2(s) \varsigma(s-t_0-\d)ds]{\bf 1}_{(t_0+\d,\infty)}(t)  \\[1mm]
 \int_{t_0}^tF(t,s)[\xi(0)+\int_{t_0}^{t_0+\d} B_2(r)\varsigma(r-t_0-\d){\bf 1}_{[t_0,s)}(r)dr]ds
\end{array}\end{bmatrix},\qq\q\ \ \notag\\
A(t,s)\equiv\begin{bmatrix}\begin{array}{ccccc}
 A_1(s) &  A_2(s) & A_3(s) \\
{\bf 1}_{(\d,\infty)}(t-s) A_1(s)  & {\bf 1}_{(\d,\infty)}(t-s) A_2(s) &  {\bf 1}_{(\d,\infty)}(t-s) A_3(s)\\
\cE(t,s)A_1(s)& \cE(t,s)A_2(s)& \cE(t,s)A_3(s)
\end{array}
\end{bmatrix},\notag\qq\qq\qq\qq\qq\q\ \ \\
B(t,s)\equiv\begin{bmatrix}\begin{array}{ccccc}
 B_1(s)+B_2(s+\d){\bf 1}_{(\d,\i)}(t-s)+\int_s^tB_3(r)\ti F(r,s)dr  \\
{\bf 1}_{(\d,\infty)}(t-s) [B_1(s)+B_2(s+\d){\bf 1}_{(2\d,\i)}(t-s)+\int_s^{t-\d}B_3(r)\ti F(r,s)dr ]\\
\cE(t,s) B_1(s)+\cE(t,s+\d)B_2(s+\d)+\int_s^t\cE(t,\th)B_3(\th)\ti F(\th,s)d\th
\end{array}
\end{bmatrix},\qq\qq\q\qq\ \ \label{B(t,s)}\\
C(t,s)\equiv\begin{bmatrix}\begin{array}{ccccc}
C_1(s) &  C_2(s)  & C_3(s)  \\
{\bf 1}_{(\d,\infty)}(t-s) C_1(s)  & {\bf 1}_{(\d,\infty)}(t-s) C_2(s) & {\bf 1}_{(\d,\infty)}(t-s) C_3(s)\\
\cE(t,s)C_1(s)& \cE(t,s)C_2(s)& \cE(t,s)C_3(s)
%
\end{array}
\end{bmatrix},\qq\qq\qq\qq\qq\q\ \ \notag\\
 D(t,s)\equiv\begin{bmatrix}\begin{array}{ccccc}
D_1(s)  \\
{\bf 1}_{(\d,\infty)}(t-s) D_1(s)  \\
\cE(t,s) D_1(s)
\end{array}\end{bmatrix},\
\ti b(t,s)\2n\equiv\2n\begin{bmatrix}\1n\begin{array}{ccccc}
b(s)  \\
{\bf 1}_{(\d,\infty)}(t-s) b(s)  \\
\cE(t,s) b(s)
\end{array}\1n\end{bmatrix},\  \ti\si(t,s)\2n\equiv\2n\begin{bmatrix}\1n\begin{array}{ccccc}
\si(s)  \\
{\bf 1}_{(\d,\infty)}(t-s) \si(s)  \\
\cE(t,s) \si(s)
\end{array}\1n\end{bmatrix}.\qq\ \ \notag
\end{eqnarray}
Based on \rf{deduced state}, \rf{pointwise delay of state} and \rf{distributed delay of state}, $X(\cd)$ satisfies the following SVIE:
\begin{eqnarray}
\ds X(t)=\f(t)+\int_{t_0}^t\big[A(t,s)X(s)+B(t,s)u(s)+\ti b(t,s)\big]ds\qq\qq\q\  \notag\\
\ns\ds \qq\q+\int_{t_0}^t\big[C(t,s)X(s)+D(t,s)u(s)+\ti \si(t,s)\big]dW(s),\q t\in(t_0,T). \label{transformed state}
\end{eqnarray}
Denote
$$\ba{ll}
\ns\ds Q(t)\equiv\begin{bmatrix}\begin{array}{ccccc}
Q_1(t) & 0&0\\
0 & Q_2(t)&0\\
0 & 0& Q_3(t)
\end{array}
\end{bmatrix},\ \  R(t)\equiv R_1(t)+R_2(t+\d){\bf 1}_{[0,T-\d)}(t).
\ea
$$
Then, the cost functional \rf{cost} becomes
\begin{eqnarray}
J(t_0,\xi(\cd),\varsigma(\cd);u(\cd))=\mathbb{E}\int_{t_0}^T\big[X(t)^\top Q(t)X(t)+u(t)^\top R(t)u(t)\big]dt\qq\q\notag\\
+\mathbb{E}\int_{t_0}^{(t_0+\d)\wedge T}\varsigma(t-t_0-\d)^\top
R_2(t)\varsigma(t-t_0-\d)dt.
\label{transformed cost}
\end{eqnarray}
To sum up the above arguments, we have the following result whose proof is given in the Appendix.
\bp{prop coefficients condition}
Let Assumption {\rm(A1)} hold. Then, $u^*\in L^2_{\dbF}(t_0,T;\dbR^m)$ is an optimal control of Problem {\rm (P)} if and only if $ u^*(\cd)$ minimizes \rf{transformed cost} subject to SVIE \rf{transformed state}.
\ep

\br{equivalence of problem}
The idea of transforming the delayed system into another one without delay is popular in the existing literature.

For example, in \cite{Delfour-1986, Delfour-Mitter-1972, Meng-Shi-Yong-2023, Meng-Shi-Zhang-Zhao-2025}, they transformed the original delayed systems into infinite-dimensional evolution control systems without delay. However, in some cases they had to treat the unbounded control operator which would bring essential difficulties. More importantly, it seems quite complex to transform the obtained operator-valued Riccati equation into the (finite dimensional) matrix-valued case.

While in this paper, we put $(x(\cd),y(\cd),z(\cd))$ together to construct a new state $X(\cd)$ and then transform the original system equivalently to a controlled SVIE. Similar procedure also appeared in \cite{Meng-Shi-Wang-Zhang-2023}. On the one hand, this helps us to utilize the developed LQ theory for SVIEs in \cite{Hamaguchi-Wang-2024-II}. On the other hand, the advantage lies in that $X(\cd)$ is still finite dimensional which helps us to bypass the complicated operator language.

Another approach to addressing delayed problems is direct decoupling the forward delayed systems and the backward anticipated systems. However, this method requires the application of It\^o formula for delayed processes, which constitutes the essential difficulty of delayed problems. In contrast, the method proposed in this paper does not need to tackle such a challenge.

We point out that the transformation to Volterra systems also appeared in e.g. \cite{Huang-Li-Wang-2016, Lee-You-1990} via the constant variation formula. Nevertheless, it is not clear whether such an approach still works when the diffusion term in \rf{state} contains the pointwise state delay and distributed state delay.
\er

\subsection{The solvability of the new Riccati system}

To study the closed-loop solvability of Problem (P), we introduce a new Riccati system based on the above transformation and the existing study in \cite{Hamaguchi-Wang-2024-II}. To this end, we need some preparations.

\ms

Firstly, we denote by $\Pi(0, T)$ the set of pairs $P=\left(P^{(1)}, P^{(2)}\right)$ with $P^{(1)}:(0, T) \rightarrow \mathbb{R}^{(3n) \times (3n)}$ and $P^{(2)}: \Box_3\left(t_0, T\right) \rightarrow \mathbb{R}^{(3n) \times (3n)}$ such that

\ms

{\rm(i)} $P^{(1)} \in L^{\infty}\left(0, T ; \mathbb{S}^{3n}\right)$;

{\rm(ii)} for a.e. $\left(s_1, s_2\right) \in(0, T)^2, t \mapsto P^{(2)}\left(s_1, s_2, t\right)$ is absolutely continuous on $\left(0, s_1 \wedge s_2\right)$;

{\rm(iii)} $(s_1,s_2)\mapsto P^{(2)}\left(s_1, s_2, s_1 \wedge s_2\right)\equiv \lim\limits _{t \uparrow (s_1 \wedge s_2)} P^{(2)}\left(s_1, s_2, t\right)$ belongs to $L^2\left((0, T)^2 ; \mathbb{R}^{(3n) \times (3n)}\right)$;

{\rm(iv)}  $(s_1, s_2, t)\mapsto \dot{P}^{(2)}\left(s_1, s_2, t\right)\equiv \frac{\partial P^{(2)}}{\partial t}\left(s_1, s_2, t\right)$ belongs to $L^{2,2,1}\left(\Box_3(0, T) ; \mathbb{R}^{(3n) \times (3n)}\right)$;

{\rm(v)} for a.e. $\left(s_1, s_2, t\right) \in \Box_3(0, T)$, it holds that $P^{(2)}\left(s_1, s_2, t\right)=P^{(2)}\left(s_2, s_1, t\right)^{\top}$.

\ms

Secondly, for later notational usefulness, we define several coefficients as follows:
\bel{Pi}\ba{ll}
\ns\ds
\Pi(s,t,\th)\equiv \begin{bmatrix}
\frac{1}{s-t}I& \1n\frac{1}{s-t}{\bf1}_{(\d,\i)}(s-t)I&\1n I
\\
\frac{1}{s-t}{\bf1}_{(\d,\i)}(s-t)I&\1n \frac{1}{s-t}{\bf1}_{(2\d,\i)}(s-t)I&\1n {\bf1}_{(\d,\i)}(s-\th)I
\\
\frac{1}{s-t}\cE(s,t)
&\1n \frac{1}{s-t}\cE(s,t+\d)&\1n \cE(s,\th)
\end{bmatrix},
\Upsilon(s,t)
\equiv\begin{bmatrix}
I\\
{\bf1}_{(\d,\i)}(s-t)I\\
\cE(s,t)
\end{bmatrix},\\
\ns\ds \cA(t)\equiv \big(A_1(t),A_2(t),A_3(t)\big),\ \ \cB(\th,t)\equiv\begin{bmatrix}
B_1(t)\\
B_2(t+\d)\\
B_3(\th)\wt F(\th,t)
\end{bmatrix},\ \ \cC(t)\equiv \big(C_1(t),C_2(t),C_3(t)\big).
\ea\ee
Given these coefficients, we introduce the following system satisfied by a pair of $(P^{(1)},P^{(2)})$:
\bel{eq_Riccati--Volterra}\5n\2n\left\{\!\!\!\ba{ll}
\ns\ds	P^{(1)}(t)=Q(t)+\cC(t)^\top \cG_1(t)\cC(t) -\cC(t)^\top \cG_1(t)D_1(t)\cR(t)^{-1} D_1(t)^\top\cG_1(t)\cC(t),0< t< T,\\
\ns\ds	P^{(2)}(s,t,r)\!=\!P^{(2)}\1n(s,t,t\wedge s)
\1n\!-\1n\!\!\int_r^{t\wedge s}\1n\!\!\!\int_\t^T\1n \!\!\!\int_\t^T\2n\! \!\cG_3(s,\t,\th_1)\cB(\th_1,\t)\cR(\t)^{-1}
\cB(\th_2,\t)^\top\\
\ns\ds \qq\qq\qq\times\cG_3(t,\t,\th_2)^\top d\th_1 d\th_2 d\t\1n,0\1n<\1n r<\1n (s\1n\wedge\1n t)\1n<\1n T,\\
\ns\ds
	P^{(2)}(\bar s,t,t)=P^{(2)}\1n(t,\bar s,t)^\top \3n=\1n\cG_2(\bar s,t)\cA(t)\1n-\! \1n\int_t^T \3n\1n\cG_3(\bar s,t,\th)\cB(\th\1n,t)
 \cR(t)^{-1}\1n   D_1(t)^\top\1n\cG_1(t)\cC(t)d\th,0\1n<\1n t\1n<\1n\bar s\1n<\1n T,\\
\ea\right.\ee
%
%
\no where for  $0<t< ( s\wedge \bar s\wedge \th)< T$,
\bel{Upsilon-cR}\ba{ll}
\ns\ds \cR(t)\equiv R_1(t)+R_2(t+\d){\bf1}_{[0,T-\d)}(t)+D_1(t)^\top \cG_1(t)D_1(t),\\
\ns\ds \cG_1(t)\!\equiv\! \cG_1\big(t;P^{(1)},P^{(2)}\big)\!\equiv\!\!\int_t^T\!\Upsilon(s,t)^\top \!P^{(1)}(s)\Upsilon(s,t) ds\1n +\!\!\int_t^T\3n\int_t^T\3n\Upsilon(s_1,t)^\top \!P^{(2)}(s_1,s_2,t)\Upsilon(s_2,t) ds_1ds_2, \\
\ns\ds \cG_2(\bar s,t)\equiv\cG_2\big(\bar s,t;P^{(1)},P^{(2)}\big)\equiv P^{(1)}(\bar s)\Upsilon(\bar s,t)+\int_t^TP^{(2)}(\bar s,r,t)
\Upsilon(r,t)dr, \\
\ns\ds  \cG_3(s,t,\th)\equiv  \cG_3\big(s,t,\th;P^{(1)},P^{(2)}\big)\equiv P^{(1)}(s){\bf1}_{(t,s)}(\theta)\Pi(s,t,\th)+\int_\th^T
P^{(2)}(s,r,t)\Pi(r,t,\th)dr.
\ea\ee
In this paper, we name \rf{eq_Riccati--Volterra} the desired Riccati system, explicitly depending on $A_i$, $B_i$, $C_i$, $D_1$, based on (at least) the following aspects. Firstly, it is consistent with the Riccati system for the stochastic Volterra system in \cite{Hamaguchi-Wang-2024-II}. Secondly, in particular cases we will show in Subsection 4.1, 4.2, 4.5 that it can reduce to
the Riccati systems in the existing literature. Thirdly, just like the existing Riccati systems, we will use the above \rf{eq_Riccati--Volterra} to construct the closed-loop control as well.
To guarantee its well-posedness, we need the following standard assumption.

\ms

{\bf(A2)} There exists a constant $\l>0$ such that for all $t\in(0,T)$, $R_1(t)+R_2(t+\d){\bf1}_{[0,T-\d)}(t)\ges \l I$, $Q_i(t)\ges0$, $i=1,2,3$.

\ms

The following proposition gives its solvability.
\bt{solvability of Riccati}
Let {\rm (A1)--(A2)} hold. Then, the Riccati system \eqref{eq_Riccati--Volterra} admits a unique solution
$(P^{(1)},P^{(2)})\in\Pi(0,T)$ such that $\cR(\cd)\ges \b I$ for some constant $\b>0$.
\et

\br{Solvability of Riccati}

We discuss the solvability results of the above Riccati system.

We first look at the particular case of the SDDE \rf{state-distributed} and compare it with the study in \cite{Meng-Shi-Zhang-Zhao-2025}. In terms of their framework, both the cost functional and the diffusion term can depend on the distributed control delay which is out of our scope. However, to derive the Riccati system they have to assume that $C_2,b,\si, R_2, Q_2=0$. In addition, to obtain the solvability, they further require that $D_1=0$ and all the coefficients are time-invariant or continuous. It is worth mentioning that these assumptions are not needed here. Moreover, the methodologies developed in both papers are essentially different.

Next we turn to another particular case of the state equation \rf{state} with $A_2, B_2,C_2,Q_2,R_2=0$, and arrive at an LQ problem for a stochastic integro-differential system. Even though there are some positive results of the Riccati system (see \cite{Pandolfi-2018, Pandolfi-2025}) in deterministic scenario, the extension to the stochastic setting is still open. Here we fill this blank in a nice manner.
 
At last we point out two interesting facts even when
\rf{state} reduces to the deterministic system. Firstly, in contrast with the relevant literature (e.g. \cite{Alekal-Brunovsky-Chyung-Lee-1971, Delfour-Mitter-1972, Ichikawa-1982, Meng-Shi-Yong-2023}), the corresponding Riccati systems and their solvability appear for the first time since both the pointwise delay and the distributed delay are allowed to appear simultaneously in the cost functional. Secondly, when all the pointwise delayed terms disappear, the corresponding integro-differential system and cost functional can cover those in \cite{Pandolfi-2018, Pandolfi-2025}, and thus, the corresponding result in the above Theorem \ref{solvability of Riccati} is also new.

\er

Next we introduce the following backward system:
\vskip-6mm
\bel{EBSVIE}\left\{\ba{ll}
\ns\ds d\eta(t,s)=-\Big\{\cG_2(t,s)b(s) +{\Gamma}^*(t,s)^\top D_1(s)^\top\cG_1(s)\si(s) +{\Gamma}^*(t,s)^\top\int_s^T\Big[D_1(s)^\top \Upsilon(r,s)^\top\zeta(r,s)\\
\ns\ds\qq\qq+\Big(\int_s^r
\cB(\th,s)^\top \Pi(r,s,\th)^\top d\th\Big)\eta(r,s)\Big]dr\Big\}ds+\zeta(t,s)dW(s),\ 0<s<t<T,\\
\ns\ds \eta(t,t)=\cC(t)^\top\cG_1(t)\si(t) +{\Xi}^*(t)^\top D_1(t)^\top\cG_1(t)\si(t)+\int_t^T\Big\{\Big[\Big(\int_t^r \Pi(r,t,\th)\cB(\th,t) d\th\Big){\Xi}^*(t)\\
\ns\ds \qq\qq +\Upsilon(r,t)
\cA(t)\Big]^\top\eta(r,t)+\Big[\Upsilon(r,t)\cC(t) +\Upsilon(r,t)D_1(t){\Xi}^*(t)\Big]^\top\zeta(r,t)\Big\}dr, \ 0<t<T,
\ea\right.\ee
\vskip-2mm
\no where $\Xi^*(\cd)$, ${\Gamma}^*(\cd,\cd)$ are defined by
\vskip-2mm
\bel{check-Xi-decompose}\ba{ll}
\ns\ds \Xi^*(t)=-\cR(t)^{-1} D_1(t)^\top \cG_1(t)\cC(t),\ 0<t<T,\\
\ns\ds {\Gamma}^*(s,t)= -\cR(t)^{-1}\int_t^T\cB(\th,t)^\top \cG_3(s,t,\th)^\top d\th,\q 0<t<s<T,
\ea\ee
and $\Pi(\cd,\cd,\cd)$, $\Upsilon(\cd,\cd)$, $\cR(\cd)$, $\cG_1(\cd)$, $\cG_3(\cd,\cd,\cd)$ are defined as \rf{Pi}, \rf{Upsilon-cR}.
In terms of \cite{Hamaguchi-Wang-2024-II}, we name it the Type-II extended backward SVIE in our scenario. To study its well-posedness, we introduce the following space. Denote by
$L_{\mathbb{F},c}^2(\triangle_2\left(0, T\right) ; \mathbb{R}^{3n})$ the set of $\eta \in L_{\mathbb{F}}^2\left(\triangle_2\left(0, T\right) ; \mathbb{R}^{3n}\right)$ such that $s \mapsto \eta(t, s)$ is uniformly continuous on $\left(0, t\right)$ with the limits defined by $\eta(t, t)\equiv\lim _{s \uparrow t} \eta(t, s)$ and $\eta\left(t, 0\right)\equiv\lim _{s \downarrow 0} \eta(t, s)$ for a.e. $t \in\left(0, T\right)$, a.s., and $\eta(\cd,\cd)$ satisfies $\mathbb{E}(\int_{0}^T \sup\limits_{s \in\left[0, t\right]}|\eta(t, s)|^2 \mathrm{~d} t)^{\frac 1 2}<\i.$

\bt{solvability of BSVIE}
Let {\rm (A1)--(A2)} hold. Then, the Type-II extended backward SVIE \rf{EBSVIE} admits a unique solution $(\eta,\zeta)\in L^2_{\dbF,\mathrm{c}}(\triangle_2(0,T);\dbR^{3n})\times L^2_\dbF(\triangle_2(0,T);\dbR^{3n})$.
\et

\subsection{The closed-loop solvability of Problem (P)}

In this part we will give an explicit form of the optimal closed-loop strategy and some sufficient conditions for the closed-loop solvability of Problem (P).

Given $\Pi(\cd,\cd,\cd)$, $\Upsilon(\cd,\cd)$, $\cR(\cd)$, $\cG_1(\cd)$, $\cG_3(\cd,\cd,\cd)$ in \rf{Pi} and \rf{Upsilon-cR}, $P=(P^{(1)},P^{(2)})$ and $(\eta,\zeta)$ being the solutions to \eqref{eq_Riccati--Volterra} and \rf{EBSVIE}, respectively, we make the following conventions.
For $i=1,3$, $t_0\les t\les T$, we denote
\begin{eqnarray}
K_i^*(t)=\1n-\cR(t)^{-1}\1n\Big\{
{K}_i^{(1)}(t)+\3n\int_t^T\2n\int_t^T \cB(\th,t)^\top \cG_3(\a,t,\th)^\top{K}_i^{(2)}(\a,t)d\a d\th\Big\},
\ \label{Ki}
\end{eqnarray}
while for $i=2,4$ and $t_0\les s< t\les T$,
\begin{eqnarray}
K_i^*(t,s)=\1n-\cR(t)^{-1}\1n\Big\{
{K}_i^{(1)}(t,s)+\3n\int_t^T\3n\int_t^T\2n \cB(\th,t)^\top\cG_3(\a,t,\th)^\top\2n {K}_i^{(2)}(t,\a,s)d\a d\th\1n\Big\}.
\ \label{Ki-}
\end{eqnarray}
In addition, for $t_0\les t\les T$, let
\begin{eqnarray}
v^*(t)
\1n=\1n-\cR(t)^{-1}\1n\Big\{v^{(1)}(t) +\2n\int_t^T\int_t^T\1n\2n\cB(\th,t)^\top
\cG_3(\a,t,\th)^\top v^{(2)}(t,\a)d\a d\th\Big\}.
\label{v*}
\end{eqnarray}
In the above, each pair of $({K}_i^{(1)},{K}_i^{(2)})$ and $(v^{(1)},v^{(2)})$ have the following representations:
\bel{K4-co}\ba{ll}
\ns\ds{K}_1^{(1)}(t)=D_1(t)^\top\cG_1(t)C_1(t),\
{K}_1^{(2)}(\a,t)=\Upsilon(\a,t),\
{K}_3^{(1)}(t)=D_1(t)^\top\cG_1(t)C_2(t),\
{K}_3^{(2)}(\a,t)=0,\\
\ns\ds {K}_2^{(1)}(t,s)=D_1(t)^\top \cG_1(t)C_3(t)F(t,s)
+\Big(\int_t^{T}\cB(\th,t)^\top \cG_3(s+\d,t,\th)^\top  d\th\Big)(0,I,0)^\top {\bf1}_{[t-\d,T-\d]}(s),  \\
\ns\ds {K}_2^{(2)}(t,\a,s)=(0,0,I)^\top F(\a,s), \ \
{K}_4^{(1)}(t,s)=0,  \\
\ns\ds {K}_4^{(2)}(t,\a,s)=\Big(I,{\bf1}_{(s+2\d,\i)}(\a) {\bf1}_{(0,T-\d)}(t)I, \int_{s+\d}^\a F(\a,\th')^\top d\th'\Big)^\top {\bf1}_{(s+\d,T)}(\a)B_2(s+\d) {\bf1}_{[t-\d,T-\d]}(s) \\
\ns\ds\q+\Big(\ti\cF(t,s,\a),
\ti\cF(t,s,\a-\d),\int_t^\a\2n \ti\cF(t,s,\b)F(\a,\b)^\top d\b\Big)^\top,\ \ \ti\cF(t,s,\t)\equiv \int_t^\t\2n\ti F(\th',s)^\top B_3(\th')^\top d\th',
\ea\ee
and
\bel{v-co2}\ba{ll}
\ns\ds v^{(1)}(t)=D_1(t)^\top \cG_1(t)\sigma(t) \!+\!\!\int_t^TD_1(t)^\top \Upsilon(\a,t)^\top\zeta(\a,t)d\a\!+\!\!\int_t^T \int_\th^T\cB(\th,t)^\top  \Pi(\a,t,\th)^\top\eta(\a,t)d\a d\th, \\
\ns\ds v^{(2)}(t,\a)\!=\1n(0,I,0)^\top  \xi(\a-\d-t_0){\bf1}_{[t_0,t_0+\d]}(\a)\!
+\!\!\int_{t-\d}^{t_0}{\bf1}_{(\th'+\d,\i)}(\a)\Big(I,{\bf1}_{(\th'+2\d,\i)}(\a){\bf1}_{(0,T-\d)}(t)I, \\
\ns\ds \qq\qq\q \int_{\th'+\d}^\a F(\a,\b)^\top d\b\Big)^\top B_2(\th'+\d) \varsigma(\th'-t_0)d\th' {\bf1}_{ [t_0,t_0+\d]}(t).
\ea\ee
\no At this moment, we present the main result of the current section.

\bt{main theorem optimal control}
Let {\rm (A1)--(A2)} hold and $t_0\in[0,T)$ be given. Then, the five-tuple $(K_1^*(\cd),K_2^*(\cd,\cd), K_3^*(\cd)$, $K_4^*(\cd,\cd), v^*(\cd))$ given by \rf{Ki}--\rf{v-co2} is an optimal closed-loop strategy, and the following $u^*(\cd)$ is the unique optimal closed-loop outcome control of Problem {\rm (P)} on $[t_0,T]$:
\begin{eqnarray}
 u^*(t)\!=\!K_1^*(t)x^*(t)
\!+\!\int_{t_0}^tK_2^*(t,s)x^*(s)ds\!+\!K_3^*(t)x^*(t-\d)
 \!+\!\int_{t_0}^tK_4^*(t,s) u^*(s)ds\!+\!v^*(t),\ \ t\in[t_0,T].
\label{transformed optimal control}
\end{eqnarray}
\et

By the above theorem and Theorem 5.4 in \cite{Hamaguchi-Wang-2024-II}, we deduce the following result for Problem (P$_0$).
\bc{coro of main theorem optimal control}
Let {\rm (A1)--(A2)} hold and $t_0\in[0,T)$ be given. Then, the five-tuple $(K_1^*(\cd),K_2^*(\cd,\cd), K_3^*(\cd)$, $K_4^*(\cd,\cd), v^*(\cd))$ with  $v^{(1)}(\cd)=0$ is an optimal closed-loop strategy, and the above \rf{transformed optimal control} is the unique optimal closed-loop outcome control of Problem {\rm(P$_0$)}. In addition, the value function is given by
\begin{eqnarray}
V_0(t_0,\xi,\varsigma)=\int_{t_0}^T\blan P^{(1)}(t)\f(t),\f(t)\bran dt+\int_{t_0}^T\int_{t_0}^T\blan P^{(2)}(t_1,t_2,t_0)\f(t_2),\f(t_1)  \bran dt_1dt_2,
  \notag
\end{eqnarray}
for any $\left(t_0, \xi,\varsigma\right) \in [0,T)\times C([-\d,0];\dbR^n)\times L^2(-\d,0;\dbR^m)$.
\ec

\br{necessity}
We see that Theorem \ref{main theorem optimal control} gives a sufficient condition of the closed-loop solvability in terms of $K_i^*$ (i=1,2,3,4) which are explicitly and clearly constructed. At this moment, we are not sure about its necessity. However, Lemma \ref{closed-loop solvability-decompose} in the Appendix actually gives a new necessary condition in terms of the so-called causal feedback strategy developed in \cite{Hamaguchi-Wang-2024-II}. On the other hand, Theorem \ref{main theorem optimal control} is also true if the standard assumption {\rm (A2)} is relaxed properly (see \cite{Meng-Shi-Zhang-Zhao-2025}). For simplicity we prefer not to pursue these generalities.
\er

\br{dimension}
To prove Theorem \ref{main theorem optimal control}, we use the equivalence between the original problem {\rm (P)} and the new control problem associated with the state equation \rf{transformed state} and the cost functional \rf{transformed cost}. Even though the dimension of \rf{transformed state} and \rf{transformed cost} is higher than that of the original one, it is still finite dimensional and is essentially different from infinite-dimensional evolution control system method. For this new problem,we borrow some new matrix products notations exemplified by \rf{CPC}--\rf{DPC} from \cite{Hamaguchi-Wang-2024-II}.Then we introduce some computational techniques to convert them into the traditionalmatrix products, and derive the optimal closed-loop strategy as in Theorem \ref{main theorem optimal control}.
\er

\br{outcome control comparison}
There have been lots of works on closed-loop outcome controls of delayed control systems. However, they either contain only state delays \cite{Chen-Zhang-2023, Lee-You-1989, Liang-Xu-Zhang-2018, Lindquist-1973}, or only control delays  \cite{Chen-Wu-2011, Ichikawa-1982, Liang-Xu-Zhang-2018, Moon-2022, Wang-Zhang-2013, Zhang-Xu-2017}, or work in deterministic systems  \cite{Ichikawa-1982, Lee-You-1990}, or have time-invariant coefficients \cite{Ichikawa-1982, Lee-You-1990, Liang-Xu-Zhang-2018, Meng-Shi-Zhang-Zhao-2025, Moon-2022, Wang-Zhang-2013, Zhang-Xu-2017}, or have no delay in the cost functional \cite{Alekal-Brunovsky-Chyung-Lee-1971, Delfour-1986, Delfour-Mitter-1972, Ichikawa-1982,Xu-Shi-Zhang-2018}, or have no solvability of the associated Riccati systems \cite{Chen-Wu-2011, Ichikawa-1982, Lee-You-1990, Liang-Xu-Zhang-2018, Lindquist-1973, Moon-2022, Wang-Zhang-2013, Zhang-Xu-2017}. In this sense, Theorem \ref{main theorem optimal control} gives a unified treatment of the existing papers with distinctive methods. Moreover, the closed-loop outcome control is explicitly constructed without any continuity or even differentiability assumptions, and does not rely on the future state and avoids complex tools of conditional expectations \cite{Xu-Shi-Zhang-2018, Zhang-Xu-2017}.

\er

\section{Important cases}

In this section, we discuss five special yet important stochastic control systems and make relevant comparisons with the existing literature.

\subsection{Case I: Stochastic control systems with control delays only}

Consider the state equation
\bel{state ii}
\left\{\begin{aligned}
d x(t)= & \big[A_1(t) x(t)+B_1(t) u(t) +B_2(t) \nu(t)+B_3(t) \mu(t)\big] d t \\
& +\big[C_1(t) x(t)+D_1(t) u(t)\big] d W(t), \quad t \in[t_0, T], \\
x(t_0)= & \xi_0, u(t)=0, \quad t \in[t_0-\delta, t_0],
\end{aligned}\right.
\ee
along with the cost functional
$$\ba{ll}
\ns\ds J(t_0,\xi_0;u(\cdot))=   \dbE\int_{t_0}^T\big[x(t)^{\top} Q_1(t) x(t)+u(t)^{\top} R_1(t) u(t)\big] d t .
\ea
$$
Also for $\th,\a\in[-\d,0]$ and $t,\th',r\in[t_0,T]$ such that $t\les (\th'\wedge r)$, we define
\begin{eqnarray}
\mathcal{P}_1(t,\th',r)\equiv\big(I,0,0\big)\Big(\int_{\th'\vee r}^TP^{(1)}(s) ds+\int_r^T \int_{\th'}^TP^{(2)}(s,\a,t) d\a ds\Big)\big(I,0,0\big)^\top,\ \qq\qq\qq\qq\qq\q\! \label{P-control delay}\\
\cS_0(t)\equiv\mathcal{P}_1(t,t,t),\qq\qq\qq\qq\qq\qq\qq\qq\qq\qq\qq
\qq\qq\qq\qq\qq\qq\qq\ \  \label{cS0} \\
\cS_1(t,\th)\equiv B_2(t\1n+\1n\d\1n+\1n\th)^\top \mathcal{P}_1(t,t,t\1n+\1n\d\1n+\1n\th)+\3n\int_t^T\3n\ti F(\th',t\1n+\1n\th)^\top B_3(\th')^\top\mathcal{P}_1(t,\th',t)^\top d\th',\qq\qq\qq\qq\qq\ \label{cS1}\\
\cS_2(t,\th,\a)\equiv B_2(t\1n+\1n\d\1n+\1n\th)^\top\2n \Big[\1n\mathcal{P}_1 (t,t\1n+\1n\d\1n+\1n\th,t\1n+\1n\d\1n+\1n\a)^\top B_2(t\1n+\1n\d\1n+\1n\a)\1n+\3n\int_t^T\3n \mathcal{P}_1(t,\1n\th',t\1n+\1n\d\1n+\1n\th)B_3(\th')
\ti F(\th',t\1n+\1n\a)d\th'\Big]\notag\\
+\3n\int_t^T\3n\ti F(\th',t\1n+\1n\th)^\top
\1nB_3(\th')^\top\2n\Big[\1n\mathcal{P}_1(t,\th',t\1n+\1n\d\1n+\1n\a)^\top
\1nB_2(t\1n+\1n\d\1n+\1n\a)\1n+\3n\int_t^T\3n\mathcal{P}_1\1n (t,\1n\th',\1n\b)^\top\2n B_3(\b)\ti F(\b,t\1n+\1n\a)d\b\Big]d\th'.\q\label{cS2}
\end{eqnarray}
In this part, we will show that $\cS_0(\cd)$ and $\cS_1(\cd,\cd)$ satisfy
\bel{dcS0}\3n\3n\1n\left\{\ba{ll}
\ns\ds \frac{d}{dt}\cS_0(t)\1n+\1nA_1(t)^\top \cS_0(t)\1n+\1n\cS_0(t)A_1(t)\1n+\1nQ_{1}(t) +C_1(t)^\top\cS_0(t)C_1(t)-\2n\big[B_1(t)^\top\cS_0(t)\1n\\ [2mm]
\ns\ds\ +\cS_1(t,\1n0) \2n+D_1(t)^\top\cS_0(t)C_1(t)\big]\1n^\top \1n \cR(t)^{-1}\1n\big[B_1(t)^\top\cS_0(t)\1n +\2n\cS_1(t,\1n0)\2n+\1nD_1(t)^\top\cS_0(t)C_1(t)\big]\2n =\1n0,\ae\ t\in[t_0,T],\\ [1mm]
\ns\ds \cS_0(T)=0,\\ [2mm]
\ea\right.\ee
\vskip-4mm
\bel{dcS1}\3n\3n\2n\left\{\ba{ll}
\ns\ds \Big(\frac{\partial}{\partial t}-\frac{\partial}{ \partial\th}\Big)\cS_1(t,\th) +\ti F(t,t+\th)^\top B_3(t)^\top\cS_0(t)+\cS_1(t,\th)A_1(t) -\big[\cS_1(t,\th)B_1(t) \\ [2mm]
\ns\ds\q +\cS_2(t,\th,0)\big] \cR(t)^{-1}\big[B_1(t)^\top \cS_0(t)+\cS_1(t,0)+D_1(t)^\top\cS_0(t)C_1(t)\big]=0,\q \ae\ t\in[t_0,T],\th\in[-\d,0],\\ [2mm]
\ns\ds \cS_1(T,\th)=0,\ \cS_1(t,-\d)=B_2(t)^\top \cS_0(t)+\int_t^T\ti F(\th',t-\d)^\top B_3(\th')^\top\mathcal{P}_1(t,\th',t)^\top d\th'.\\
\ea\right.\ee
Moreover, $\cS_2(\cd,\cd,\cd)$ satisfies
\bel{dcS2}\3n\3n\left\{\ba{ll}
\ns\ds \Big(\frac{\partial}{\partial t}-\frac{\partial}{\partial\th} -\frac{\partial}{\partial\a}\Big)\cS_2(t,\th,\a)+\ti F(t,t+\th)^\top B_3(t)^\top \cS_1(t,\a)^\top+\cS_1(t,\th)B_3(t)^\top\ti F(t,t+\a)\\ [2mm]
\ns\ds\q -\1n\big[\cS_1\1n(t,\th)B_1\1n(t)\1n +\1n\cS_2(t,\th,0)\1n\big]\1n\cR\1n(t)^{-1}\1n \big[B_1(t)^\top\1n \cS_1(t,\a)^\top\1n+\1n\cS_2(t,0,\a)\1n \big]\1n=\1n0, \ae\ t\1n\in\1n[t_0,T],\th,\a\1n\in\1n[-\d,0],\\ [2mm]
\ns\ds \cS_2(t,\th,-\d)=\cS_1(t,\th)B_2(t)+B_2(t+\d+\th)^\top \int_t^T\mathcal{P}_1(t,\a,t+\d+\th)B_3(\a)\ti F(\a,t-\d)d\a\\
\ns\ds \q +\int_t^T\int_t^T\ti F(\a,t+\th)^\top B_3(\a)^\top\mathcal{P}_1(t,\a,\b)^\top  B_3(\b)
\ti F(\b,t-\d) d\b d\a , \\ [2mm]
\ns\ds \cS_2(t,-\d,\th)=\cS_2(t,\th,-\d)^\top,\ \cS_2(T,\th,\a)=0.
\ea\right.\ee
In addition,  the optimal closed-loop outcome control \rf{transformed optimal control}  is represented as follows:
\begin{eqnarray}
\1n{u}^*(t)=\1n-\cR(t)^{-1}\1n \bigg\{\big[B_1(t)^\top \cS_0(t)\1n+D_1(t)^\top\cS_0(t)C_1(t)+\1n\cS_1(t,0)\big] {x}^*(t)+\1n\int_{t\vee(t_0+\d)}^{(t+\d)\wedge T}\big[B_1(t)^\top \cS_1(t,r-\d-t)^\top \2n\notag\\
\3n\3n+\cS_2(t,0,r-\d-t)\big] {u}^*(r-\d)dr\bigg\}
+\cR(t)^{-1}\Big(\int_{t\vee(t_0+\d)}^{(t+\d)\wedge T}-\int_{t_0+\d}^{t+\d}\Big)\int_t^T\bigg\{ B_1(t)^\top\cP_1(t,t,\th')\qq\q\notag\\
+B_2(t+\d) ^\top\cP_1(t,t+\d,\th')^\top+\int_t^T\ti F(\th,t)^\top B_3(\th)^\top\cP_1(t,\th,\th')^\top d\th\bigg\}\notag\qq\qq\qq\qq\qq\qq\\
\times B_3(\th')\ti F(\th',r-\d) u^*(r-\d)d\th' dr,
\ \ae\ t\in[t_0,T],\qq\qq\qq\qq\qq\qq\qq\qq\qq\q\ \label{optimal-control-with-control-delay}
\end{eqnarray}
where $\cR(t)=R_1(t)+D_1(t)^\top\cS_0(t)D_1(t)$.
\ss

We state the main result of this subsection as follows.
\bc{optimal control with control delay}
Let {\rm (A1)-(A2)} hold and $A_2,A_3,C_2,C_3,Q_2,Q_3,R_2,b,\si=0$. Then, $\cS_0(\cd)$, $\cS_1(\cd,\cd)$ and $\cS_2(\cd,\cd,\cd)$, defined by \rf{cS0}--\rf{cS2}, satisfy the coupled Riccati equations \rf{dcS0}--\rf{dcS2}, and the process in \rf{optimal-control-with-control-delay} is the optimal closed-loop outcome control.
\ec

\br{remark for case ii}
When Problem {\rm (P)} contains only control delays, we obtain the optimal closed-loop outcome control  \rf{optimal-control-with-control-delay} by the coupled Riccati equations \rf{dcS0}--\rf{dcS2}. Furthermore, if the diffusion term disappears in \rf{state ii}, then \rf{dcS0}--\rf{dcS2} essentially reduce to (2.33)--(2.38) in \cite{Ichikawa-1982}. Compared with  \cite{Chen-Wu-2011, Ichikawa-1982, Liang-Xu-Zhang-2018, Moon-2022, Wang-Zhang-2013, Zhang-Xu-2017}, we successfully obtain the solvability of the desired Riccati system.
\er

\subsection{Case II: Stochastic control systems with state delays only}

Consider the state equation
 $$
\left\{\ba{ll}
\ns\ds  dx(t)=\big[A_1x(t)+A_2 y(t)+B_1u(t)\big]dt+\big[C_1x(t)+C_2 y(t)+D_1u(t)\big]dW(t),\q t\in[t_0,T],\\
\ns\ds x(t)=\xi(t-t_0),\ \ t\in[t_0-\d,t_0],
\ea\right.
$$
along with the cost functional
$$
\begin{aligned}
J(t_0,\xi(\cd);u(\cdot))= & \dbE\int_{t_0}^T\left[x(t)^{\top} Q_1 x(t)\right. \left.+u(t)^{\top} R_1 u(t)\right] d t .
\end{aligned}
$$
Here the coefficients are time-invariant. For $t_0\les t\les\th\les T$, we define
\begin{eqnarray}
\cP_2(t)\equiv\int_t^T\big(I,{\bf1}_{(\d,\i)}(s-t)I,0\big)P^{(1)}(s) \big(I,{\bf1}_{(\d,\i)}(s-t)I,0\big)^\top ds\qq\qq\qq\qq\qq\qq\qq \notag\\
+\int_t^T\int_t^T\big(I,{\bf1}_{(\d,\i)}(s_1-t)I,0\big)P^{(2)}(s_1,s_2,t) \big(I,{\bf1}_{(\d,\i)}(s_2-t)I,0\big)ds_1ds_2,\qq\qq\qq
\label{cP2}\\
{\cP}_3(t,\th)\equiv\Big[\big(I,{\bf1}_{(\d,\i)}(\th-t)I,0\big)P^{(1)} (\th)^\top+\int_t^T\big(I,{\bf1}_{(\d,\i)}(r-t)I,0\big)P^{(2)} (\th,r,t)^\top dr\Big](0,I,0)^\top.\q
\label{cP3}
\end{eqnarray}
\no In this part, we will show that $\cP_2(\cd)$ and $\cP_3(\cd,\cd)$ satisfy the following coupled Riccati equations. More precisely, for $t \in(T-\d, T], \theta \in(t, T]$, we have
\bel{Riccati-1-state delay}\5n\left\{\ba{ll}
\ns\ds -\dot{\cP_2}(t)= \cP_2(t) A_1+A_1^\top \cP_2(t)+C_1^\top \cP_2(t) C_1+Q_1 \\
\ns\ds \qq-\big(B_1^\top \cP_2(t)+D_1^\top \cP_2(t) C_1\big)^\top \big(R_1+D_1^\top\cP_2(t)D_1\big)^{-1}\big(B_1^\top \cP_2(t)+D_1^\top \cP_2(t) C_1\big), \\
\ns\ds -\frac{\partial \cP_3(t, \theta)}{\partial t}=  A_1^{\top} \cP_3(t, \theta)-\big(B_1^\top \cP_2(t)+D_1^\top \cP_2(t) C_1\big)^\top \big(R_1+D_1^\top\cP_2(t)D_1\big)^{-1} B_1^\top \cP_3(t, \theta), \\
\ns\ds \cP_3(t, t)=  \cP_2(t) A_2+C_1^\top\cP_2(t) C_2-\big(B_1^\top \cP_2(t)+D_1^\top \cP_2(t) C_1\big)^\top \big(R_1+D_1^\top\cP_2(t)D_1\big)^{-1} D_1^\top \cP_2(t) C_2,\\
\ns\ds \cP_2(T)=0,
\ea\right.\ee
while for $t \in[0,T-\d], \theta \in(t, t+\d]$,
\bel{Riccati-2-state delay}\left\{\ba{ll}
\ns\ds -\dot{\cP_2}(t)= \cP_2(t) A_1+A_1^\top \cP_2(t)+C_1^\top \cP_2(t) C_1+C_2^\top \cP_2(t+\d) C_2+Q_1 \\
\ns\ds \qq+\cP_3(t, t+\d)+\cP_3(t, t+\d)^\top-\big(B_1^\top \cP_2(t)+D_1^\top \cP_2(t) C_1\big)^\top \big(R_1+D_1^\top\cP_2(t)D_1\big)^{-1}\big(B_1^\top \cP_2(t) \\
\ns\ds  \qq+D_1^\top \cP_2(t) C_1\big)-\big(D_1^\top \cP_2(t+\d) C_2\big)^\top \big(R_1+D_1^\top\cP_2(t)D_1\big)^{-1} D_1^\top \cP_2(t+\d) C_2, \\
\ns\ds -\frac{\partial \cP_3(t, \theta)}{\partial t}=    A_1^\top \cP_3(t, \theta)-\big(B_1^\top \cP_2(t)+D_1^\top \cP_2(t) C_1\big)^\top \big(R_1+D_1^\top\cP_2(t)D_1\big)^{-1} B_1^\top \cP_3(t, \theta)\\
\ns\ds \qq+\cP_3(\th, t+\d)^\top A_2-\big(B_1^\top \cP_3(\th, t+\d) \big)^\top  \big(R_1+D_1^\top\cP_2(t)D_1\big)^{-1} D_1 \cP_2(t) C_2 \\
\ns\ds \qq-\int_t^\theta \big(B_1^\top \cP_3(s, t+\d)\big)^\top \big(R_1+D_1^\top\cP_2(t)D_1\big)^{-1} B_1^\top \cP_3(s, \th) d s, \\
\ns\ds \cP_3(t, t)=  \cP_2(t) A_2\1n+\1nC_1^\top\cP_2(t) C_2\1n-\1n\big(B_1^\top \cP_2(t)\1n+\1nD_1^\top \cP_2(t) C_1\big)^\top \big(R_1\1n+\1nD_1^\top\cP_2(t)D_1\big)^{-1} D_1^\top \cP_2(t) C_2.\\
\ea\right.\ee
We state the main result of this subsection as follows.
\bc{optimal control with state delay}
Let {\rm (A1)--(A2)} hold with $A_3,B_2,B_3,C_3,Q_2,Q_3,R_2,b,\si=0$.
Then, $\cP_2(\cd)$ and $\cP_3(\cd,\cd)$, defined by \rf{cP2}--\rf{cP3}, satisfy the coupled Riccati equation \rf{Riccati-1-state delay}--\rf{Riccati-2-state delay}. In this case, the following process is the optimal closed-loop outcome control
\begin{eqnarray}
u^*(t)= -\big(R_1+D_1^\top\cP_2(t)D_1\big)^{-1}\Big\{ \big(B_1^\top \cP_2(t)+D_1^\top \cP_2(t) C_1\big) x^*(t)
+ \int_{t\vee{(t_0+\d)}}^{ (t+\d)\wedge T} B_1^\top \cP_3(t, s) x^*(s-\d) d s \notag\\
+ D_1^\top \cP_2(t) C_2 x^*(t-\d)+\int_t^{T\wedge(t_0+\d)}B_1^\top \cP_3(t,s)\xi(s-t_0-\d)ds\Big\}, \ t\in[t_0,T].\qq\qq\qq\q
 \label{u*-state delay}
\end{eqnarray}
\ec

\br{remark iii}
When Problem {\rm (P)} contains only state delays, the coupled Riccati equations \rf{Riccati-1-state delay}--\rf{Riccati-2-state delay} are the same as (3)--(12) in \cite{Liang-Xu-Zhang-2018}, the optimal closed-loop outcome control \rf{u*-state delay} coincides with (13) in \cite{Liang-Xu-Zhang-2018}, and Corollary \ref{optimal control with state delay} is similar to Theorem 1 in \cite{Liang-Xu-Zhang-2018}. Compared with \cite{Lee-You-1990, Liang-Xu-Zhang-2018, Lindquist-1973}, one advantage lies in the positive result on the solvability of the Riccati system. References \cite{Delfour-1977-1,Delfour-1986} focused on deterministic systems involving both pointwise delay and distributed delay of state. They not only derived the corresponding optimal control results but also established the solvability of the Riccati equations. In contrast, the solvability conclusion proposed in this paper only requires the coefficients to be integrable, without the need for differentiability, thus relaxing the constraints imposed by \cite{Delfour-1977-1,Delfour-1986}.
\er

\subsection{Case III: Stochastic control systems with pointwise delays only}
Consider the state equation
 $$
\left\{\ba{ll}
\ns\ds  dx(t)=\big[A_1(t)x(t)+A_2(t) y(t)+B_1(t)u(t)+B_2(t)\nu(t)\big]dt\\
\ns\ds\qq\qq +\big[C_1(t)x(t)+C_2(t) y(t)+D_1(t)u(t)\big]dW(t),\q t\in(t_0,T),\\
\ns\ds x(t)=\xi(t-t_0),\ u(t)=\varsigma(t-t_0),\ t\in[t_0-\d,t_0],
\ea\right.
$$
along with the cost functional
$$
\begin{aligned}
J(t_0,\xi(\cd),\varsigma(\cd);u(\cdot))= & \dbE\int_{t_0}^T\left[x(t)^{\top} Q_1(t) x(t)\right. \left.+y(t)^{\top} Q_2(t) y(t)+u(t)^{\top} R_1(t) u(t)+\nu(t)^\top R_2(t)\nu(t)\right] d t .
\end{aligned}
$$
\no For $0<t<(s\wedge\bar s)<T$, we define $\cG_1(\cd)$, $\cG_2(\cd,\cd)$ as in \rf{Upsilon-cR}, and
\vskip-3mm
\bel{Upsilon-cR-TI}\ba{ll}
\ns\ds  Q(t)=\begin{bmatrix}Q_1(t) & 0\\ 0 & Q_2(t)\end{bmatrix},\q
\cC(t)
=\big(
C_1(t),
C_2(t)
\big), \ \   \cA(t)
=\big(
A_1(t),
A_2(t)
\big),\ \
\cB(t)
=\begin{bmatrix}
B_1(t)\\
B_2(t+\d)
\end{bmatrix},\\
\ns\ds
\Upsilon( s,t)
=\begin{bmatrix}
I\\
{\bf1}_{(\d,\i)}( s-t)I
\end{bmatrix},\q
\ \cR(t)=R_1(t)+R_2(t+\d){\bf1}_{[0,T-\d)}(t)+D_1(t)^\top \cG_1(t)D_1(t),\\
%
\ns\ds \cG_3( s,t)= P^{(1)}( s)\Big(\Upsilon( s,t), {\bf1}_{(\d,\i)} ( s-t) \Upsilon(s,t+\d)\Big)\\
\ns\ds\qq\qq\ +\int_t^TP^{(2)}( s,r,t)\Big(\Upsilon(r,t),{\bf1}_{(\d,\i)} (r-t) \Upsilon(r,t+\d)\Big) dr.
\ea\ee
\vskip-1mm
\no We consider the Riccati system \rf{eq_Riccati--Volterra}, where $P^{(2)}(\cd,\cd,\cd)$ satisfies
\vskip-5mm
\bel{eq_Riccati--Volterra-TI}\5n\2n\left\{\2n\ba{ll}
\ns\ds P^{(2)}(s,t,r)\!=\!P^{(2)}(s,t,t\wedge s)
\!-\!\!\int_r^{t\wedge s}\! \!\cG_3(s,\t)\cB(\t) \cR(\t)^{-1}
\cB(\t)^\top \cG_3(t,\t)^\top d\t,\ 0<r<(s\wedge t)<T,\\
\ns\ds
	P^{(2)}(\bar s,t,t)=P^{(2)}(t,\bar s,t)^\top =\cG_2(\bar s,t)\cA(t)- \cG_3(\bar s,t)\cB (t)\cR(t)^{-1} D_1(t)^\top\cG_1(t)\cC(t),\ 0<t<\bar s<T.
\ea\right.\ee
\vskip-1mm
\no With notations in \rf{Upsilon-cR-TI}, we consider the closed-loop strategy:
\vskip-3mm
\bel{v-co2-TI}\ba{ll}
\ns\ds K_i^*(t)=\1n-\cR(t)^{-1}\1n\Big\{
{K}_i^{(1)}(t)+\cB(t)^\top \!\! \int_t^T\cG_3(\a,t)^\top {K}_i^{(2)}(\a,t)d\a  \Big\}, \ i=1,3, \ t_0\les t\les T,\\
\ns\ds K_i^*(t,s)=\1n-\cR(t)^{-1}\1n\Big\{
{K}_i^{(1)}(t,s)\2n+\2n\cB(t)^\top\!\!
 \int_t^T \cG_3(\a,t)^\top {K}_i^{(2)}(t,\a,s)d\a\Big\}, \ i=2,4, \ \  t_0\les s< t\les T,\\
\ns\ds v^*(t)
\1n=\1n-\cR(t)^{-1}\1n\Big\{v^{(1)}(t)+\cB(t)^\top \int_{t}^T \cG_3(\a,t)^\top
 v^{(2)}(t,\a)d\a \Big\},\ t_0\les t\les T,
\ea
\ee
where ${K}_1^{(1)}, {K}_1^{(2)},
{K}_3^{(1)},
{K}_3^{(2)}, {K}_4^{(1)}$ have the same forms as in \rf{K4-co}, and
$$\ba{ll}
%
\ns\ds {K}_2^{(1)}(t,s)=
\cB(t)^\top \cG_3(s+\d,t)^\top (0,I)^\top {\bf1}_{[t-\d,T-\d]}(s),\qq\qq {K}_2^{(2)}\1n(t,\a,s)=0,\\
\ns\ds {K}_4^{(2)}\1n(t,\a,s)=\big(I\1n,{\bf1}_{(s+2\d,\i)}(\a) {\bf1}_{(0,T-\d)}(t)I\big)^\top\1n {\bf1}_{(s+\d,T)}(\a)B_2 (s+\d) {\bf1}_{[t-\d,T-\d]}\1n(s),\\
\ns\ds v^{(1)}(t)=0,\q
v^{(2)}(t,\a)=\1n(0,I)^\top  \xi(\a-\d-t_0){\bf1}_{[t_0,t_0+\d]}(\a)
+\int_{t-\d}^{t_0}{\bf1}_{(\th'+\d,\i)}(\a)\notag\\
\times \big(I,{\bf1}_{(\th'+2\d,\i)}(\a){\bf1}_{(0,T-\d)}(t)I\big)^\top B_2(\th'+\d) \varsigma(\th'-t_0)d\th' {\bf1}_{ [t_0,t_0+\d]}(t).
\ea$$
\no As a result of Theorem \ref{main theorem optimal control}, we derive the following result.

\bc{optimal control with time-invariant coefficients}
Let {\rm (A1)--(A2)} hold with $A_3,B_3,C_3,Q_3,b,\si=0$. Then, all the strategies in \rf{v-co2-TI} and the process $u^*(\cd)$ in the same form of \rf{transformed optimal control} are optimal.
\ec

\br{remark-TI}
We make some comparisons with the existing literature.
Firstly, in contrast with \cite{Liang-Xu-Zhang-2018, Meng-Shi-Zhang-Zhao-2025}, our advantages lie in that our coefficients are allowed to be time-variant and the cost functional depends on both the pointwise state delay and pointwise control delay. Secondly, even when these features disappear, our framework is still general than that in \cite{Liang-Xu-Zhang-2018} (except the disappearance of pointwise control delay in the diffusion term), and the solvability of the Riccati system is given. In addition, \cite{Meng-Shi-Zhang-Zhao-2025} requires the diffusion term to be independent of the pointwise state delay and control variables, while we drop this assumption here.
Thirdly, even though the coefficients in \cite{Chen-Zhang-2023} are time-variant, both their state equations and cost functional are still particular cases of ours. Eventually, even for deterministic systems, our coefficients are allowed to be integrable, but not necessarily differentiable as in \cite{Delfour-1977-1}.
\er

\subsection{Case IV: Stochastic control systems with distributed delays only}

Consider the state equation
 $$
\left\{\ba{ll}
\ns\ds  dx(t)=[A_1(t)x(t)+A_3(t) z(t)+B_1(t)u(t)+B_3(t)\mu(t)]dt\\
\ns\ds\qq\qq +[C_1(t)x(t)+C_3(t) z(t)+D_1(t)u(t)]dW(t),\q t\in(t_0,T),\\
\ns\ds x(t)=\xi(t-t_0),\ u(t)=\varsigma(t-t_0),\ t\in[t_0-\d,t_0],
\ea\right.
$$
along with the cost functional
$$
\begin{aligned}
J(t_0,\xi(\cd),\varsigma(\cd);u(\cdot))= & \dbE\int_{t_0}^T\left[x(t)^{\top} Q_1(t) x(t)\right. \left.+z(t)^{\top} Q_3(t) z(t)+u(t)^{\top} R_1(t) u(t)\right] d t .
\end{aligned}
$$
In this case, we consider the Riccati system \rf{eq_Riccati--Volterra},
\no where for $0<t<(s\wedge\bar s\wedge \th)<T$,
\vskip-3mm
\bel{Upsilon-cR-TI-}\ba{ll}
\ns\ds  Q(t)=\begin{bmatrix}Q_1(t) & 0\\ 0 & Q_3(t)\end{bmatrix},\q
\cC(t)
=\big(
C_1(t),
C_3(t)
\big), \ \   \cA(t)
=\big(
A_1(t),
A_3(t)
\big),\ \
\cB(\th,t)
=\begin{bmatrix}
B_1(t)\\
B_3(\th)\ti F(\th,t)
\end{bmatrix},\\
\ns\ds
\Upsilon( s,t)
=\begin{bmatrix}
I\\
\cE(s,t)
\end{bmatrix},\q
\ \cR(t)=R_1(t)+D_1(t)^\top \cG_1(t)D_1(t), \  \
\Pi(s,t,\th)=\begin{bmatrix}\frac{1}{s-t} I & I\\ \frac{1}{s-t}\cE(s,t) & \cE(s,\th)\end{bmatrix},

 \\
\ea\ee
\vskip-1mm
\no and $\cG_1(\cd)$, $\cG_2(\cd,\cd)$, $\cG_3(\cd,\cd,\cd)$ are defined in \rf{Upsilon-cR}.
With the notations in \rf{Upsilon-cR-TI-}, we consider the closed-loop strategy \rf{Ki}--\rf{v-co2},
where ${K}_1^{(1)}, {K}_1^{(2)}, {K}_3^{(2)},{K}_4^{(1)}$ have the same forms as in \rf{K4-co}, $v^{(1)}, \ v^{(2)}=0$, and
\bel{K4-co-2}\ba{ll}
\ns\ds
{K}_3^{(1)}(t)=0, \ \ {K}_2^{(1)}(t,s)=D_1(t)^\top \cG_1(t)C_3(t)F(t,s),  \  \ {K}_2^{(2)}(t,\a,s)=(0,I)^\top F(\a,s), \ \
  \\
\ns\ds {K}_4^{(2)}(t,\a,s)= \Big(\ti\cF(t,s,\a),\int_t^\a\2n \ti\cF(t,s,\b)F(\a,\b)^\top d\b\Big)^\top,\ \ \ti\cF(t,s,\t)\equiv \int_t^\t\2n\ti F(\th',s)^\top B_3(\th')^\top d\th'.
\ea\ee
\no Then, as a result of Theorem \ref{main theorem optimal control}, we derive the following result.

\bc{optimal control with time-invariant coefficients-II}
Let {\rm (A1)--(A2)} hold and $A_2,B_2,C_2,Q_2,R_2, b,\si=0$.
Then, all the closed-loop strategies \rf{Ki}--\rf{Ki-}  with
\rf{K4-co-2}  and the following process are optimal:
\bel{SVIDE-feedback}\ba{ll}
\ns\ds
 u^*(t)=K_1^*(t)x^*(t)
+\int_{t_0}^tK_2^*(t,s)x^*(s)ds
 +\int_{t_0}^tK_4^*(t,s) u^*(s)ds,\ t\in[t_0,T].
\ea\ee
\ec
\br{remark-svide}
LQ problems for deterministic integro-differential equations were treated in \cite{Pandolfi-2018}, see also
\cite[Subsection 5.5]{Gong-Wang-2025}. Both of them derived the Riccati systems and the closed-loop optimal controls in the spirit of \rf{SVIDE-feedback}.
However, they require that $A_1,B_3\equiv0$, $A_3\equiv 1$, $F(\cd,\cd)$ is the convolution kernel, $Q_3\equiv0$, $Q_1$, $R_1$ are time-invariant.

As to the case of stochastic integro-differential equations (SVIDEs), let us point out the discussion in \cite[Subsection 5.3]{Gong-Wang-2025} where the SVIDE is directly regarded as a special SVIE. However, there are no closed-loop controls in the spirit of \rf{SVIDE-feedback} and the corresponding Riccati systems. To our best knowledge, the above Corollary \ref{optimal control with time-invariant coefficients-II} appears for the first time.
\er

\subsection{Case V: Stochastic control systems without delay}

In this subsection, let us look at the particular SDEs case. To this end, we consider
\bel{Riccati-without delay}\left\{\ba{ll}
\ns\ds \dot{\cP}(t)+\cP(t) A_1(t)+A_1(t)^{\top} \cP(t)+C_1(t)^{\top} \cP(t) C_1(t)+Q_1(t)-\left(\cP(t) B_1(t)+C_1(t)^{\top} \cP(t) D_1(t)\right) \\
\q\times\cR(t)^{-1}\left(B_1(t)^{\top} \cP(t)+D_1(t)^{\top} \cP(t) C_1(t)\right)=0, \quad 0<t<T,\\
\ns\ds \cP(T)=0,
\ea\right.\ee
and the following backward stochastic differential equation:
\bel{BSDE-without delay}
\left\{\begin{aligned}
d \widetilde{\eta}(t)= & -\left\{\left[A_1^{\top}(t)-\left(\mathcal{P}(t) B_1(t)+C_1^{\top}(t) \mathcal{P}(t) D_1(t)\right)\cR(t)^{-1} B_1^{\top}(t)\right] \widetilde{\eta}(t)\right. \\
& +\left[C_1^{\top}(t)-\left(\mathcal{P}(t) B_1(t)+C_1^{\top}(t) \mathcal{P}(t) D_1(t)\right)\cR(t)^{-1} D_1^{\top}(t)\right] \tilde{\zeta}(t) \\
& +\left[C_1^{\top}(t)-\left(\mathcal{P}(t) B_1(t)+C_1^{\top}(t) \mathcal{P}(t) D_1(t)\right)\cR(t)^{-1} D_1^{\top}(t)\right] \mathcal{P}(t) \sigma(t) \\
& +\mathcal{P}(t) b(t)\big\} d t+\tilde{\zeta}(t) d W(t), \quad 0<t<T, \\
\widetilde{\eta}(T)= & 0,
\end{aligned}\right.\ee
where $\cR(\cd)= R_1(\cd)+D_1(\cd)^\top\cP(\cd)D_1(\cd)$.

On the other hand, given $P=(P^{(1)},P^{(2)})$, $(\eta,\zeta)$ satisfying the Riccati system \rf{eq_Riccati--Volterra} and the backward SVIE \rf{EBSVIE}, respectively, we define
\bel{Definition-cp-et-al}\ba{ll}
\ns\ds \cP(t)=(I,0,0)\Big(\int_t^TP^{(1)}(s) ds+\int_t^T\int_t^TP^{(2)}
(s_1,s_2,t) ds_1ds_2\Big)(I,0,0)^\top,\ 0< t< T, \\
\ns\ds \widetilde{\eta}(t)=\int_t^T (I,0,0)\eta(s,t)ds,\qq \tilde{\zeta}(t)=\int_t^T (I,0,0)\zeta(s,t)ds,\ 0< t< T.
\ea\ee
Then, we have the following result.
\bc{closed-loop solvability without delay}
Let {\rm (A1)-(A2)} hold such that $A_2,A_3,B_2,B_3,C_2,C_3,Q_2,Q_3,R_2=0$. Then, the above $\cP(\cd)$ and $(\wt\eta(\cd),\wt\zeta(\cd))$ are the unique solutions to \rf{Riccati-without delay} and \rf{BSDE-without delay}, respectively. In addition, the five-tuple $(K_1^*(\cd), 0,0,0,v^*(\cd))$ is the optimal closed-loop strategy of Problem {\rm (P)}, where
$$\ba{ll}
\ns\ds K_1^*(t)=-\cR(t)^{-1}\left( B_1(t)^\top\mathcal{P}(t)+D_1(t)^\top \mathcal{P}(t)C_1(t) \right),\  t\in[t_0,T],\\
\ns\ds v^*(t)=-\cR(t)^{-1}\Big(B_1^{\top}(t) \widetilde{\eta}(t)+D_1^{\top}(t) \tilde{\zeta}(t)+D_1^{\top}(t) \mathcal{P}(t) \sigma(t)\Big),\  t\in[t_0,T].
\ea$$
\ec

\br{rela of cl solv}
When delays appear in Problem {\rm (P)}, the optimal closed-loop strategy $(K_1^*(\cd),v^*(\cd))$, the equation \rf{Riccati-without delay} and the equation \rf{BSDE-without delay} are reduced to that in \cite{Sun-Yong-2020}.
\er

\section{Concluding remarks}
In this paper we study a general stochastic LQ optimal control problem, where the coefficients are time-varying, and both state delay and control delay can appear in the state equation and the cost functional. We put the original state process and its delay processes together for dimension expansion, and use the Volterra integral system without delay to describe the new process, then transform the original delayed problem into the control problem without delay. Based on the equivalent problem, we propose the closed-loop solvability of the delayed problem, and assure it by the solvability of the Riccati system and the extended backward SVIEs. Furthermore, we derive the optimal closed-loop outcome control and obtain the solvability of the associated Riccati system. Finally, we study several important stochastic systems and find that our results are consistent with those in the existing literature.

\bibliographystyle{plain}
\bibliography{references}

\section*{Appendix}
\renewcommand{\theequation}{A.\arabic{equation}}
\setcounter{theorem}{0}
\renewcommand{\thetheorem}{A.\arabic{theorem}}

{\bf The proof of Proposition \ref{Solvability of SDDE}:}
\bpf
Let $\b>0$, $T_0\in(t_0,T]$ be undetermined constants, and
$$
\cM_\b[t_0,T_0]\equiv L_\dbF^2(\O;C([t_0,T_0];\dbR^n)),
$$
equipped with the norm
$$||x(\cd)||_{\cM_{\b}[t_0,T_0]}\equiv\Big(\dbE\big[\sup\limits_{t_0\les t\les T_0}e^{\b h(t)}|x(t)|^2\big]\Big)^{\frac 1 2}.$$
Here
$$h(t)\equiv \int_{t_0}^t\Big(|\bar A_1(s)|^2+|\bar A_2(s)|^2+|\bar A_3(s)|^2\Big)ds,\ \ t\in[t_0,T_0].$$
For any $\bar x(\cd)\in {\cM_{\b}[t_0,T_0]}$ with $\bar x(t)=\bar \xi(t-t_0)$, $t\in[t_0-\d,t_0]$,
consider the mapping
$$\ba{ll}
\cT:\cM_{\b}[t_0,T_0]\rightarrow\cM_\b[t_0,T_0],\\
\qq \qq\ \bar x(\cd)\mapsto \bar X(\cd),
\ea$$
where $\bar X(\cd)$ is the solution to the following SDDE:
\bel{SDDE-proof}\left\{\ba{ll}
\ns\ds  d\bar X(t)=[\bar A_1(t)\bar x(t)+\bar A_2(t) \bar y(t)+\bar A_3(t)\bar z(t)+\bar b(t)]dt\\
\ns\ds\qq\qq+[\bar C_1(t)\bar x(t)+\bar C_2(t)\bar y(t)+\bar C_3(t)\bar z(t)+\bar \si(t)]dW(t),\q t\in(t_0,T),\\
\ns\ds \bar X(t)=\bar \xi(t-t_0),\ \ t\in[t_0-\d,t_0].
\ea\right.
\ee
For any $\bar x_1(\cd), \bar x_2(\cd)\in {\cM_{\b}[t_0,T_0]}$ and $\bar x_1(t),\bar x_2(t)=\bar \xi(t-t_0)$, $t\in[t_0-\d,t_0]$, denote
$$
\bar X_1(\cd)=\cT(\bar x_1(\cd)), \bar X_2(\cd)=\cT(\bar x_2(\cd)),\\
\hat X(\cd)=\bar X_1(\cd)-\bar X_2(\cd), \hat x(\cd)=\bar x_1(\cd)-\bar x_2(\cd).
$$
Then, applying It\^o formula to $s\mapsto e^{-\b h(s)}|\hat X(s)|^2$ on $[t_0,t]$, we have
\vskip-2mm
$$\ba{ll}
\ns\ds  e^{-\b h(t)}|\hat X(t)|^2
 \les \int_{t_0}^te^{-\b h(s)}\Big[(1-\b)|\hat X(s)|^2\big(|\bar A_1(s)|^2+|\bar A_2(s)|^2+|\bar A_3(s)|^2\big)\\
\ns\ds \qq+\big(2|\bar C_1(s)|^2+1\big)|\hat x(s)|^2+
\big(2|\bar C_2(s)|^2+1\big)|\hat y(s)|^2+\big(2|\bar C_3(s)|^2+1\big)|\hat z(s)|^2\Big]ds\\
\ns\ds \qq+\int_{t_0}^t2\blan e^{-\b h(s)}\hat X(s),
\bar C_1(s)\hat x(s)+\bar C_2(s)\hat y(s)+\bar C_3(s)\hat z(s)\bran dW(s).
\ea$$
Notice that
$$\ba{ll}
\ns\ds \dbE\int_{t_0}^{T_0}e^{-\b h(s)}\big(2 |\bar C_3(s)|^2+1\big)|\hat z(s)|^2ds\\
\ns\ds \les \dbE\Big[\sup\limits_{t_0\les s\les T_0}\Big\{e^{-\b h(s)}|\hat x(s)|^2\Big(2|\bar C_3(s)|^2+1\Big)\Big\}\int_{t_0}^{T_0} \Big(\int_{t_0}^s|\bar F(s,r)|dr\Big)^2 ds\Big],
\ea$$
and
$$\ba{ll}
\ns\ds \dbE\Big[\sup\limits_{t_0\les t\les T_0}\Big|2\int_{t_0}^t\blan e^{-\b h(s)}\hat X(s),
\bar C_1(s)\hat x(s)+\bar C_2(s)\hat y(s)+\bar C_3(s)\hat z(s)\bran dW(s)\Big|\Big]\\
\ns\ds \les\dbE\Big(\int_{t_0}^{T_0}4e^{-2\b h(s)}|\hat X(s)|^2\big|\bar C_1(s)\hat x(s)+\bar C_2(s)\hat y(s)+\bar C_3(s)\hat z(s)\big|^2ds\Big)^{\frac 1 2}\\
\ns\ds \les\frac 1 2\dbE\sup\limits_{t_0\les t\les T_0}\Big\{e^{-\b h(s)}|\hat X(s)|^2\Big\}+\Big(
4(T_0-t_0)\sup\limits_{0\les s\les T}\Big(|\bar C_1(s)|^2+|\bar C_2(s)|^2\Big)\\
\ns\ds \q+4\sup\limits_{0\les s\les T}|\bar C_3(s)|^2\int_{t_0}^{T_0}\Big(\int_{t_0}^s|\bar F(s,r)|dr \Big)^2ds\Big)\dbE\sup\limits_{t_0\les s\les T_0}\Big\{e^{-\b h(s)}|\hat x(s)|^2\Big\}.
\ea$$
Let $\b>1$. Then, we deduce
\vskip-2mm
$$\ba{ll}
\ns\ds \dbE\sup\limits_{t_0\les t\les T_0}\big\{e^{-\b h(t)}|\hat X(t)|^2\big\}\les 2\Big[6(T_0-t_0)\sup\limits_{t_0\les s\les T_0}\Big\{1+|\bar C_1(s)|^2+|\bar C_2(s)|^2\Big\}\\
\ns\ds \qq+3\sup\limits_{t_0\les s\les T_0}\Big(2|\bar C_3(s)|^2+1\Big)\int_{t_0}^{T_0}\Big(\int_{t_0}^s|\bar F(s,r)| dr\Big)^2ds\Big]\dbE\sup\limits_{t_0\les s\les T_0}\big[e^{-\b h(s)}|\hat x(s)|^2\big].
\ea$$
Since $\bar C_1(\cd), \bar C_2(\cd),\bar C_3(\cd)\in L^\i(0,T;\dbR^{n\times n})$ and $\bar F(\cd,\cd)\in L^{2,1}(\D_2(0,T);\dbR^{n\times n})$, we can choose $T_0$ to be small enough such that
\vskip-8mm
$$12(T_0-t_0)\sup\limits_{t_0\les s\les T_0}\Big\{1+|\bar C_1(s)|^2+|\bar C_2(s)|^2\Big\}+6\sup\limits_{t_0\les s\les T_0}\Big(2|\bar C_3(s)|^2+1\Big)\int_{t_0}^{T_0}\Big(\int_{t_0}^s|\bar F(s,r)| dr\Big)^2ds<1.$$
Then, $\cT$ is a contraction mapping on $[t_0,T_0]$ and \rf{SDDE-proof} admits a unique solution on $[t_0,T_0]$. Since terminal time $T$ is finite, by repeating the above steps on $[t_0,T]$, we obtain a unique solution of \rf{SDDE-proof} on $[t_0,T]$.
\epf

{\bf The proof of Proposition \ref{solvability of closed-loop system}:}
\bpf
Without loss of generality, assume that $T=t_0+n\d$ with an integer $n$. Then, on $[t_0,t_0+\d]$, the closed-loop system \rf{closed-loop system} becomes
\vskip-4mm
\begin{eqnarray}
x(t)=\xi(0)+\int_{t_0}^{t }B_2(s)\varsigma(s-\d-t_0) ds+\int_{t_0}^t\Big[A_1(s)x(s) +A_2(s)\xi(s-t_0-\d)
\qq\qq\qq\qq\q\notag\\
+A_3(s)z(s)+\Big(B_1(s)+B_2(s+\d){\bf1}_{[t_0,t-\d)}(s)
+\int_s^tB_3(r)\tilde{F}(r,s)dr\Big)u(s)+b(s)\Big]ds\qq \qq\ \notag\\
+\int_{t_0}^t\big[C_1(s)x(s) +C_2(s)\xi(s-t_0-\d)+C_3(s)z(s)+D_1(s)u(s)
+\sigma(s)\big]dW(s),\qq\qq\q\ \ \label{deduced state-ax}
\end{eqnarray}
and
\begin{eqnarray}
u(t)=K_1(t)\xi(0)+K_1(t)\int_{t_0}^{t }\3nB_2(s)\varsigma(s-\d-t_0) ds +K_1(t)\2n\int_{t_0}^t\Big[A_1(s)x(s)+A_2(s)\xi(s-t_0-\d)
\qq\qq\qq\notag\\
+A_3(s)z(s)+\Big(B_1(s)+B_2(s+\d){\bf1}_{[t_0,t-\d)}(s)
+\int_s^tB_3(r)\tilde{F}(r,s)dr\Big)u(s)+b(s)\Big]ds \qq\qq\qq\qq\notag\\
+K_1(t)\int_{t_0}^t\big[C_1(s)x(s) +C_2(s)\xi(s-t_0-\d)+C_3(s)z(s)+D_1(s)u(s)
+\sigma(s)\big]dW(s)\qq\qq\qq\q\ \notag\\
+\int_{t_0}^tK_2(t,s)x(s)ds+K_3(t) \xi(t-\d-t_0)+\int_{t_0}^tK_4(t,s)u(s)ds +v(t).\qq\qq\qq\qq\q\qq\qq\ \ \!\label{deduced-contol-ax}
\end{eqnarray}
By Fubini theorem, we have
\begin{eqnarray}
\int_{t_0}^t\int_{t_0}^sF(t,s)B_2(r)\nu(r)drds
\qq\qq\qq\qq\qq\qq\qq\qq\qq\qq\qq\qq\qq\qq \notag\\
=\int_{t_0}^{(t_0+\d)\wedge t }\Big[\int_r^tF(t,s)ds\Big]B_2(r)\varsigma(r-\d-t_0)dr
+\int_{t_0\wedge (t-\d)}^{t-\d}\Big[\int_{r+\d}^tF(t,s)ds\Big]B_2(r+\d) u(r)dr,\label{Z2-ax}
\end{eqnarray}
and
\begin{eqnarray}
\int_{t_0}^t\int_{t_0}^sF(t,s)B_3(r)\mu(r)drds
=\int_{t_0}^t\int_{t_0}^sF(t,s)B_3(r)\Big[\int_{t_0}^r\tilde{F}(r,\a)
u(\a)d\a\Big] drds\qq\qq\qq\q\notag
\\
=\int_{t_0}^t\int_{t_0}^s\int_\a^s F(t,s)B_3(r)\tilde{F}(r,\a)
u(\a)drd\a ds=\int_{t_0}^t\Big[\int_\a^t\int_\a^sF(t,s) B_3(r)\tilde{F}(r,\a)drds\Big] u(\a)d\a.\label{Z3-ax}
\end{eqnarray}
Hence, from \rf{Z2-ax}--\rf{Z3-ax} and by adding some indicative functions, we derive
\begin{eqnarray}
z(t)=\int_{t_0}^t F(t,s)\Big(\xi(0)+\int_{t_0}^{s} B_2(r)\varsigma(r-\d-t_0)dr\Big)ds +\int_{t_0}^t \cE(t,s)[A_1(s)x(s) \qq\qq\qq\notag\\[-1mm]
+A_2(s) \xi(s-t_0-\d)+A_3(s)z(s)]ds
+\int_{t_0}^t\int_s^tF(t,r) \big(B_1(s)+B_2(s+\d){\bf1}_{[0,r-\d)}(s) \qq\q
\notag\\[-1mm]
+\int_s^rB_3(\th)\ti F(\th,s)d\th\big)dru(s)ds+\int_{t_0}^t\cE(t,s)b(s)ds
+\int_{t_0}^t\cE(t,s)
[C_1(s)x(s)\qq\qq\qq\q \notag\\
\ns\ds +C_2(s) \xi(s-t_0-\d) +C_3(s)z(s)+D_1(s)u(s)+\si(s)]dW(s),\qq\qq\qq\qq\qq\qq\! \ \label{distributed delay of state-ax}
\end{eqnarray}
where $\cE(t,s)\equiv{\bf1}_{[0,t)}(s)\int_s^tF(t,r)dr$.
Let $\BX(\cd)\equiv\begin{bmatrix}x(\cd)\\z(\cd)\\u(\cd)\end{bmatrix}$. Then, by \rf{deduced state-ax}, \rf{deduced-contol-ax} and \rf{distributed delay of state-ax}, we have
\bel{BX}
\BX(t)=\BBf(t)+\int_{t_0}^t\BA(t,s)\BX(s)ds+\int_{t_0}^t\BC(t,s) \BX(s)dW(s),\q t\in(t_0,t_0+\d),
\ee
where
$$\ba{ll}\5n
\BA(t,s)\2n\equiv\2n\begin{bmatrix}\3n\begin{array}{ccc}
\ns\ds A_1(s) \1n&\1n  \1n A_3(s) \1n&\1n B_1(s)\1n+\1nB_2(s\1n+\1n\d){\bf1}_{[t_0,t-\d)}(s)
\1n+\1n\int_s^t\1nB_3(r)\tilde{F}(r,s)dr \\
\ns\ds \cE(t,s)A_1(s) \1n&\1n \cE(t,s)A_3(s) \1n&\1n \int_s^tF(t,r) (B_1(s)\1n+\1nB_2(s+\d){\bf1}_{[0,r-\d)}(s) \1n+\1n\int_s^rB_3(\th)\ti F(\th,s)d\th)dr\\
\ns\ds K_1(t)A_1(s)\1n+\1nK_2(t,s)  \1n&\1n  K_1(t)A_3(s) \1n&\1n K_1(t)(B_1(s)\1n+\1nB_2(s+\d){\bf1}_{[t_0,t-\d)}(s)
\1n+\1n\int_s^tB_3(r)\tilde{F}(r,s)dr)\1n+\1nK_4(t,s)
\end{array}
\3n\end{bmatrix},
\ea$$
$$\ba{ll}
\BC(t,s)\equiv\begin{bmatrix}\begin{array}{ccc}
\ns\ds C_1(s)  & C_3(s) & D_1(s) \\
\ns\ds  \cE(t,s)C_1(s) &  \cE(t,s)C_3(s) & \cE(t,s)D_1(s)\\
K_1(t)C_1(s)  & K_1(t)C_3(s) & K_1(t)D_1(s)
\end{array}
\end{bmatrix},
\ea$$
and
$$\ba{ll}
\BBf(t)\equiv\begin{bmatrix}
\ns\ds \xi(0)+\int_{t_0}^t\big( A_2(s)\xi(s-\d-t_0)+B_2(s)\varsigma(s-\d-t_0) +b(s)\big)ds\\
\ns\ds +\int_{t_0}^t\big(C_2(s) \xi(s-t_0-\d)+\si(s)\big)dW(s)\\
\ns\ds \\[-2mm]
\ns\ds \int_{t_0}^t F(t,s)\Big(\xi(0)+\int_{t_0}^{s} B_2(r)\varsigma(r-\d-t_0)dr\Big)ds\\
\ns\ds
+\int_{t_0}^t \cE(t,s)(A_2(s)\xi(s-\d-t_0)+b(s))ds
+\int_{t_0}^t \cE(t,s)(C_2(s)\xi(s-\d-t_0)+\si(s))dW(s)\\
\ns\ds \\[-2mm]
\ns\ds K_1(t)\xi(0)+K_1(t)\int_{t_0}^t\big( A_2(s)\xi(s-\d-t_0)+B_2(s)\varsigma(s-\d-t_0)+b(s))ds\\
\ns\ds +K_1(t)\int_{t_0}^t\big(C_2(s) \xi(s-t_0-\d)+\si(s)\big)dW(s)+K_3(t)\xi(t-\d-t_0) +v(t)
\end{bmatrix}.
\ea$$
Notice that $K_1(\cd)\in L^2(t_0,T;\dbR^{m\times n})$ and $\si(\cd)\in L^2_{\dbF}(0,T;\dbR^n)$. Then, we obtain
$$\ba{ll}
\ns\ds \dbE\int_{t_0}^T|K_1(t)\int_{t_0}^t\si(s)dW(s)|^2dt\les
\int_{t_0}^T|K_1(t)|^2\dbE\big|\int_{t_0}^t\si(s)dW(s) \big|^2dt\\
\ns\ds \qq\qq\qq\qq \les\int_{t_0}^T|K_1(t)|^2dt\dbE\int_{t_0}^T|\si(s)|^2ds <\i.
\ea$$
In addition, from (A1), $\xi(\cd)\in C([-\d,0];\dbR^n)$ and $\varsigma(\cd)\in L^2(-\d,0;\dbR^m)$, we have $\BBf(\cd)\in L^2_{\dbF}(t_0,T;\dbR^{2n+m})$. Similarly, by (A1), the Holder inequality and $(K_1(\cd),K_2(\cd,\cd),K_3(\cd),K_4(\cd,\cd), v(\cd))\in  \dbL$, we obtain  $\BA(\cd,\cd)\in L^2(\D_2(t_0,T);\dbR^{(2n+m)\times (2n+m)})$, $\BC(\cd,\cd)\in\sL^2(\triangle_2(t_0,T);$ $\mathbb{R}^{(2n+m)\times (2n+m)})$.
Thus, by Proposition \ref{Solvability of SVIE}, \rf{BX} admits a unique solution $\BX(\cd)\in L_{\dbF}^2(t_0,t_0+\d;\dbR^{2n+m})$, which implies that $x(\cd)\in L_{\dbF}^2(t_0,t_0+\d;\dbR^{n})$ and $u(\cd)\in L_{\dbF}^2(t_0,t_0+\d;\dbR^{m})$. Furthermore, by the definition of the solution of the first equation in \rf{closed-loop system}, $x(\cd)\in L_{\dbF}^2(\Omega;C([t_0,t_0+\d];\dbR^n))$, which implies the existence of the solution to \rf{closed-loop system}. As for the uniqueness, since \rf{BX} and \rf{closed-loop system} are equivalent, it can be obtained from the uniqueness of the solution to \rf{BX}. Hence, the closed-loop system \rf{closed-loop system} admits a unique solution on $[t_0,t_0+\d]$. Then, the same steps are repeated on $[t_0+\d,t_0+2\d]$, $[t_0+2\d,t_0+3\d]$ and so on. The terminal time $T$ is finite. Thus, \rf{closed-loop system} admits a unique solution on $[t_0,T]$.

The following mainly proves the estimate \rf{pre-estimate-of-controlled-SDDE}.
By the estimate of the SVIE \rf{BX}, we derive
\begin{eqnarray}
\dbE\int_{t_0}^{t_0+\d}| \BX(t)|^2dt\les
L\dbE\int_{t_0}^{t_0+\d}|\BBf(t)|^2dt\qq\qq\notag\\
\les L\Big\{\sup\limits_{t_0-\d\les t\les t_0}|\xi(t-t_0)|^2+\int_{t_0-\d}^{t_0}
|\varsigma(t-t_0)|^2dt\q\notag
\end{eqnarray}
\vskip-5mm
\begin{eqnarray}
+\dbE\2n\int_{t_0}^{t_0+\d}\2n\big(| b(t)|^2\1n+\1n| \si(t)|^2\1n+\1n|v(t)|^2\big)dt
\Big\},\label{estimate-of-mathring-X}
\end{eqnarray}
here and hereafter, $L$ is a generic constant. Hence we deduce
\vskip-6mm
\begin{eqnarray}
\dbE\int_{t_0}^{t_0+\d}| u(t)|^2dt\les\1n L\Big\{\1n\sup\limits_{t_0-\d\les t\les t_0}\2n|\xi(t-t_0)|^2+\1n\int_{t_0-\d}^{t_0}
\1n|\varsigma(t-t_0)|^2dt\ \notag\\
+\dbE\2n\int_{t_0}^{t_0+\d}\5n\big(| b(t)|^2\1n+\1n| \si(t)|^2\1n+\1n|v(t)|^2\big)dt
\Big\}.\ \label{estimate-of-mathring-u}
\end{eqnarray}
By \rf{estimate-of-mathring-u} and the estimate of the SDDE \rf{closed-loop system}, we obtain
\begin{eqnarray}
\dbE\sup\limits_{t_0\les t\les t_0+\d}| x(t)|^2
\les L\Big\{\dbE\int_{t_0}^{t_0+\d}\big(| u(t)|^2+| b(t)|^2+
|\si(t)|^2\big)dt
+\sup\limits_{t_0-\d\les t\les t_0}|\xi(t-t_0)|^2+\int_{t_0-\d}^{t_0}|
\varsigma(t-t_0)|^2dt
\Big\}\notag\\
\les L\Big\{\sup\limits_{t_0-\d\les t\les t_0}|\xi(t-t_0)|^2+\int_{t_0-\d}^{t_0}
|\varsigma(t-t_0)|^2dt
+\dbE\int_{t_0}^{t_0+\d}\5n\big(| b(t)|^2+| \si(t)|^2+|v(t)|^2\big)dt
\Big\}.\label{estimate-of-mathring-x}\qq\qq
\end{eqnarray}
Notice that $T$ is finite, and we can repeat the above steps on $[t_0+\d,t_0+2\d)$, $[t_0+2\d,t_0+3\d)$, $\cds$. Then, the estimates
\rf{estimate-of-mathring-x} and \rf{estimate-of-mathring-u} also hold on the above intervals, thus we complete the proof of \rf{pre-estimate-of-controlled-SDDE}.

\epf

{\bf The proof of Proposition \ref{prop coefficients condition}:}
\bpf
To show the equivalence, it is sufficient to prove that under (A1), the $X(\cd)$ defined in \rf{B(t,s)} is the unique solution of SVIE \rf{transformed state}. In terms of Proposition \ref{Solvability of SDDE} and Proposition \ref{Solvability of SVIE}, we only need to prove that
$$\left\{\ba{ll}
\ns\ds A(\cd,\cd)\in L^2(\D_2(0,T);\dbR^{{(3n)}\times {(3n)}}),\ B(\cd,\cd)\in L^2(\D_2(0,T);\dbR^{{(3n)}\times {m}}), \ C(\cd,\cd)\in\sL^2\left(\triangle_2(0,T); \mathbb{R}^{{(3n)}\times {(3n)}}\right),\\
\ns\ds D(\cd,\cd)\in \sL^2\left(\triangle_2(0,T);\mathbb{R}^{{(3n)}\times m}\right), \ \ \ Q(\cd) \in L^{\infty}\left(0,T;\mathbb{S}^{3n}\right),\q R(\cd)\in L^{\infty}\left(0,T; \mathbb{S}^{m}\right),\\
\ns\ds \ti b(\cd,\cd)\in L_{\mathbb{F}}^{2,1}\left(\triangle_2(0,T); \mathbb{R}^{3n}\right), \qq\  \ \ti \sigma(\cd,\cd)\in L_{\mathbb{F}}^2\left(\triangle_2(0,T); \mathbb{R}^{3n}\right).
\ea\right.$$
Notice that
$$\ba{ll}
\ns\ds \int_0^T\int_0^t|\cE(t,s)A_1(s)|^2dsdt =\int_0^T\int_0^t|\int_s^tF(t,r)drA_1(s)|^2dsdt\les \int_0^T|A_1(s)|^2ds\int_0^T\(\int_0^t|F(t,r)|dr\)^2dt.
\ea$$
Then, by $A_1(\cd)\in L^2(0,T;\dbR^{n\times n})$ and
$F(\cd,\cd)\in L^{\i}(\D_2(0,T);\dbR^{n\times n})$, we have
$\cE(\cd,\cd)A_1(\cd)\in L^2(\D_2(0,T);\dbR^{n\times n})$. In terms of Assumption (A1), $A(\cd,\cd)\in L^2(\D_2(0,T);\dbR^{{(3n)}\times {(3n)}})$. Similarly, one has $C(\cd,\cd)\in L^2(\D_2(0,T)$; $\dbR^{{(3n)}\times {(3n)}})$. By the fact of
$C_1(\cd),C_2(\cd),C_3(\cd)\in L^\i(0,T;\dbR^{n\times n})$, we deduce
$$\ba{ll}
\ns\ds \underset{s\in\left(0,T\right)}{\operatorname{ess}\sup} \left(\int_s^T|\cE(t,s)C_1(s)|^2 dt\right)^{\frac 1 2}\les M\Big(\int_0^T(\int_0^t|F(t,r)|dr)^2dt\Big)^{\frac 1 2}<\infty,
\ea$$
where $M$ is a generic constant. For any $\varepsilon>0$, by $F(\cd,\cd)\in L^{\i}(\D_2(0,T);\dbR^{n\times n})$, there exists a finite partition $\left\{a_i\right\}_{i=0}^m$ of $\left(0,T\right)$ with $0=a_0<a_1<\cdots<a_m=T$ such that
$$\ba{ll}
\ns\ds \underset{t \in\left(a_i, a_{i+1}\right)}{\operatorname{ess} \sup }\left(\int_t^{a_{i+1}}|\cE(s,t)C_1(t)|^2 ds\right)^{\frac 1 2}\les M\(\int_{a_i}^{a_{i+1}}\(\int_{a_i}^s|F(s,r)|dr\)^2ds\)^ {\frac 1 2}<\varepsilon,
\ea$$
which implies that $C(\cd,\cd)\in\sL^2(\triangle_2(0,T); \mathbb{R}^{{(3n)}\times {(3n)}})$. Using the Holder inequality and Assumption (A1), we can prove the integrability of $B(\cd,\cd)$, $D(\cd,\cd)$, $Q(\cd)$, $R(\cd)$, $\tilde b(\cd,\cd)$ and $\tilde{\si}(\cd,\cd)$, thus the proof is completed. \epf

{\bf The proof of Theorem \ref{solvability of Riccati}:}
\bpf
Inspired by \cite{Hamaguchi-Wang-2024-II}, let us introduce the following Riccati system:
\bel{eq_Riccati--Volterra-0}\left\{\ba{ll}
\ns\ds	P^{(1)}(t)=Q(t)+(C^\top\lint P\rint C)(t)\\
\ns\ds	\hspace{1.5cm}-(C^\top\lint P\rint D)(t)(R(t)+(D^\top\lint P\rint D)(t))^{-1}(D^\top\lint P\rint C)(t),\hspace{0.3cm}0<t<T,\\
\ns\ds
	P^{(2)}(s,t,t)=P^{(2)}(t,s,t)^\top\\
\ns\ds
	=(P\rint A)(s,t)-(P\rint B)(s,t)(R(t)+(D^\top\lint P\rint D)(t))^{-1}(D^\top\lint P\rint C)(t),\q 0<t<s<T,\\
\ns\ds	\dot{P}^{(2)}(s_1,s_2,t)=(P\rint B)(s_1,t)(R(t)+(D^\top\lint P\rint D)(t))^{-1}(B^\top\lint P)(s_2,t),\q 0<t<(s_1\wedge s_2)<T,
\ea\right.\ee
%
%
where for each $M_1: \triangle_2(0, T) \rightarrow \mathbb{R}^{d_1 \times (3n)}$, $M_2: \triangle_2(0, T) \rightarrow$ $\mathbb{R}^{(3n) \times d_2}$ with any positive integers $d_1, d_2$,
\begin{eqnarray}
\left(M_1 \ltimes P\right)(s, t)\equiv M_1(s, t) P^{(1)}(s)+\int_t^T M_1(r, t) P^{(2)}(r, s, t) \mathrm{d} r, \quad 0<t<s<T,\label{left product}\qq\qq\qq\qq\q\\
\left(P \rtimes M_2\right)(s, t)\equiv P^{(1)}(s) M_2(s, t)+\int_t^T P^{(2)}(s, r, t) M_2(r, t) \mathrm{d} r, \quad 0<t<s<T,\label{right product}\qq\qq\qq\qq\q\\
\left(M_1\2n \ltimes \2n P \2n\rtimes \2n M_2\right)\1n(t)\equiv\2n\int_t^T \3n\2n M_1(s, t) P^{(1)}(s) M_2(s, t) d s\1n+\3n\int_t^T\3n\2n \int_t^T \3n\1n M_1\1n\left(s_1, t\right) P^{(2)}\2n\left(s_1, s_2, t\right) \1nM_2\1n\left(s_2, t\right) d s_1 d s_2, 0\1n<\1nt\1n<\1nT.\ \label{left-right product}
\end{eqnarray}
By Corollary 6.7 in \cite{Hamaguchi-Wang-2024-II}, the equation \eqref{eq_Riccati--Volterra-0} admits a unique solution $(P^{(1)},P^{(2)})\in\Pi(0,T)$ such that $R\mathalpha{+}(D^\top\lint P\rint D)>\b I$ for some constant $\b>0$.
\no Next we will equivalently transform \eqref{eq_Riccati--Volterra-0} into \rf{eq_Riccati--Volterra}. To this end, we first show that
\begin{eqnarray}
(C^\top\lint P\rint C)(t) =(C_1(t),C_2(t),C_3(t))^\top \cG_1(t)(C_1(t),C_2(t),C_3(t)),\label{CPC}\\
\ \ (C^\top\lint P\rint D)(t) =(C_1(t),C_2(t),C_3(t))^\top\cG_1(t)D_1(t),\qq\qq\q\ \ \ \label{CPD}\\
(D^\top\lint P\rint D)(t)=D_1(t)^\top \cG_1(t)D_1(t),\qq\qq\qq\qq\qq\qq \label{DPD}\\
(D^\top\lint P\rint C)(t)
=D_1(t)^\top\cG_1(t)(C_1(t),C_2(t),C_3(t)).\qq\qq\qq\ \!\label{DPC}
\end{eqnarray}
Notice that
$$\ba{ll}
C(s,t)=\Upsilon(s,t)(C_1(t),C_2(t),C_3(t)),
\ea$$
and
$$(C^\top\lint P\rint C)(t)=\int_t^TC(s,t)^\top P^{(1)}(s)C(s,t)ds
+\int_t^T\int_t^TC(s_1,t)^\top P^{(2)}(s_1,s_2,t)C(s_2,t)ds_1ds_2.$$
Then, \rf{CPC} holds. Similarly, by $D(s,t)=\Upsilon(s,t)D_1(t)$, we obtain \rf{CPD}--\rf{DPC}.

Next it is time to treat the terms $(P\rtimes A)$ and $(P\rtimes B)$.
From \rf{right product} and $A(s,t)=\Upsilon(s,t)(A_1(t),A_2(t),A_3(t))$, we have
\bel{PA-}
(P\rtimes A)(s,t)=\Big[P^{(1)}(s)\Upsilon(s,t)+\int_t^TP^{(2)}(s,r,t)
\Upsilon(r,t)dr\Big](A_1(t),A_2(t),A_3(t)).
\ee
On the other hand, we observe that
$$\ba{ll}
\ns\ds B(t,s)=\int_s^t\Pi(t,s,\th)(B_1(s)^\top,B_2(s+\d)^\top, (B_3(\th) \tilde{F}(\th,s))^\top)^\top d\th.
\ea$$
Then, it is easy to see,
\begin{eqnarray}
(P\rtimes B)(s,t)=\int_t^sP^{(1)}(s)\Pi(s,t,\th)(B_1(t)^\top, B_2(t+\d)^\top,(B_3(\th)\tilde{F}(\th,t))^\top)^\top d\th\qq\q\notag\\
+\int_t^T\int_t^rP^{(2)}(s,r,t)\Pi(r,t,\th)(B_1(t)^\top\2n, B_2(t+\d)^\top\2n,(B_3(\th)\tilde{F}(\th,t))^\top\1n)^\top\2n d\th dr. \label{PB-}
\end{eqnarray}
 Hence, after some direct calculations, together with \rf{CPC}--\rf{PB-}, we see that \rf{eq_Riccati--Volterra-0} can be written as \rf{eq_Riccati--Volterra}, which completes the proof of Theorem \ref{solvability of Riccati}.
\epf

{\bf The proof of Theorem \ref{solvability of BSVIE}:}
\bpf
Introduce the following Type-II extended backward SVIE:
\bel{EBSVIE-0}\left\{\ba{ll}
\ns\ds d\eta(t,s)=-\Big\{(P\rtimes\ti b)(t,s)+ \G^*(t,s)^\top(D^\top\ltimes P\rtimes \ti\si)(s)+ \G^*(t,s)^\top\int_s^TB(r,s)^\top\eta(r,s)dr\\
\ns\ds \qq\qq+\check\G(t,s)^\top\int_s^T D(r,s)^\top\zeta(r,s)dr\Big\}ds+\zeta(t,s)dW(s),\q 0<s<t<T,\\
\ns\ds \eta(t,t)=(C^\top\ltimes P\rtimes \ti\si)(t)+\Xi^*(t)^\top(D^\top\ltimes P\rtimes \ti\si)(t) +\int_t^T(A(r,t)+B(r,t)\Xi^*(t))^\top\eta(r,t)dr \\
\ns\ds \qq\qq+\int_t^T(C(r,t)+D(r,t)\Xi^*(t))^\top\zeta(r,t)dr,\q 0<t<T,
\ea\right.\ee
where
\bel{check-Xi}\ba{ll}
\ns\ds  \Xi^*(t)\equiv-(R(t)+(D^\top\ltimes P\rtimes D)(t))^{-1}(D^\top\lint P\rint C)(t),\ 0<t<T,
\ea\ee
\bel{check-G}\ba{ll}
 \Gamma^*(s,t)\equiv-(R(t)+(D^\top\lint P\rint D)(t))^{-1}(B^\top\lint P)(s,t),\ 0<t<s<T.
\ea\ee
Under Assumptions (A1)--(A2), it follows from Theorem 3.2 in \cite{Hamaguchi-Wang-2024-II} that Equation \rf{EBSVIE-0} admits a unique solution $(\eta,\zeta)\in L^2_{\dbF,\mathrm{c}}(\triangle_2(0,T);\dbR^{3n})\times L^2_\dbF(\triangle_2(0,T);\dbR^{3n})$.
After a careful observation, we have that
$$\ba{ll}
\ns\ds (P\rtimes \ti b)(t,s)=\Big[P^{(1)}(t)\Upsilon(t,s)+\int_s^TP^{(2)} (t,r,s)
\Upsilon(r,s)dr\Big]b(s),\\
\ns\ds (C^\top\ltimes P\rtimes\ti\si)(t)=(C_1(t),C_2(t),C_3(t))^\top\cG_1(t)\si(t),
\ea$$
and
$$\ba{ll}
\ns\ds (D^\top\ltimes P\rtimes\ti\si)(s)=D_1(s)^\top\cG_1(s)\si(s),\q C(t,s)=\Upsilon(t,s)(C_1(s),C_2(s),C_3(s)),\\
\ns\ds D(t,s)=\Upsilon(t,s)D_1(s),\q B(t,s)=\int_s^t\Pi(t,s,\th)(B_1(s)^\top,B_2(s+\d)^\top, (B_3(\th) \tilde{F}(\th,s))^\top)^\top d\th.
\ea$$
Then, we see that \rf{EBSVIE-0} can be rewritten as \rf{EBSVIE}. It is then easy to see the conclusion of Theorem \ref{solvability of BSVIE}.

\epf

\ss

Before proving Theorem \ref{main theorem optimal control}, according to \cite{Hamaguchi-Wang-2024-II}, we need some auxiliary results about the new control problem with the state equation \rf{transformed state} and the cost functional \rf{transformed cost}.

To begin with, for any given $t_0\in[0,T)$, we consider the following system of the new control problem:
\vskip-3mm
\bel{integral closed-loop system}\ba{ll}
\left\{\begin{aligned}
& X^{t_0,\xi,\varsigma}(t)=\f(t)+ \int_{t_0}^t\big[A(t, s) X^{t_0,\xi,\varsigma}(s)+B(t,s) u^{t_0,\xi,\varsigma}(s)+\ti b(t,s)\big] d s \\
& \quad+\int_{t_0}^t\big[C(t, s) X^{t_0,\xi,\varsigma}(s)+D(t, s) u^{t_0,\xi,\varsigma}(s)+\ti \sigma(t, s)\big] d W(s), \ \ t_0<t<T, \\
& \Theta^{t_0,\xi,\varsigma}(s, t)=\f(s)+ \int_{t_0}^t\big[A(s, r) X^{t_0,\xi,\varsigma}(r)+B(s, r) u^{t_0,\xi,\varsigma}(r)+\ti b(s, r)\big] d r \\
& \quad+\int_{t_0}^t\big[C(s, r) X^{t_0,\xi,\varsigma}(r)+D(s, r) u^{t_0,\xi,\varsigma}(r)+\ti \sigma(s, r)\big] d W(r),\ \ t_0\1n<\1nt\1n<\1ns\1n<\1nT, \\
& u^{t_0,\xi,\varsigma}(t)=\Xi(t) X^{t_0,\xi,\varsigma}(t)+\int_t^T \Gamma(s, t) \Theta^{t_0,\xi,\varsigma}(s, t) \mathrm{d} s+\o(t), \ \  t_0<t<T.
\end{aligned}\right.
\ea\ee
In terms of \cite{Hamaguchi-Wang-2024-II}, we call any triplet $(\Xi, \Gamma, \o) \in \mathcal{S}(t_0, T)\equiv L^{\infty}\left(t_0, T ; \mathbb{R}^{m \times n}\right) \times L^2\left(\triangle_2(t_0, T) ; \mathbb{R}^{m \times n}\right) \times L_{\dbF}^2(t_0,T;\dbR^m)$ the \emph{causal feedback strategy}.
For any $(\Xi, \Gamma, \o) \in \mathcal{S}(t_0, T)$,  $\xi\in C([-\d,0];\dbR^n)$ and $\varsigma\in L^2(-\d,0;\dbR^m)$, let the triplet $(X^{t_0,\xi,\varsigma}, \Theta^{t_0,\xi,\varsigma}$, $u^{t_0,\xi,\varsigma})$ be the solution to the system \rf{integral closed-loop system} and write $u^{t_0,\xi,\varsigma}=(\Xi, \Gamma, w)\left[t_0, \xi,\varsigma\right]$.
A causal feedback strategy $(\Xi^*, \Gamma^*, \o^*) \in \mathcal{S}(t_0, T)$ is called a \emph{causal feedback optimal strategy} of the new control problem if
$$
J\big(t_0,\xi,\varsigma;(\Xi^*, \Gamma^*, \o^*)\left[t_0, \xi,\varsigma\right]\big) \les J\big(t_0,\xi,\varsigma;u\big),
$$
for any $\left( \xi,\varsigma\right) \in   C([-\d,0];\dbR^n)\times L^2(-\d,0;\dbR^m)$ and any $u(\cd)\in L^2_{\dbF}(t_0,T;\dbR^m)$.
%

\ss

The following result gives the closed-loop solvability for the new control problem.

\bl{closed-loop solvability-decompose}
Let Assumptions {\rm (A1)--(A2)} hold. Then, for any given $t_0\in[0,T)$, the causal feedback optimal strategy $(\Xi^*,\Gamma^*,\o^*)$ of the new control problem on $[t_0,T]$ with \rf{transformed state} and  \rf{transformed cost}, is given by \rf{check-Xi-decompose} and for $t_0<t<T$,
\begin{eqnarray}
{\o}^*(t)=-\cR(t)^{-1} \(D_1(t)^\top \cG_1(t)\si(t)+\int_t^T\[\int_t^s\cB(\th,t)^\top \Pi(s,t,\th)^\top\eta(s,t)d\th +D_1(t)^\top\Upsilon(s,t)^\top\zeta(s,t)\]ds\), \q \label{check-o-decompose}
\end{eqnarray}
\no where $\Pi(\cd,\cd,\cd)$, $\Upsilon(\cd,\cd)$, $\cR(\cd)$, $\cG_1(\cd)$ and $\cG_3(\cd,\cd,\cd)$ are defined by \rf{Pi},\rf{Upsilon-cR}, respectively.
\el

\bpf
First, given $(\eta,\zeta)$ satisfying \rf{EBSVIE}, we define
\bel{k}\ba{ll}
\ns\ds k(t)\mathalpha{\equiv}\big(D^\top\lint P\rint\ti \sigma\big)(t)\mathalpha{+}\int^T_tB(s,t)^\top\eta(s,t)\ d s\mathalpha{+}\int^T_tD(s,t)^\top\zeta(s,t) d s,\\
\ns\ds \o^*(t)\equiv-\big(R(t)\mathalpha{+}(D^\top\lint P\rint D)(t)\big)^{-1} k(t),\ 0<t<T.
\ea
\ee
Under (A1)--(A2), by Theorem 5.4 in \cite{Hamaguchi-Wang-2024-II}, and Theorem \ref{solvability of Riccati} and Theorem \ref{solvability of BSVIE} in the current paper, we see that $(\Xi^*,\Gamma^*,\o^*)$ is the optimal causal feedback strategy where
$(\Xi^*,\Gamma^*)$ is in \rf{check-Xi-decompose}. Based on this fact, it is sufficient to show that $\o^*(\cd)$ defined in \rf{k} can be rewritten in the form of \rf{check-o-decompose}.

In fact, by \rf{B(t,s)} and \rf{Pi}, we obtain
$$
B(s,t)=\int_t^s\Pi(s,t,\th)\big(B_1(t)^\top,B_2(t+\d)^\top,\ti F(\th,t)^\top B_3(\th)^\top\big)^\top d\th{\bf1}_{(0,s)}(t).
$$
Recall that
$$
\left(B^\top \ltimes P\right)(s, t)\equiv B(s, t)^\top P^{(1)}(s)+\int_t^T B(r, t)^\top P^{(2)}(r, s, t) \mathrm{d} r, \quad0<t<s<T.
$$
Then, we deduce
\begin{eqnarray}
(B^{\top} \ltimes P)(s, t)=\int_t^s\big(B_1(t)^\top,B_2(t+\d)^\top,\ti F(\th,t)^\top B_3(\th)^\top\big)\Pi(s,t,\th)^\top d\th P^{(1)}(s)\qq\qq\qq\qq\notag\\
\ns\ds +\int_t^T\3n\int_t^r\big(B_1(t)^\top\3n,B_2(t+\d)^\top,\ti F(\th,t)^\top B_3(\th)^\top\big)\Pi(r,t,\th)^\top  P^{(2)}(r,s,t)d\th dr.\label{BP-prove}
\end{eqnarray}
Notice that $D(t,s)=\Upsilon(t,s)D_1(s)$ and
$$
(D^\top\lint P\rint \ti\si)(t)
=D_1(t)^\top\cG_1(t)\si(t),\ \ t\in[t_0,T].
$$
Combining it with \rf{k}, we then have that
\begin{eqnarray}
k(t)=D_1(t)^\top \cG_1(t)\si(t)+\int_t^T \[\int_t^s\big(B_1(t)^\top,B_2(t+\d)^\top,\ti F(\th,t)^\top B_3(\th)^\top\big) \Pi(s,t,\th)^\top\eta(s,t)d\th \qq\q\notag\\
\ns\ds +D_1(t)^\top\Upsilon(s,t)^\top\zeta(s,t)\]ds. \qq\qq\qq\qq\qq\qq\qq\qq\qq\qq\qq\qq\q\label{k-prove}
\end{eqnarray}
This naturally implies the desired conclusion of Lemma \ref{closed-loop solvability-decompose}.
\epf

{\bf The proof of Theorem \ref{main theorem optimal control}:}
\bpf
The basic idea of the following arguments is to explicitly construct the desired five-tuple closed-loop strategy by the causal feedback strategy $(\Xi^*,\G^*,\o^*)$ in Lemma \ref{closed-loop solvability-decompose}. To this end, we divide the proof into six steps.

Step 1:  \no Given $\Xi^*(\cd),\Gamma^*(\cd,\cd)$ in \rf{check-Xi-decompose}, we decompose them as follows for later convenience:
$$\ba{ll}
\ns\ds \Xi^*(t) =\big[\Xi^*_1(t), \Xi_2^*(t), \Xi_3^*(t)\big],
 \G^*(s,t) =\big[\G_1^*(s,t),\G_2^*(s,t), \G_3^*(s,t)\big].
\ea
$$
In this step, we prove that the following process is an optimal closed-loop outcome control of Problem (P) on $[t_0,T]$:
\begin{eqnarray}
 u^*(t)=K_1^*(t)x^*(t)
+\int_{t_0}^tK_2^*(t,s)x^*(s)ds+K_3^*(t)x^*(t-\d)
 +\int_{t_0}^tK_4^*(t,s) u^*(s)ds+ v^*(t),
\label{transformed optimal control-ap}
\end{eqnarray}
where
\begin{eqnarray}
K_1^*(t)={\Xi}_1^*(t)+\int_t^T\[{\G}_1^*(s,t) +{\G}_2^*(s,t){\bf1}_{\big([t+\d]\wedge T,\i\big)}(s)+\int_t^s{\G}_3^*(s,t)F(s,\th)d\th\] ds,\qq\qq\qq\q\q\  \label{K1-}\\
K_2^*(t,s)={\Xi}_3^*(t)F(t,s) +{\G}_2^*(s+\d,t){\bf1}_{\big[t-\d,T-\d\big]}(s) +\int_t^T{\G}_3^*(\th,t)F(\th,s)d\th,\qq\qq\qq \qq\qq\qq\label{K2-}\\
K_3^*(t)={\Xi}_2^*(t),\qq\qq\qq\qq\qq\qq\qq\qq\qq\qq\qq \qq\qq\qq\qq\qq\qq\qq\ \ \ \ \ \label{K3-}\\
K_4^*(t,s)\1n=\1n\Big\{\2n\int_{s+\d}^T\2n \Big[{\G}_1^*(r,\1nt) \1n+\1n{\G}_2^*(r,\1nt) {\bf1}_{(s+2\d,\i)}(r) {\bf1}_{[0,T-\d)}(t) \2n+\3n\int_{s+\d}^r\5n{\G}_3^*(r,t)F(r,\1n\th) d\th\Big]\1ndr B_2({s\1n+\1n\d})\Big\} \qq\qq\qq\notag\\
\times{\bf1}_{[t-\d,T-\d]}(s)+\int_t^T \3n\int_{\th'}^T\3n\Big({\G}_1^*(r,t)+{\G}_2^*(r,t) {\bf1}_{[t_0,T-\d)}(\th'){\bf1}_{[t_0,T-\d)}(t) {\bf1}_{(\th'+\d,\i)}(r) \qq\qq\qq\qq\notag\\
+\int_{\th'}^r{\G}_3^*(r,t)F(r,\th)d\th\Big) B_3(\th')\ti F(\th',s)drd\th',\q\qq\qq\qq\qq\qq\qq\qq\qq\qq\qq\ \label{K4-}\\
v^*(t)=\o^*(t)\2n+\3n\int_t^{T\wedge (t_0+\d)}\3n{\Gamma}^*_2(s,t) \xi(s-t_0-\d)ds +\int_{t-\d}^{t_0}\Big\{ \int_{s+\d}^T\Big[{\G}_1^*(r,t) +{\G}_2^*(r,t) {\bf1}_{(s+2\d,\i)}(r) \ \ \qq\qq \notag\\
\times {\bf1}_{[0,T-\d)}(t)+\int_{s+\d}^r{\G}_3^*(r,t)F(r,\th)d\th\Big]drB_2(s+\d)\Big\} \varsigma(s-t_0)ds{\bf1}_{[t_0,t_0+\d]}(t),  \qq\qq\qq\qq\ \ \label{v-}
\end{eqnarray}
\no and $\o^*(\cd)$ is defined by  \rf{check-o-decompose}.
%

According to Proposition \ref{prop coefficients condition}, given $(\Xi^*,\G^*,\o^*)$ in Lemma  \ref{closed-loop solvability-decompose},  the following process is an optimal control of Problem (P):
\bel{causal feedback}\ba{ll}
\ns\ds u^*(t)=\Xi^*(t)X^*(t)+\int_t^T\G^*(s,t)\Th^*(s,t)ds +\o^*(t),\ t_0<t<T,
\ea\ee
where $X^*(\cd)\equiv \big[x^*(\cd)^\top,y^*(\cd)^\top,z^*(\cd)^\top\big]^\top$, and
\begin{eqnarray}
\ds \Th^*(s,t)=\f(s)+\int_{t_0}^t\bigg[A(s,r)X^*(r) +B(s,r)u^*(r)+\ti b(s,r)\bigg]dr\qq\qq\q\  \notag\\
\ns\ds \qq\q+\int_{t_0}^t\bigg[C(s,r)X^*(r)+D(s,r)u^*(r) +\ti \si(s,r)\bigg]dW(r),\q t_0<t<s<T.\notag
\end{eqnarray}
Our next idea is to rewrite \rf{causal feedback} into \rf{transformed optimal control-ap}. For later convenience, let $%
\Th^*(\cd,\cd)=\big[\Th_1^*(\cd,\cd)^\top,
\Th_2^*(\cd,\cd)^\top,\Th_3^*(\cd,\cd)^\top\big]^\top.$ It then yields
\begin{eqnarray}
u^*(t)=\Xi_1^*(t)x^*(t)+\Xi_2^*(t)x^*(t-\d)
+\Xi_3^*(t)\int_{t_0}^tF(t,s)x^*(s)ds\qq\qq\qq\qq\qq\notag\\
+\int_t^T\Big[\G_1^*(s,t)\Th_1^*(s,t) +\G_2^*(s,t)\Th_2^*(s,t) +\G_3^*(s,t)\Th_3^*(s,t) \Big]ds+w^*(t),\qq\q \label{decomposed causal feedback}
\end{eqnarray}
where
\begin{eqnarray}
\Th_1^*(s,t)=\xi(0)+\int_{t_0}^{(t_0+\d)\wedge s}B_2(r)\varsigma(r-t_0-\d)dr+\int_{t_0}^t \Big[A_1(r)x^*(r)+A_2(r)y^*(r)+A_3(r)z^*(r)
\qq\qq\notag\\
+\ \! B_1(r)u^*(r)+B_2(r+\d){\bf1}_{[t_0,s-\d)}(r)u^*(r) +\int_r^sB_3(\th)\ti F(\th,r)d\th u^*(r)+b(r)\Big]dr
\qq\qq\qq\ \notag\\
+\int_{t_0}^t\Big[C_1(r)x^*(r)+C_2(r)y^*(r) +C_3(r)z^*(r)+D_1(r)u^*(r)
+\si(r)\Big]dW(r),\ t_0<t<s<T,\q\label{th1}
\end{eqnarray}
and
\begin{eqnarray}
\Th_2^*(s,t)=\xi(s-\d-t_0){\bf1}_{[t_0,t_0+\d]}(s) +{\bf1}_{(t_0+\d,\i)}(s)
\Big\{\xi(0)+\int_{t_0}^{(t_0+\d)\wedge(s-\d)} B_2(r)\varsigma(r-t_0-\d)dr\Big\}\qq\
\notag\\
+\int_{t_0}^{t\wedge(s-\d)}\Big[A_1(r)x^*(r) +A_2(r)y^*(r)+A_3(r)z^*(r) +B_1(r)u^*(r)+B_2(r+\d)u^*(r){\bf1}_{[t_0,s-2\d)}(r) \!\notag\\
+
\int_r^{s-\d}B_3(\th)\ti F(\th,r)d\th u^*(r)+b(r)\Big]dr{\bf1}_{(t_0+\d,\i)}(s) +\int_{t_0}^{t\wedge(s-\d)}\Big[C_1(r)x^*(r) +C_2(r)y^*(r)\qq\ \!\notag\\
+C_3(r)z^*(r)+D_1(r)u^*(r)
+\si(r)\Big]dW(r){\bf1}_{{(t_0+\d,\i)}}(s),\  t_0<t<s<T,\qq\qq\qq\qq\qq \label{th2}
\end{eqnarray}
and
\begin{eqnarray}
\Th_3^*(s,t)=\int_{t_0}^sF(s,r)\Big(\xi(0) +\int_{t_0}^{t_0+\d} B_2(\th)\varsigma(\th-t_0-\d) {\bf1}_{[t_0,r)}(\th)d\th\Big)dr+\int_{t_0}^t \Big[\cE(s,r) \big(A_1(r)x^*(r)+A_2(r)y^*(r)\notag\\
+A_3(r)z^*(r) +B_1(r)u^*(r)\big)+\cE(s,r+\d)B_2(r+\d)u^*(r)
+\int_r^s\cE(s,\th)B_3(\th)
\ti F(\th,r)d\th u^*(r)\qq\ \ \!
\notag\\
+\cE(s,r)b(r)\Big]dr+\int_{t_0}^{t}\cE(s,r) \Big[C_1(r)x^*(r)
+C_2(r)y^*(r)
+C_3(r)z^*(r)\qq\qq\qq\qq\qq\q\ \notag\\
+D_1(r)u^*(r)
+\si(r)\Big]dW(r),\ t_0<t<s<T.\qq\qq\qq\qq\qq\qq\qq\qq\qq\qq\ \ \ \notag 
\end{eqnarray}
As to $\Th_1^*$, $\Th_2^*$, $\Th_3^*$, by the Fubini theorem, we obtain
\bel{th1-}\ba{ll}
\ns\ds \Th_1^*(s,t)=x^*(t)+\int_t^{s\wedge(t+\d)}B_2(r) u^*(r-\d)dr
+\int_t^sB_3(r)\int_{t_0}^t\ti F(r,\th)u^*(\th)d\th dr,
\ea\ee
and
\begin{eqnarray}
 \Th_2^*(s,t)=\xi(s-\d-t_0){\bf1}_{[t_0,t_0+\d]}(s)
+{\bf1}_{(\d+t_0,\i)}(s)\Big\{x^*(t\wedge(s-\d))
\notag\qq\q\ \\
+\int_{t\wedge(s-\d)}^{(t+\d)\wedge(s-\d)} B_2(r)u^*(r-\d)dr+\int_{t\wedge(s-\d)}^{s-\d}B_3(r)\int_{t_0}^{ t}\ti F(r,\th)u^*(\th)d\th dr\Big\},\label{th2-}
\end{eqnarray}
and
\begin{eqnarray}
\Th_3^*(s,t)=\int_{t_0}^sF(s,r)x^*(t\wedge r)dr\qq\qq\qq
\qq\qq\qq\qq\qq\qq\qq\qq\qq\qq\notag\\
+\int_{t_0}^sF(s,r)\Big[\int_{t\wedge r}^{(t+\d)\wedge r}B_2(\th)u^*(\th-\d)d\th+\int_{t_0}^{t\wedge r}\int_{t\wedge r}^r
B_3(\th)\ti F(\th,\th')u^*(\th')d\th d\th'\Big]dr.\label{th3-}
\end{eqnarray}
By \rf{decomposed causal feedback}, \rf{th1-}--\rf{th3-}, for $t\in[t_0,T]$, we derive
\vskip-4mm
\begin{eqnarray}
u^*(t)=\Big[\Xi_1^*(t)+\int_t^T\G_1^*(s,t)ds +\int_{(t+\d)\wedge T}^T\G_2^*(s,t)ds+\int_t^T\int_r^T \G_3^*(s,t)F(s,r)dsdr\Big]x^*(t)\qq\qq\qq\notag\\
+\Xi_2^*(t)x^*(t-\d) +\int_t^T\G_2^*(s,t) \xi(s-t_0-\d){\bf1}_{[t_0,t_0+\d]}(s)ds \qq\qq\qq\qq\qq\qq\qq\qq\q\!\notag\\
+\int_{t_0}^t\Big[\Xi_3^*(t)F(t,s) +\G_2^*(s+\d,t){\bf1}_{[t-\d,T-\d]}(s) +\int_t^T\G_3^*(r,t)F(r,s)dr\Big]x^*(s)ds\ \ \q\qq\qq\qq\ \!\notag\\[-1mm]
+\int_t^{(t+\d)\wedge T}\Big[\int_r^T\G_1^*(s,t)B_2(r)ds +\int_{r+\d}^T\G_2^*(s,t) B_2(r)ds{\bf1}_{[t_0,T-\d)}(t)\qq\qq\qq\qq\qq\q\ \ \ \ \notag\\[-1mm]
+\1n\int_r^T\3n\int_r^{s\wedge(t+\d)}\3n\G_3^*(s,t) F(s,\th)B_2(r)d\th ds\2n+\3n\int_t^T\3n\1n\int_{s\wedge(t+\d)}^s \3n\G_3^*(s,t)F(s,\th)B_2(r)d\th ds{\bf1}_{[0,T-\d)}(t)\1n\Big]\1nu^*\1n(r\1n-\1n\d)dr \ \ \! \notag\\[-1mm]
+\int_{t_0}^t\Big[\int_t^T\int_\th^T \G_1^*(s,t)B_3(\th)\ti F(\th,r)dsd\th+\int_{t+\d}^T\int_t^{s-\d} \G_2^*(s,t)B_3(\th)\ti F(\th,r)d\th ds{\bf1}_{[t_0,T-\d)}(t)\qq\ \ \notag\\[-1mm]
+\int_t^T\int_t^s\int_t^{\th'}\G_3^*(s,t) F(s,\th')B_3(\th)\ti F(\th,r)d\th d\th'ds\Big]u^*(r)dr+ w^*(t).\qq\qq\qq\qq\qq\qq\ \ \!\!
\notag
\end{eqnarray}
Notice that
$$\ba{ll}
\ns\ds \int_t^{(t+\d)\wedge T}\Big[\int_r^T\G_1^*(s,t)B_2(r)ds +\int_{r+\d}^T\G_2^*(s,t) B_2(r)ds{\bf1}_{[0,T-\d)}(t)\qq\qq\qq\qq\qq\q\ \ \ \ \notag\\
\ns\ds +\1n\int_r^T\3n\int_r^{s\wedge(t+\d)}\3n\G_3^*(s,t) F(s,\th)B_2(r)d\th ds\2n+\3n\int_t^T\3n\1n\int_{s\wedge(t+\d)}^s \3n\G_3^*(s,t)F(s,\th)B_2(r)d\th ds{\bf1}_{[0,T-\d)}(t)\1n\Big]\1nu^*\1n(r\1n-\1n\d)dr\\
\ns\ds=\int_{t}^{(t+\d)\wedge T}\Big[\int_r^T\G_1^*(s,t)ds +\int_{r+\d}^T\G_2^*(s,t)ds
{\bf1}_{[t_0,T-\d)}(t) +\int_r^T\int_r^s\G_3^*(s,t)F(s,\th)d\th ds\Big]B_2(r)u^*(r-\d)dr\\
\ns\ds =\int_{t}^{(t+\d)\wedge T}\int_r^T\Big[\G_1^*(s,t) +\G_2^*(s,t)
{\bf1}_{(r+\d,T]}(s){\bf1}_{[t_0,T-\d)}(t) +\int_r^s\G_3^*(s,t)F(s,\th)
d\th\Big]ds B_2(r)u^*(r-\d)dr,
\ea$$
and
$$\ba{ll}
\ns\ds \int_{t_0}^t\Big[\int_t^T\int_\th^T \G_1^*(s,t)B_3(\th)\ti F(\th,r)dsd\th+\int_{t+\d}^T\int_t^{s-\d} \G_2^*(s,t)B_3(\th)\ti F(\th,r)d\th ds{\bf1}_{[t_0,T-\d)}(t)\qq\ \ \notag\\
\ns\ds +\int_t^T\int_t^s\int_t^{\th'}\G_3^*(s,t) F(s,\th')B_3(\th)\ti F(\th,r)d\th d\th'ds\Big]u^*(r)dr\\
\ns\ds =\1n\int_{t_0}^t\2n\int_t^T\3n\1n \Big(\1n\int_\th^T\3n
\G_1^*(s,t)ds\1n
+\3n\int_{\th+\d}^T\3n\2n\G_2^*(s,t)ds{\bf1}_{[t_0,T-\d)}
(\th){\bf1}_{[t_0,T-\d)}(t)
\1n+\3n\int_\th^T\3n\2n\G_3^*(s,t)\3n
\int_{\th}^s\3nF(s,\th')d\th' ds\1n\Big)B_3(\th)\ti F(\th,r)d\th u^*(r)dr\\
\ns\ds =\int_{t_0}^t\int_t^T\int_\th^T \Big(\G_1^*(s,t)
+\G_2^*(s,t){\bf1}_{[0,s-\d)}(\th) {\bf1}_{[t_0,T-\d)}(t)
+\G_3^*(s,t)\int_\th^sF(s,\th')d\th'\Big)ds B_3(\th)\ti F(\th,r)d\th u^*(r)dr.
\ea$$
Then, the optimal closed-loop outcome control \rf{causal feedback} of Problem (P) becomes \rf{transformed optimal control-ap}, and $K_1^*(\cd),K_2^*(\cd,\cd), K_3^*(\cd)$, $K_4^*(\cd,\cd), v^*(\cd)$ are given by \rf{K1-}--\rf{v-}.

\ss

Step 2: In this step, we give furthermore explicit representation of $K_1^*$ by means of $(P^{(1)},P^{(2)})$ and other given coefficients of the optimal control problem.
Let $P=(P^{(1)},P^{(2)})$ be the solution to the Riccati--Volterra equation \eqref{eq_Riccati--Volterra}, and decompose them as $P^{(1)}(\cd)=\big(P^{(1)}_{ij}(\cd)\big)_{1\les i,j\les 3}$ and $P^{(2)}(\cd)=\big(P^{(2)}_{ij}(\cd)\big)_{1\les i,j\les 3}$. Then, by \rf{check-Xi-decompose} and \rf{Pi}, for $i=1,2,3$, we have
\bel{Ga_i}\ba{ll}
\ns\ds \Xi_i^*(t)=-\cR(t)^{-1} D_1(t)^\top \cG_1(t)C_i(t),\q t_0\1n<\1nt\1n<\1nT, \\
\ns\ds \G_i^*(s,t)=-\cR(t)^{-1} \int_t^T\cB(\th,t)^\top \(\Pi(s,t,\th)^\top{\bf1}_{(t_0,s)}(\th) \big[P^{(1)}_{1i}(s)^\top,P^{(1)}_{2i}(s)^\top,P^{(1)}_{3i}(s)^\top\big]^\top \\
\ns\ds\qq\qq\!+\!\int_\th^T\Pi(r,t,\th)^\top \big[P^{(2)}_{1i}(r,s,t)^\top,P^{(2)}_{2i}(r,s,t)^\top,
P^{(2)}_{3i}(r,s,t)^\top\big]^\top dr\)d\th,\ \  t_0\1n<\1nt\1n<\1ns\1n<\1nT.
\ea\ee

\no By \rf{K1-} and \rf{Ga_i}, we derive
\begin{eqnarray}
K_1^*(t)\1n=\1n-\cR(t)^{-1}\1n\Big\{\1nD_1(t)^\top \1n\cG_1(t)C_1(t)
\2n+\3n\int_t^T\3n\int_t^s\3n\cB(\th,t)^\top\Pi(s,t,\th)^\top P^{(1)}(s)^\top (I,0,0)^\top d\th ds \qq\qq\qq\notag\\
+\int_{(t+\d)\wedge T}^T\int_t^s\cB(\th,t)^\top\Pi(s,t,\th)^\top P^{(1)}(s)^\top (0,I,0)^\top d\th ds
+\int_t^T\int_t^s\cB(\th,t)^\top\Pi(s,t,\th)^\top P^{(1)}(s)^\top \notag\\
\qq \times (0,0,I)^\top \int_t^sF(s,\th')d\th' d\th ds
+\int_t^T\int_t^T\int_t^r\cB(\th,t)^\top\Pi(r,t,\th)^\top P^{(2)}(s,r,t)^\top (I,0,0)^\top d\th drds\ \  \notag\\
+\int_{(t+\d)\wedge T}^T\int_t^T\int_t^r\cB(\th,t)^\top\Pi(r,t,\th)^\top P^{(2)}(s,r,t)^\top (0,I,0)^\top d\th drds\qq\qq\qq\qq\qq\q\ \  \notag\\
+\int_t^T\2n\int_t^T\2n\int_t^r\2n\cB(\th,t)^\top \Pi(r,t,\th)^\top P^{(2)}(s,r,t)^\top (0,0,I)^\top\2n\int_t^sF(s,\th')d\th' d\th drds\Big\}.\qq\qq\qq\qq\!\notag
\end{eqnarray}
It then follows from some basic calculations that
$$\ba{ll}
\ns\ds
K_1^*(t)\1n=\1n-\cR(t)^{-1}\1n\Big\{\1nD_1(t)^\top \1n\cG_1(t)C_1(t)
+\int_t^T\int_t^s\cB(\th,t)^\top\Pi(s,t,\th)^\top \notag\\
\ns\ds \qq\qq\times P^{(1)}(s)^\top \Big(I,I{\bf1}_{((t+\d)\wedge T,\i)}(s),\int_t^sF(s,\th')^\top d\th'\Big)^\top d\th ds
+\2n\int_t^T\3n\int_t^T\3n\int_t^r\3n\cB(\th,t)^\top \ \notag\\
\ns\ds \qq\qq \times\Pi(r\1n,t,\th)^\top \1nP^{(2)}(s,r,t)^\top \Big(I,I{\bf1}_{((t+\d)\wedge T,\i)}(s),\int_t^s\3nF(s,\th')^\top d\th'\Big)^\top\1n d\th drds\Big\},\notag
\ea
$$
which then implies
$$\ba{ll}
\ns\ds K_1^*(t)\1n=\1n-\cR(t)^{-1}\1n\Big\{
D_1(t)^\top \1n\cG_1(t)C_1(t)+\3n\int_t^T\2n\cB(\th,t)^\top \Big[\int_\th^T\Pi(\a,t,\th)^\top P^{(1)}(\a)^\top\Upsilon(\a,t)d\a\\
\ns\ds \qq\qq +\int_t^T\int_\th^T \Pi(\a,t,\th)^\top P^{(2)}(r,\a,t)^\top\Upsilon(r,t)d\a dr\Big]d\th\Big\}.
\ea$$

Step 3: In this step, we turn to calculate $K_2^*$ and $K_3^*$.
From \rf{K2-} and \rf{Ga_i}, we deduce
$$\ba{ll}
\ns\ds  K_2^*(t,s)=-\cR(t)^{-1} \Big\{D_1(t)^\top\cG_1(t)C_3(t)F(t,s)
+\int_t^{s+\d}\cB(\th,t)^\top\Pi(s+\d,t,\th)^\top P^{(1)}(s+\d)^\top \notag\\
\ns\ds \qq\times (0,I,0)^\top{\bf1}_{[t-\d,T-\d]}(s)d\th +\int_t^T\3n\int_t^r\2n\cB(\th,t)^\top\Pi(r,t,\th)^\top \1nP^{(2)}(s\1n+\1n\d,r,t)^\top\1n(0,I,0)^\top\1n d\th dr{\bf1}_{[t-\d,T-\d]}(s)\notag\\
\ns\ds \qq+\3n\int_t^T\2n\int_t^\th \cB(\th',t)^\top \Pi(\th,t,\th')^\top P^{(1)}(\th)^\top(0,0,I)^\top F(\th,s)d\th'd\th \notag\\
\ns\ds \qq +\2n\int_t^T\3n\int_{t}^T\3n\int_t^r\1n\cB(\th',t)^\top \Pi(r,t,\th')^\top P^{(2)}(\th,r,t)^\top(0,0,I)^\top  F(\th,s)d\th'drd\th\Big\}.
\ea
$$
Applying Fubini theorem, we obtain
$$\ba{ll}
\ns\ds K_2^*(t,s)=-\cR(t)^{-1}\Big\{D_1(t)^\top \cG_1(t)C_3(t)F(t,s)
+\int_t^{T}\cB(\th,t)^\top \Big[\Big(\Pi(s+\d,t,\th)^\top P^{(1)}(s+\d)^\top {\bf1}_{[t,s+\d)} (\th)  \notag\\
\ns\ds \qq\q +\int_\th^T\1n\Pi(r,t,\th)^\top P^{(2)}(s+\d,r,t)^\top dr\Big)(0,I,0)^\top{\bf1}_{[t-\d,T-\d]}(s)
+ \int_\th^T\Pi(\a,t,\th)^\top P^{(1)}(\a)^\top(0,0,I)^\top F(\a\1n,s)d\a \\
\ns\ds\qq\q +\int_{t}^T\int_\th^T\Pi(\a,t,\th)^\top P^{(2)}(r,\a,t)^\top(0,0,I)^\top F(r,s)d\a dr\Big]d\th\Big\}.
\ea$$
Finally by \rf{Upsilon-cR}, \rf{K3-}, we derive %
$$\ba{ll}
\ns\ds K_3^*(t)=-\cR(t)^{-1} D_1(t)^\top\cG_1(t)C_2(t).
\ea$$

Step 4: In this step, we turn to look at the case of $K_4^*$.
From \rf{K4-}, we have
\bel{I1I2}\ba{ll}
\ns\ds K_4^*(t,s)=\1n\Big\{\2n\int_{s+\d}^T\2n \Big[{\G}_1^*(r,\1nt) \1n+\1n{\G}_2^*(r,\1nt) {\bf1}_{(s+2\d,\i)}(r) {\bf1}_{[t_0,T-\d)}(t) \2n\\
\ns\ds\qq\qq\q +\3n\int_{s+\d}^r\5n{\G}_3^*(r,t)F(r,\1n\th)d\th\Big]\1ndr B_2({s\1n+\1n\d})\Big\} {\bf1}_{[t-\d,T-\d]}(s)+\!\int_t^T \3n\int_{\th'}^T\3n\Big({\G}_1^*(r,t)+\!{\G}_2^*(r,t) \!\\
\ns\ds\qq\qq\q \times{\bf1}_{[t_0,r-\d)}(\th'){\bf1}_{[t_0,T-\d)}(t)+\!\!\int_{\th'}^r{\G}_3^*(r,t)F(r,\th)d\th\Big) B_3(\th')\ti F(\th',s)drd\th' \\
\ns\ds \qq\qq \equiv I_1(t,s+\d){\bf1}_{[t-\d,T-\d]}(s)+I_2(t,s),
\ea\ee
which implies that
$$
I_1(t,s)=\2n\int_{s}^T\2n\Big[{\G}_1^*(r,\1nt) \1n+\1n{\G}_2^*(r,\1nt) {\bf1}_{(s+\d,\i)}(r) {\bf1}_{[t_0,T-\d)}(t) \2n+\3n\int_{s}^r\3n{\G}_3^*(r,t)F(r,\1n\th)d\th\Big]\1ndr B_2({s}).
$$
By \rf{Ga_i}, we get
$$\ba{ll}
\ns\ds I_1(t,s)=-\cR(t)^{-1}\Big\{\int_s^T\int_t^r\cB(\th,t)^\top\Pi(r,t,\th)^\top P^{(1)}(r)^\top (I,0,0)^\top d\th dr+\int_{s+\d}^T\int_t^r\cB(\th,t)^\top\Pi(r,t,\th)^\top\\
\ns\ds \q\times P^{(1)}(r)^\top (0,I,0)^\top d\th dr{\bf1}_{[0,T-\d]}(t)
+\int_s^T\int_s^r\int_t^r\cB(\th',t)^\top\Pi(r,t,\th')^\top P^{(1)}(r)^\top (0,0,I)^\top F(r,\th) d\th'd\th dr \\
\ns\ds\q +\int_s^T\int_t^T\int_t^\a\cB(\th,t)^\top\Pi(\a,t,\th)^\top P^{(2)}(r,\a,t)^\top (I,0,0)^\top d\th d\a dr
+\3n\int_{s+\d}^T\1n\int_t^T\3n\1n\int_t^\a\3n\cB(\th,t)^\top\Pi(\a,t,\th)^\top\1n P^{(2)}(r,\a,t)^\top  \\
\ns\ds \q\times(0,I,0)^\top \1nd\th d\a dr{\bf1}_{[0,T-\d]}(t)+\3n\int_{s}^T\3n\1n\int_s^r\3n\1n\int_t^T\3n\1n\int_t^\a\3n\1n \cB(\th',t)^\top\Pi(\a,t,\th')^\top \1nP^{(2)}(r,\a,t)^\top \1n(0,0,I)^\top \1nF(r,\th)d\th' d\a d\th dr\1n\Big\}B_2(s).
\ea
$$
It then follows from some calculations that
\bel{I1}\ba{ll}
\ns\ds
I_1(t,s)=-\cR(t)^{-1}\Big\{\int_t^T\cB(\th,t)^\top \Big[\int_{s\vee\th}^T\Pi(r,t,\th)^\top
P^{(1)}(r)^\top\Big(I,{\bf1}_{(s+\d,\i)}(r) \\
\ns\ds \qq\qq\times {\bf1}_{[0,T-\d]}(t)I, \int_s^rF(r,\th')^\top d\th'\Big)^\top dr+\int_s^T\int_\th^T\Pi(\a,t,\th)^\top P^{(2)}(r,\a,t)^\top \\
\ns\ds \qq\qq \times \Big(I,{\bf1}_{(s+\d,\i)}(r){\bf1}_{[0,T-\d]}(t)I,\int_s^rF(r,\th')^\top d\th'\Big)^\top d\a dr\Big]d\th\Big\}B_2(s).
\ea\ee
Similarly, by \rf{Ga_i}, we obtain
\begin{eqnarray}
I_2(t,s)=-\cR(t)^{-1}\Big\{\int_t^T\cB(\th,t)^\top \Big\{\int_\th^T\Pi(r,t,\th)^\top P^{(1)}(r)^\top(\1n\int_t^r\3n\ti F(\th'\1n,s)^\top\1n B_3(\th')^\top \1n d\th',\1n\int_t^{r-\d}\3n\ti F(\th'\1n,s)^\top\1n B_3(\th')^\top d\th',\notag\\
 \int_t^r\3n\int_t^{\b }\3n\ti F(\th'\1n,s)^\top\1n B_3(\th')^\top\1n F(r,\b)^\top\1n d\th'd\b)^\top dr
+\int_t^T\int_\th^T\Pi(\a,t,\th)^\top P^{(2)}(r,\a,t)^\top (\int_t^r\ti F(\th',s)^\top B_3(\th')^\top d\th', \notag\\
\int_t^{r-\d}\ti F(\th',s)^\top B_3(\th')^\top d\th',\int_t^r\int_t^{\b }\ti F(\th',s)^\top B_3(\th')^\top F(r,\b)^\top d\th'd\b)^\top d\a dr\Big\}d\th\Big\}.\q\qq\qq\qq\qq\label{I2}
\end{eqnarray}
Substituting \rf{I1} and \rf{I2} into \rf{I1I2},  we deduce
$$\ba{ll}
\ns\ds K_4^*(t,s)=\1n-\cR(t)^{-1}\1n\Big\{\1n \int_t^T\1n\2n\cB(\th,t)^\top \Big[\int_{\th}^T\1n \Pi(\a,t,\th)^\top \1nP^{(1)}(\a)^\top\1n\Big\{\Big(I,{\bf1}_{(s+2\d,\i)}(\a)  {\bf1}_{(0,T-\d)}(t)I,\int_{s+\d}^\a F(\a,\th')^\top d\th'\Big)^\top \\
\ns\ds\qq\q \times{\bf1}_{(s+\d,T)}(\a)B_2(s+\d) {\bf1}_{[t-\d,T-\d]}(s)+\Big(\int_t^\a\2n\big(B_3(\th')\ti F(\th',s)\big)^\top d\th',\int_t^{\a-\d}\big(B_3(\th')\ti F(\th',s)\big)^\top d\th',\\
\ns\ds\qq\q \int_t^\a\2n\int_t^{\b}\3n\big(F(\a,\b)B_3(\th')\ti F(\th',s)\big)^\top d\th'd\b\Big)^\top\Big\} d\a +\int_{t}^T\int_\th^T\Pi(\a,t,\th)^\top P^{(2)}(r,\a,t)^\top
\ea$$
$$\ba{ll}
\ns\ds\qq\q \times \Big\{\Big(I,{\bf1}_{(s+2\d,\i)}(r){\bf1}_{(0,T-\d)}(t)I, \int_{s+\d}^r F(r,\th')^\top d\th'\Big)^\top{\bf1}_{(s+\d,T)}(r)
B_2(s+\d){\bf1}_{[t-\d,T-\d]}(s)\1n \notag\q\\
\ns\ds\qq\q +\Big(\int_t^r\2n\big(B_3(\th')\ti F(\th',s)\big)^\top\2n d\th'\1n\1n,\int_t^{r-\d}\5n\3n
\big(B_3(\th')\ti F(\th',s)\big)^\top d\th'\1n\1n,\int_t^r\1n\int_t^{\b}\1n\big(F(r,\b)B_3(\th')\ti F(\th',s)\big)^\top d\th'd\b)^\top\1n \Big\}d\a dr\Big]d\th\Big\}.
\ea
$$

Step 5: In this step, we calculate $v^*$. Recalling \rf{Ga_i}, we have
%
$$\ba{ll}
\ns\ds \int_{t-\d}^{t_0}\Big\{ \int_{s+\d}^T\Big[{\G}_1^*(r,t) +{\G}_2^*(r,t) {\bf1}_{(s+2\d,\i)}(r) {\bf1}_{[0,T-\d)}(t) \\
\ns\ds\qq +\int_{s+\d}^r{\G}_3^*(r,t)F(r,\th)d\th\Big]dr B_2(s+\d)\Big\} \varsigma(s-t_0)ds{\bf1}_{[t_0,t_0+\d]}(t)\notag\\
\ns\ds =
-\cR(t)^{-1}\int_t^T\cB(\th,t)^\top \Big[\int_{\th}^T\Pi(r,t,\th)^\top P^{(1)}(r)^\top \int_{t-\d}^{t_0}\2n{\bf1}_{(s+\d,\i)}(r)\\
\ns\ds \q\times\Big(I, {\bf1}_{(s+2\d,\i)}(r){\bf1}_{(0,T-\d)}(t)I, \int_{s+\d}^rF(r,\th')^\top d\th'\Big)^\top B_2(s+\d)\varsigma(s-t_0)ds{\bf1}_{[t_0,t_0+\d]}(t) dr\notag\
\\
\ns\ds \q+\2n\int_{t}^T\int_\th^T\2n\Pi(\a,t,\th)^\top \2nP^{(2)}(r,\a,t)^\top\int_{t-\d}^{t_0} {\bf1}_{(s+\d,\i)}(r)\Big(I,{\bf1}_{(s+2\d,\i)}(r) {\bf1}_{(0,T-\d)}(t)I, \\
\ns\ds \q\int_{s+\d}^r\2n F(r,\th')^\top\2n d\th'\Big)^\top B_2(s\1n+\1n\d)\varsigma(s-t_0)ds {\bf1}_{[t_0,t_0+\d]}(t) d\a dr\Big]d\th.
\ea$$
Similarly, we obtain
\bel{v*-2}\ba{ll}
\ns\ds \int_t^T\3n{\Gamma}^*_2(s,t) \xi(s-t_0-\d){\bf1}_{[t_0,t_0+\d]}(s)ds  \\
\ns\ds =-\cR(t)^{-1}\int_t^T\(\int_t^s \cB(\th,t)^\top \Pi(s,t,\th)^\top \1n d\th P^{(1)}\1n(s)^\top(0,I,0)^\top\\
\ns\ds \qq\3n\2n+\3n\int_t^T\3n\1n\int_t^r\3n \cB(\th,t)^\top \Pi(r,t,\th)^\top P^{(2)}(r,s,t)(0,I,0)^\top\1n d\th dr \) \xi(s\1n-\1n\d\1n-t_0){\bf1}_{[t_0,t_0+\d]}(s)ds.
\ea\ee
Finally, by the Fubini theorem,
\begin{eqnarray}
\3n\int_t^T\3n\int_t^T\3n\int_t^r\cB(\th,t)^\top \Pi(r,t,\th)^\top P^{(2)}(r,s,t)(0,I,0)^\top d\th dr \xi(s-\d-\1nt_0){\bf1}_{[t_0,t_0+\d]}(s)ds\notag\\
=\int_t^T\3n\int_t^T\3n\int_\th^T\cB(\th,t)^\top\Pi(\a,t,\th)^\top P^{(2)}(\a,r,t)(0,I,0)^\top d\a  d\th \xi(s-\d-\1nt_0){\bf1}_{[t_0,t_0+\d]}(s)ds,\notag
\end{eqnarray}
which and \rf{check-o-decompose}, \rf{v-}, \rf{v*-2} imply that
\begin{eqnarray}
v^*(t)
\1n=\1n-\cR(t)^{-1}\1n\Big\{\1nD_1(t)^\top\1n \cG_1(t)\sigma(t) \2n+\1n\2n\int_t^T\2n\1n\1nD_1(t)^\top \1n\Upsilon(\a,t)^\top\1n\zeta(\a,t)d\a\1n+\1n\2n\int_t^T \1n\2n\1n\int_\th^T\1n\2n\cB(\th,t)^\top\Pi(\a,t,\th)^\top\1n\eta(\a,t)d\a d\th
\1n\1n\qq\q\notag\\
+\2n\int_t^T\1n\2n\cB(\th,t)^\top\Big[\1n\int_{\th}^T\3n\2n \Pi(\a,t,\th)^\top\1n P^{(1)}\1n(\a)^\top\1n\Big(\1n(0,I,0)^\top \1n \xi(\a\1n-\1n\d-t_0){\bf1}_{[t_0,t_0+\d]}(\a)
\1n\1n+\1n\2n\int_{t-\d}^{t_0}\1n \3n{\bf1}_{(\th'+\d,\i)}\1n(\a) (I,{\bf1}_{(\th'+2\d,\i)}\1n(\a)\notag\\
\times{\bf1}_{(0,T-\d)}(t)I,\1n \int_{\th'+\d}^\a\3n\1nF(\a\1n,\b)^\top\1n d\b)^\top \1n\1nB_2(\th'\1n\1n+\1n\d) \varsigma(\th'\1n\1n-\1nt_0)d\th' {\bf1}_{ [t_0,t_0+\d]}(t)\1n\Big)d\a
\1n+\3n\int_t^T\1n\3n\int_\th^T\3n\Pi(\a,t, \th)^\top \1n P^{(2)}(r,\a,t)^\top \q\ \ \notag\\
\times\1n\Big(\1n(0,I,0)^\top
\1n\xi(r\1n-\1n\d\1n-\1nt_0){\bf1}_{[t_0,t_0+\d]}(r) \1n+\1n\1n\int_{t-\d}^{t_0}\1n{\bf1}_{(\th'+\d,\i)}(r) (I,{\bf1}_{(\th'+2\d,\i)}(r) {\bf1}_{(0,T-\d)}(t)I, \qq\qq\qq\qq\  \notag\\
\int_{\th'+\d}^r\2n F(r,\b)^\top\2n d\b)^\top B_2(\th'\1n+\1n\d)\varsigma(\th'-t_0)d\th' {\bf1}_{[t_0,t_0+\d]}(t) \Big)d\a dr\Big]d\th\Big\}.\ \qq\qq\qq\qq\qq\qq\qq\qq
\notag
\end{eqnarray}

Step 6: In this step, we show that $(K_1^*(\cd),K_2^*(\cd,\cd),K_3^*(\cd)$, $K_4^*(\cd,\cd), v^*(\cd))$ is the optimal closed-loop strategy of Problem {\rm (P)} on $[t_0,T]$ in terms of Definition \ref{def optimal closed-loop strategy}.

In fact, by the optimality in Step 1, it is sufficient to prove that $(K_1^*(\cd),K_2^*(\cd,\cd),K_3^*(\cd)$, $K_4^*(\cd,\cd), v^*(\cd))\in\dbL$.
Next we prove that $K_1(\cd)\in L^2(t_0,T;\dbR^{n\times m})$. Since $(P^{(1)},P^{(2)})\in\Pi(0,T)$, we have
\bel{P-integrability}\ba{ll}
\ns\ds
\operatorname*{ess\,sup}\limits_{t\in(0,T)}|P^{(1)}(t)| %
+\big(\int_0^T\int_0^T\sup\limits_{t\in[0,s_1\wedge s_2]}|P^{(2)}(s_1,s_2,t)|^2 ds_1ds_2\big)^{\frac 1 2}<\i.
\ea\ee
By the boundedness of $F(\cd,\cd)$, we obtain
\bel{P1-integrability}\ba{ll}
\ns\ds \sup\limits_{t\in(0,T)}\int_t^T|\Upsilon(s,t)^\top P^{(1)}(s)\Upsilon(s,t)|ds 
\les M\operatorname*{ess\,sup}\limits_{t\in(0,T)} |P^{(1)}(t)| \int_0^T\Big|1+\int_0^s|F(s,r)|dr\Big|^2ds<\i,
\ea\ee
and
\bel{P2-integrability}\ba{ll}
\ns\ds \sup\limits_{t\in(0,T)}\int_t^T\int_t^T|\Upsilon(s_1,t)^\top P^{(2)}(s_1,s_2,t)\Upsilon(s_2,t)|ds_1ds_2\\
\ns\ds \q\les 
M\big(\int_0^T\int_0^T\sup\limits_{t\in[0,s_1\wedge s_2]}|P^{(2)}(s_1,s_2,t)|^2 ds_1ds_2\big)^{\frac 1 2}<\i,
\ea\ee
where $M$ is a generic constant. Thus \rf{P1-integrability} and \rf{P2-integrability} imply the boundedness of $\cG_1(\cd)$. Recall Theorem \ref{solvability of Riccati}, $\cR(\cd)\ges \b I$ for some constant $\b>0$. Then, the boundedness of $R_1(\cd)$, $R_2(\cd)$, $D_1(\cd)$ imply that $\cR(\cd)^{-1}$ is bounded.  Notice that
\bel{K1-term}\ba{ll}
\ns\ds \int_{t_0}^T\Big|\int_t^T\int_t^T\cB(\th,t)^\top\int_\th^T\Big( P^{(2)}(\a,r,t)\Pi(r,t,\th)\Big)^\top dr\Upsilon(\a,t)d\a d\th\Big|^2dt\\
\ns\ds \les \int_{t_0}^T\Big|\int_t^T|\cB(\th,t)|\int_{\th}^T\int_t^T \Big(\frac{1}{r-t}+1+\frac{1}{r-t} |\cE(r,t)|+\frac{1}{r-t} |\cE(r,t+\d)|+|\cE(r,\th)|\Big)\\
\ns\ds \qq
\times|P^{(2)}(\a,r,t)|\Big(1+ |\cE(\a,t)|\Big)d\a drd\th\Big|^2dt.
\ea\ee
By the boundedness of $F(\cd,\cd)$, $B_1(\cd)$, $B_2(\cd)$ and $B_3(\cd)$, we deduce
\bel{K1-1-term}\ba{ll}
\ns\ds \int_{t_0}^T\Big|\int_t^T|\cB(\th,t)|\int_\th^T\int_t^T \frac{1}{r-t}|P^{(2)}(\a,r,t)||\cE(\a,t)|d\a drd\th\Big|^2dt\\
\ns\ds =\int_{t_0}^T\Big|\int_t^T\int_t^r|\cB(\th,t)|\int_t^T \frac{1}{r-t}|P^{(2)}(\a,r,t)||\cE(\a,t)|d\a d\th dr\Big|^2dt\\
\ns\ds \les \int_{t_0}^T\int_t^T\Big(\int_t^r|\cB(\th,t)|d\th\frac{1}{r-t} \Big)^2dr\int_t^T\Big(\int_t^T|P^{(2)}(\a,r,t)\cE(\a,t)| d\a\Big)^2drdt\\
\ns\ds \les M\int_0^T\int_0^T\sup\limits_{t\in[0,\a\wedge r]}|P^{(2)}(\a,r,t)|^2 d\a dr<\i,
\ea\ee
and
\bel{K1-2-term}\ba{ll}
\ns\ds \int_{t_0}^T\Big|\int_t^T|\cB(\th,t)|\int_\th^T\int_t^T \frac{1}{r-t}|\cE(r,t)||P^{(2)}(\a,r,t)||\cE(\a,t)|d\a drd\th\Big|^2dt\\
\ns\ds =\int_{t_0}^T\Big|\int_t^T\int_t^r|\cB(\th,t)|\int_t^T \frac{1}{r-t}|\cE(r,t)||P^{(2)}(\a,r,t)||\cE(\a,t)|d\a d\th dr\Big|^2dt\\
\ns\ds \les M\int_0^T\int_0^T\sup\limits_{t\in[0,\a\wedge r]}|P^{(2)}(\a,r,t)|^2 d\a dr<\i.
\ea\ee
Similar to \rf{K1-1-term} and \rf{K1-2-term}, we can deal with other terms in \rf{K1-term} and prove that
$K_1(\cd)\in L^2(t_0,T;\dbR^{n\times m})$. Furthermore, we can verify that $(K_1^*(\cd),K_2^*(\cd,\cd),K_3^*(\cd)$, $K_4^*(\cd,\cd), v^*(\cd))\in\dbL$.
\epf

\ss

Before proving Corollary \ref{optimal control with control delay} and Corollary \ref{optimal control with state delay}, we firstly decompose the Riccati system \rf{eq_Riccati--Volterra} and the Type-II extended backward SVIE \rf{EBSVIE}.

Let $P=(P^{(1)},P^{(2)})$ be the solution to the Riccati--Volterra equation \eqref{eq_Riccati--Volterra} and decompose them as follows:
$$\ba{ll}
P^{(1)}(t)=
\begin{bmatrix}
\ns\ds P^{(1)}_{11}(t)\q P^{(1)}_{12}(t)\q P^{(1)}_{13}(t)\\
\ns\ds P^{(1)}_{21}(t)\q P^{(1)}_{22}(t)\q P^{(1)}_{23}(t)\\
\ns\ds P^{(1)}_{31}(t)\q P^{(1)}_{32}(t)\q P^{(1)}_{33}(t)
\end{bmatrix},\q
P^{(2)}(s,t,r)=
\begin{bmatrix}
\ns\ds P^{(2)}_{11}(s,t,r)\q P^{(2)}_{12}(s,t,r)\q P^{(2)}_{13}(s,t,r)\\
\ns\ds P^{(2)}_{21}(s,t,r)\q P^{(2)}_{22}(s,t,r)\q P^{(2)}_{23}(s,t,r)\\
\ns\ds P^{(2)}_{31}(s,t,r)\q P^{(2)}_{32}(s,t,r)\q P^{(2)}_{33}(s,t,r)
\end{bmatrix}.
\ea$$
Recall that \rf{eq_Riccati--Volterra} is equivalent to \rf{eq_Riccati--Volterra-0}. Then, we decompose its coefficients as follows:
\vskip-6mm
\begin{eqnarray}
(P\rtimes A)(s,t)
=\begin{bmatrix}
\cT_{PA}^{(11)}(s,t)& \cT_{PA}^{(12)}(s,t)& \cT_{PA}^{(13)}(s,t)\\
\cT_{PA}^{(21)}(s,t)& \cT_{PA}^{(22)}(s,t)& \cT_{PA}^{(23)}(s,t)\\
\cT_{PA}^{(31)}(s,t)& \cT_{PA}^{(32)}(s,t)& \cT_{PA}^{(33)}(s,t)
\end{bmatrix},\q
(P\rtimes B)(s,t)
=\begin{bmatrix}\cT_{PB}^{(1)}(s,t)\\
\cT_{PB}^{(2)}(s,t)\\
\cT_{PB}^{(3)}(s,t)
\end{bmatrix}.\notag
\end{eqnarray}
In this case, we have
\begin{eqnarray}
\cT_{PA}^{(ij)}(s,t)=\Big[P_{i1}^{(1)}(s)+P_{i2}^{(1)}(s)
{\bf1}_{(t+\d,\i)}(s)
+P_{i3}^{(1)}(s)\cE(s,t)\Big]A_j(t)\qq\qq\qq\qq\qq\qq\qq\notag\\
+\3n\int_t^T\3n\Big[P_{i1}^{(2)}(s,r,t)
\1n +\1nP_{i2}^{(2)}(s,r,t){\bf1}_{(t+\d,\i)}(r) \1n+\1nP_{i3}^{(2)}(s,r,t)\cE(r,t)
\Big]drA_j(t),\ 0\1n<\1nt\1n<\1ns\1n<\1nT,i,j\1n=\1n1,2,3,\label{PA}
%
\end{eqnarray}
and
\begin{eqnarray}
\cT_{PB}^{(i)}(s,t)=\Big\{P_{i1}^{(1)}(s)+P_{i2}^{(1)}(s){\bf1}_{[0,s-\d)}(t)
+P_{i3}^{(1)}(s)\cE(s,t)+\int_t^T\Big[P_{i1}^{(2)}(s,r,t)
+P_{i2}^{(2)}(s,r,t){\bf1}_{[0,r-\d)}(t)\notag\\
+P_{i3}^{(2)}(s,r,t)\cE(r,t)\Big]dr
\Big\}B_1(t)
+\Big\{P_{i1}^{(1)}(s){\bf1}_{[0,s-\d)}(t)+P_{i2}^{(1)}(s){\bf1}_{[0,s-2\d)}(t)
+P_{i3}^{(1)}(s)\cE(s,t+\d)\notag\\
+\int_t^T\Big[P_{i1}^{(2)}(s,r,t)
{\bf1}_{[0,r-\d)}(t)+P_{i2}^{(2)}(s,r,t) {\bf1}_{[0,r-2\d)}(t) +P_{i3}^{(2)}(s,r,t)\cE(r,t+\d)\Big]dr
\Big\}B_2(t+\d)\ \1n\notag\\
+\int_t^s\Big[P_{i1}^{(1)}(s)+P_{i2}^{(1)}(s)
{\bf1}_{[0,s-\d)}(\th)+P_{i3}^{(1)}(s)\cE(s,\th)\Big] B_3(\th)\ti F(\th,t)d\th+\int_t^T\int_t^r\Big[P_{i1}^{(2)}(s,r,t)
\q \notag\\
+P_{i2}^{(2)} (s,r,t)
{\bf1}_{[0,r-\d)}(\th)+P_{i3}^{(2)}(s,r,t)\cE(r,\th)\Big] B_3(\th)\ti F(\th,t)d\th dr, i=1,2,3.\label{PB}\qq\qq\qq\qq\
\end{eqnarray}
%
%
Hence, by \rf{CPC}--\rf{DPC}, we derive the following decomposed coupled Riccati equation:
\bel{coupled Riccati}\left\{\ba{ll}
\ns\ds P_{ii}^{(1)}(t)=Q_i(t)+C_i(t)^\top\cG_1(t)C_i(t) -C_i(t)^\top\cG_1(t)D_1(t)
\cR(t)^{-1} D_1(t)^\top\cG_1(t)C_i(t),\\
 \ns\ds \qq\qq\qq\qq\qq\qq\qq\qq\qq 0<t<T,\ i=1,2,3,\\
\ns\ds P_{ij}^{(1)}(t)=C_i(t)^\top \cG_1(t)C_j(t)-C_i(t)^\top\cG_1(t)D_1(t)
\cR(t)^{-1} D_1(t)^\top\cG_1(t)C_j(t),\\
\ns\ds \qq\qq\qq\qq\qq\qq\qq 0<t<T,\q i\neq j,\q i,j=1,2,3,\\
\ns\ds \dot{P}_{ij}^{(2)}(s_1,s_2,t) =\cT_{PB}^{(i)}(s_1,t)\cR(t)^{-1} \cT_{PB}^{(j)}(s_2,t)^\top,\q 0<t<(s_1\wedge s_2)<T,\q i,j=1,2,3,\\
\ns\ds P_{ij}^{(2)}(s,t,t)=P_{ji}^{(2)}(t,s,t)^\top
=\cT_{PA}^{(ij)}(s,t)
-\cT_{PB}^{(i)}(s,t)\cR(t)^{-1} D_1(t)^\top\cG_1(t)C_j(t),\q 0<t<s<T,
\ea\right.\ee
where $\cR(\cd)$, $\cG_1(\cd)$, $\cT_{PA}^{(ij)}(\cd,\cd)$ and $\cT_{PB}^{(i)}(\cd,\cd)$ are defined by \rf{Upsilon-cR}, \rf{PA}, \rf{PB}, and the decomposed $\cG_1(\cd)$ has a specific representation in the subsequent special cases, which is omitted here.
 Similarly one can decompose the Type-II extended backward SVIE \rf{EBSVIE}. We will not go into the details.

\ms

{\bf The proof of Corollary \ref{optimal control with control delay}:}
\bpf
In this special case, by \rf{coupled Riccati} we get
\begin{eqnarray}
P_{11}^{(1)}(t)=Q_1(t)+C_1(t)^\top\cG_1(t)C_1(t) -C_1(t)^\top\cG_1(t)D_1(t)
\cR(t)^{-1} D_1(t)^\top\cG_1(t)C_1(t),\ 0<t<T,\qq\qq\q\  \label{P11-1-control delay}\notag
%
%
\end{eqnarray}
and
\begin{eqnarray}
P_{11}^{(2)}(s,t,t)=\cT_{PA}^{(11)}(s,t)
-\cT_{PB}^{(1)}(s,t)\cR(t)^{-1} D_1(t)^\top \cG_1(t)C_1(t),\q 0<t<s<T,\qq\qq\qq\q \q \label{P11-2-control delay-1}\\
\ns\ds \dot{P}_{11}^{(2)}(s_1,s_2,t)=\cT_{PB}^{(1)}(s_1,t) \cR(t)^{-1} \cT_{PB}^{(1)}(s_2,t)^\top,\ 0<t<(s_1\wedge s_2)<T,\ \ \qq\qq\qq\qq\qq\q \ \label{P11-2-control delay-2}\\
P_{ij}^{(1)}(\cd), P_{ij}^{(2)}(\cd,\cd,\cd)=0, \q (i,j)=(1,2),(1,3),(2,1),(2,2),(2,3),(3,1),(3,2),(3,3), \q 0<t<T,\q\notag
\end{eqnarray}
where $\cT_{PB}^{(1)}(\cd,\cd)$ is given as follows:
$$\ba{ll}
\ns\ds
\cT_{PB}^{(1)}(s,t)=\Big\{P_{11}^{(1)}(s)+\int_t^T P_{11}^{(2)}(s,r,t) dr
\Big\}B_1(t)
+\Big\{P_{11}^{(1)}(s){\bf1}_{[0,s-\d)}(t)+\int_t^T P_{11}^{(2)}(s,r,t)
{\bf1}_{[0,r-\d)}(t)dr
\Big\} \\
\ns\ds \qq\times B_2(t+\d)+P_{11}^{(1)}(s) \int_t^s B_3(\th)\ti F(\th,t)d\th\ \1n +\int_t^T P_{11}^{(2)}(s,r,t)\Big[\int_t^r
 B_3(\th)\ti F(\th,t)d\th\Big] dr, i=1,2,3.
\ea
$$
In addition,
$$\ba{ll}
\ns\ds \cG_1(t)=\int_t^TP_{11}^{(1)}(s)ds +\int_t^T\int_t^TP_{11}^{(2)}
(s_1,s_2,t)ds_1ds_2=\cP_1(t,t,t)=\cS_0(t), \\
\ns\ds \cR(t)=R_1(t) +D_1(t)^\top \cG_1(t)D_1(t), \\
\ns\ds \cT_{PA}^{(11)}(s,t)=\Big[P_{11}^{(1)}(s)
+\int_t^T P_{11}^{(2)}(s,r,t) dr\Big]A_1(t).
\ea
$$
From \rf{P-control delay} and \rf{cS1}, for $r\in[t-\d,t]$, we have
\begin{eqnarray}
\cS_1(t,r-t)=B_2(r+\d)^\top\cP_1(t,t,r+\d)+\int_t^T\ti F(\th',r)^\top B_3(\th')^\top \cP_1(t,t,\th')d\th'.
\notag%
\end{eqnarray}
Notice that for $t\les \th'$, we have
$$\ba{ll}
\ns\ds \cP_1(t,t,\th')=\int_{\th'}^TP^{(1)}_{11}(s)ds+\int_{\th'}^T\int_{t}^T P^{(2)}_{11}(s,\a,t)d\a ds.
\ea
$$
Therefore,

$$\ba{ll}
\ns\ds \frac{\partial}{\partial t}\cP_1(t,t,\th')=-\int_{\th'}^T P^{(2)}_{11}(s,t,t)ds-\int_{\th'}^T
\int_t^T \dot{P}_{11}^{(2)}(s,\a,t) d\a ds\\
\ns\ds =-\int_{\th'}^T P^{(2)}_{11}(s,t,t)ds-\Big[\int_{\th'}^T \cT_{PB}^{(1)}(s,t) ds\Big]
\cR(t)^{-1} \Big[\int_t^T \cT_{PB}^{(1)}(\a,t)^\top d\a\Big].
\ea
$$
It then follows from some calculations that
\begin{eqnarray}
\frac{\partial}{\partial t}\cP_1(t,t,\th')=-\int_{\th'}^TP_{11}^{(1)}(s)dsA_1(t) -\int_{\th'}^T\int_t^TP_{11}^{(2)} (s,r,t)drdsA_1(t)\qq \qq\qq\qq\qq\qq\qq\qq\qq \notag\\
+\Big(\cP_1(t,t,\th')B_1(t)+\cP_1(t,t+\d,\th')B_2(t+\d) +\int_t^T\cP_1(t,\th,\th')B_3(\th)\ti F(\th,t)d\th\Big) \qq\qq\qq\qq \qq\qq\notag\\
\times\cR(t)^{-1} D_1(t)^\top\cG_1(t)C_1(t) \2n+\2n\Big(\cP_1(t,t,\th')B_1(t)+\cP_1(t,t+\d,\th') B_2(t+\d)+\int_t^T\cP_1(t,\th,\th')B_3(\th)\ti F(\th,t)d\th\Big)\q \notag\\
\times\cR(t)^{-1}\Big(\cP_1(t,t,t)B_1(t)
\2n+\1n\cP_1(t,t+\d,t) B_2(t+\d)+\int_t^T\cP_1(t,\th,t)B_3(\th)\ti F(\th,t)d\th \Big)^\top,\ \ae\ t.\ \ \!\qq\qq\qq\
\notag
\end{eqnarray}
Then, we deduce
\begin{eqnarray}
\frac {\partial} {\partial t}\cS_1(t,r-t)=
-B_2(r+\d)^\top\int_{r+\d}^TP_{11}^{(1)}(s)dsA_1(t)-B_2(r+\d)^\top \int_{r+\d}^T\int_t^TP_{11}^{(2)}(s,\t,t)d\t dsA_1(t)\qq\qq\notag\\
+B_2(r+\d)^\top\Big( \cP_1(t,t,r+\d)B_1(t)+\cP_1(t,t+\d,r+\d)B_2(t+\d)+\int_t^T\cP_1(t,\th,r+\d) B_3(\th)\ti F(\th,t)d\th\Big)\q\notag\\
\times\cR(t)^{-1} D_1(t)^\top\cG_1(t)C_1(t)
+B_2(r+\d)^\top\Big(\cP_1(t,t,r+\d)B_1(t) +\cP_1(t,t+\d,r+\d)B_2(t+\d) \qq\qq\q\!\ \notag\\
+\2n\int_t^T\cP_1(t,\th,r+\d)B_3(\th)\ti F(\th,t)d\th\Big)\cR(t)^{-1}\Big(\cP_1(t,t,t)B_1(t)+\cP_1(t,t+\d,t)B_2(t+\d) \qq\qq\qq\q\ \ \ \notag\\
+\1n\int_t^T\3n\cP_1(t,\th,t)B_3(\th)\ti F(\th,t)d\th\Big)^\top\3n
-\1n\ti F(t,r)^\top\1n B_3(t)^\top\1n \cP_1(t,t,t)\1n-\3n\int_t^T\1n\ti F(\th',r)^\top \1n B_3(\th')^\top\3n\int_{\th'}^T\3nP_{11}^{(1)}(s)dsA_1(t)d\th'\q\ \notag
\end{eqnarray}
\begin{eqnarray}
-\int_t^T\ti F(\th',r)^\top B_3(\th')^\top\int_{\th'}^T\int_t^TP_{11}^{(2)}(s,\t,t)d\t dsd\th'A_1(t)
+\int_t^T\ti F(\th',r)^\top B_3(\th')^\top\Big(\cP_1(t,t,\th')B_1(t)\q\ \ \ \notag\\
+\cP_1(t,t+\d,\th')B_2(t+\d) +\int_t^T\cP_1(t,\th,\th')B_3(\th)\ti F(\th,t)d\th\Big)\cR(t)^{-1} D_1(t)^\top \cG_1(t)C_1(t)d\th'\qq\qq\qq\ \ \ \notag\\
+\int_t^T\ti F(\th',r)^\top B_3(\th')^\top \Big(\cP_1(t,t,\th')B_1(t)+\cP_1(t,t+\d,\th')B_2(t+\d)
+\int_t^T\cP_1(t, \th,\th')B_3(\th)\ti F(\th,t)d\th\Big)\ \ \notag\\
\times\cR(t)^{-1}\Big(\cP_1(t,t,t)B_1(t)+\cP_1(t,t+\d,t) B_2(t+\d)+\int_t^T\cP_1(t, \th,t)B_3(\th)\ti F(\th,t)d\th\Big)^\top d\th'.\ \ae\ t.\q\qq\
\notag
\end{eqnarray}
It yields
\begin{eqnarray}
\frac{\partial}{\partial t}\cS_1(t,r-t) +\ti F(t,r)^\top B_3(t)^\top\cS_0(t)
+\cS_1(t,r-t)A_1(t)-\big[\cS_1(t,r-t)B_1(t) \q\qq\qq\qq\q\notag\\
+\cS_2(t,r-t,0)\big] \cR(t)^{-1}\big[B_1(t)^\top \cS_0(t)+\cS_1(t,0)+D_1(t)^\top\cS_0(t)C_1(t)\big]=0,\q \ae\ t\in[t_0,T],
\notag
\end{eqnarray}
which implies that \rf{dcS1} holds. The proofs for \rf{dcS0} and \rf{dcS2} are similar.

Recalling \rf{Ki}--\rf{transformed optimal control}, in this case the optimal closed-loop outcome control is as follows:
\begin{eqnarray}
u^*(t)=-\cR(t)^{-1}\Big\{\Big(D_1(t)^\top\cS_0(t)C_1(t)+B_1(t)^\top \cS_0(t)+B_2(t+\d)^\top\cP_1(t,t+\d,t)^\top \qq\qq\qq\qq\qq\q \notag\\
+\int_t^T\ti F(\th,t)^\top B_3(\th)^\top\cP_1(t,\th,t)^\top d\th\Big)x^*(t)+\int_{(t-\d)\vee t_0}^{t\wedge(T-\d)}\Big[B_1(t)^\top \cP_1(t,t,r+\d)^\top\qq\qq\qq\qq\q \notag\\
+B_2(t+\d)^\top\cP_1(t,t+\d,r+\d)^\top+\int_t^T\ti F(\th,t)^\top B_3(\th)^\top\cP_1(t,\th,r+\d)^\top d\th\Big]B_2(r+\d)u^*(r)dr\ \q\ \qq\notag\\
+\int_{t_0}^t\int_t^T\Big[B_1(t)^\top\cP_1(t,t,\th')+B_2(t+\d)^\top \cP_1(t,t+\d,\th')^\top+\int_t^T\ti F(\th,t)^\top B_3(\th)^\top\cP_1(t,\th,\th')^\top d\th\Big]\q\ \notag\\
\times B_3(\th')\ti F(\th',r)u^*(r)d\th' dr\Big\}\notag.\qq\qq\qq\qq\qq\qq\qq\qq\qq\qq\qq\qq\qq \qq\q\ \ \
\notag
\end{eqnarray}
By the definitions of \rf{cS0}--\rf{cS2}, we deduce \rf{optimal-control-with-control-delay} and thus complete the proof of Corollary \ref{optimal control with control delay}.
\epf

{\bf The proof of Corollary \ref{optimal control with state delay}:}
\bpf
In this case, by \rf{coupled Riccati} we get
\begin{eqnarray}
P_{i3}^{(j)}(t),
P_{3i}^{(j)}(t),P_{3i}^{(j)}(t)=0,\ i=1,2,3,\ j=1,2, \ t\in(0,T), \qq\qq\qq\qq\q\notag
\end{eqnarray}
and other $P^{(1)}_{ij}$ and $P^{(2)}_{ij}$ satisfy \rf{coupled Riccati}, $i,j=1,2$.
By the definition \rf{cP3} of $\cP_3(\cd,\cd)$, we have
\begin{eqnarray}
\frac{\partial\cP_3(t,\th)}{\partial t}=-P_{21}^{(2)}(\th,t,t)^\top-P_{22}^{(2)}(\th,t+\d,t)^\top +\int_t^T\Big[P_{11}^{(1)}(r)+P_{12}^{(1)}(r) {\bf1}_{[0,r-\d)}(t)\qq\qq\qq\notag\\
+\int_t^T\Big(P_{11}^{(2)}(r,\a,t) +P_{12}^{(2)}(r,\a,t){\bf1}_{[0,\a-\d)}(t)\Big)d\a \Big]dr B_1(R_1+D_1^\top\cP_2(t)D_1)^{-1} B_1^\top\Big[P_{21}^{(1)}(\th)^\top\!\notag\\
+P_{22}^{(1)}(\th)^\top {\bf1}_{[0,\th-\d)}(t)+\int_t^T\Big(P_{21}^{(2)}(\th,\a,t) ^\top+P_{22}^{(2)}(\th,\a,t)^\top{\bf1}_{[0,\a-\d)}(t)\Big) d\a\Big]\qq\qq\qq\ \ \!\notag\\
+\int_{t+\d}^T\Big[P_{21}^{(1)}(r)+P_{22}^{(1)}(r) {\bf1}_{[0,r-\d)}(t)+\int_t^T\Big(P_{21}^{(2)}(r,\a,t) +P_{22}^{(2)}(r,\a,t){\bf1}_{[0,\a-\d)}(t)\Big)d\a\Big]dr  \ \ \notag\\
\times B_1(R_1+D_1^\top\cP_2(t)D_1)^{-1} B_1^\top\Big[P_{21}^{(1)}(\th)^\top+P_{22}^{(1)}(\th)^\top {\bf1}_{[0,\th-\d)}(t)+\int_t^T\Big(P_{21}^{(2)} (\th,\a,t)^\top\qq\ \ \notag\\
+P_{22}^{(2)}(\th,\a,t)^\top {\bf1}_{[0,\a-\d)}(t) \Big)d\a\Big].\qq\qq\qq\qq\qq\qq\qq\qq\qq\qq\qq\qq\q
\notag
\end{eqnarray}
Notice that
\begin{eqnarray}
-P_{21}^{(2)}(\th,t,t)^\top=-A_1^\top\Big[ P_{21}^{(1)}(\th)^\top+P_{22}^{(2)}(\th)^\top {\bf1}_{[0,\th-\d)}(t)\Big]-\int_t^TA_1^\top\Big[ P_{21}^{(2)}(\th,\a,t)^\top\q \notag\\
+P_{22}^{(2)}(\th,\a,t)^\top {\bf1}_{[0,\a-\d)}(t)\Big]d\a +C_1^\top\cP_2(t)D_1(R_1+D_1^\top\cP_2(t)D_1)^{-1} B_1^\top\Big[P_{21}^{(1)}(\th)^\top\3n\2n \notag\\
+P_{22}^{(1)}(\th)^\top {\bf1}_{[0,\th-\d)}(t) +\int_t^T\Big(P_{21}^{(2)}(\th,\a,t)^\top +P_{22}^{(2)}(\th,\a,t)^\top{\bf1}_{[0,\a-\d)}(t)\Big) d\a\Big],\ \ \
\notag
\end{eqnarray}
and
\begin{eqnarray}
-P_{22}^{(2)}(\th,t+\d,t)^\top=-P_{22}^{(2)}(\th,t+\d, \th)^\top+\int_t^\th\dot P_{22}^{(2)}(t+\d,\th,s)ds\qq\qq\qq\qq\qq\qq\qq\qq\q\  \notag\\
=-\Big[P_{21}^{(1)}(t+\d)+P_{22}^{(1)}(t+\d) {\bf1}_{(\th,\i)}(t) \Big]A_2\1n-\3n\int_\th^T\3n\Big(P_{21}^{(2)}(t+\d,\a,\th) +P_{22}^{(2)}(t+\d,\a,\th){\bf1}_{(\th+\d,\i)}(\a)\Big) A_2d\a\q \notag\\
+\Big\{P_{21}^{(1)}(t+\d)+P_{22}^{(1)}(t+\d){\bf1}_{(\th, \i)}(t)+\int_\th^T\Big(P_{21}^{(2)}(t+\d,\a,\th) +P_{22}^{(2)}(t+\d,\a,\th){\bf1}_{(\th+\d,\i)}(\a)\Big) d\a\Big\}\qq\notag\\
\times B_1(R_1+D_1^\top \cP_2(t)D_1)^{-1} D_1^\top\cP_2(t)C_2 +\int_t^\th\Big[P_{21}^{(1)}(t+\d)+P_{22}^{(1)}(t+\d) {\bf1}_{(s,\i)}(t)\qq\qq\qq\qq\q\ \notag\\
+\int_s^T\Big(P_{21}^{(2)}(t+\d,\a,s) +P_{22}^{(2)}(t+\d,\a,s){\bf1}_{(s+\d,\i)}(\a)\Big)d\a \Big]B_1(R_1+D_1^\top\cP_2(t)D_1)^{-1} B_1^\top\Big[P_{21}^{(1)}(\th)^\top\qq\notag\\
+P_{22}^{(1)}(\th)^\top {\bf1}_{(t,\th-\d)}(s)+\int_s^T\Big(P_{21}^{(2)}(\th,\a,s) ^\top+P_{22}^{(2)}(\th,\a,s)^\top{\bf1}_{(t,\a-\d)}(s) \Big)d\a\Big]B_1ds.\qq\qq\qq\qq\ \ \ \!
\notag
\end{eqnarray}
Then, for $t\in(T-\d,T]$, $\th\in(t,T]$, with some calculations we deduce
\begin{eqnarray}
-\frac{\partial \cP_3(t, \theta)}{\partial t}=A_1^{\top} \cP_3(t, \theta)-\big(B_1^\top \cP_2(t)+D_1^\top \cP_2(t) C_1\big)^\top \big(R_1+D_1^\top\cP_2(t)D_1\big)^{-1} B_1^\top \cP_3(t, \theta), \label{P3-1}
\end{eqnarray}
and for $t\in[0,T-\d]$, $\th\in(t,t+\d]$, we get
\begin{eqnarray}
-\frac{\partial \cP_3(t, \theta)}{\partial t}=    A_1^\top \cP_3(t, \theta)-\big(B_1^\top \cP_2(t)+D_1^\top \cP_2(t) C_1\big)^\top \big(R_1+D_1^\top\cP_2(t)D_1\big)^{-1} B_1^\top \cP_3(t, \theta)\notag\\
\ns\ds \qq+\cP_3(\th, t+\d)^\top A_2-(B_1^\top \cP_3(\th, t+\d) )^\top  \big(R_1+D_1^\top\cP_2(t)D_1\big)^{-1} D_1 \cP_2(t) C_2 \q\notag\\
\ns\ds \qq-\int_t^\theta \big(B_1^\top \cP_3(s, t+\d)\big)^\top \big(R_1+D_1^\top\cP_2(t)D_1\big)^{-1} B_1^\top \cP_3(s, \th) d s. \label{P3-2}\qq\qq\q\ \
\end{eqnarray}
For $\cP_3(t,t)$, by \rf{coupled Riccati} we have
\begin{eqnarray}
\cP_3(t,t)=C_1^\top\cP_2(t)C_2-C_1^\top\cP_2(t)D_1(R_1 +D_1^\top\cP_2(t)D_1)^{-1} D_1^\top\cP_2(t)C_2\notag\\
+\Big[A_2\cP_2(t)^\top-C_2^\top\cP_2(t)^\top D_1(R_1+D_1^\top\cP_2(t)D_1)^{-1} B_1^\top\cP_2(t)\Big]^\top. \label{P3-3}
\end{eqnarray}
Finally, \rf{P3-1}--\rf{P3-3} imply that  the parts for $\cP_3(\cd,\cd)$ in \rf{Riccati-1-state delay}--\rf{Riccati-2-state delay} hold,  similarly we can also verify the parts for $\cP_2(\cd)$. In this case, \rf{Ki}--\rf{transformed optimal control} show that \rf{u*-state delay} is the optimal closed-loop outcome control.
\epf
%


\setcounter{equation}{0}   
\renewcommand\theequation{A.\arabic{equation}}  

\end{document}